\documentclass[12pt,fullpage]{article}

\usepackage{amsfonts,xr,graphicx,amsmath,amssymb,colortbl,epsfig,subfigure}
\def\R{\mathbb{R}}
\def\N{\mathbb{N}}

\def\epsilon{\varepsilon}

\def\beq{\begin{equation}}
\def\eeq{\end{equation}}

\def\div{\rm{div}}

\newcommand{\SE}{\setcounter{equation}{0} \section}
\newcommand{\baa}{\begin{array}}
\newcommand{\eaa}{\end{array}}
\newcommand{\ba}{\begin{eqnarray}}
\newcommand{\ea}{\end{eqnarray}}
\newtheorem{theo}{\bf Theorem}[section]
\newtheorem{lem}[theo]{\bf Lemma}
\newtheorem{pro}[theo]{\bf Proposition}

\newtheorem{conj}[theo]{\bf Conjecture}

\newtheorem{rem}[theo]{\bf Remark}

\begin{document}

\title{\bf{Circular flows for the Euler equations in two-dimensional annular domains, and related free boundary problems}}
\author{Fran{\c{c}}ois Hamel and Nikolai Nadirashvili$\,$\thanks{This work has been carried out in the framework of Archim\`ede Labex of Aix-Marseille University. The project leading to this publication has received funding from Excellence Initiative of Aix-Marseille University~-~A*MIDEX, a French ``Investissements d'Avenir'' programme, from the European Research Council under the European Union's Seventh Framework Programme (FP/2007-2013) ERC Grant Agreement n.~321186~- ReaDi~- Reaction-Diffusion Equations, Propagation and Modelling and from the ANR NONLOCAL (ANR-14-CE25-0013) and RESISTE (ANR-18-CE45-0019) projects. The authors are grateful to Jean-Michel Coron for suggesting this research topic and to Boyan Sirakov for pointing out his reference~\cite{si}.}\\
\\
\footnotesize{Aix Marseille Univ, CNRS, Centrale Marseille, I2M, Marseille, France}}
\date{}
\maketitle

\begin{abstract}
In this paper, we firstly consider steady Euler flows in two-dimensional bounded annuli, as well as in complements of disks, in punctured disks and in the punctured plane. We prove that, if the flow does not have any stagnation point and satisfies rigid wall boundary conditions together with further conditions at infinity in the case of unbounded domains and at the center in the case of punctured domains, then the flow is circular, namely the streamlines are concentric circles. In other words, the flow then inherits the radial symmetry of the domain. We secondly show two classification results for the steady Euler equations in simply or doubly connected bounded domains with free boundaries. Here, the flows are further assumed to have constant norm on each connected component of the boundary, and the domains are then proved to be disks or annuli. On the one hand, the proofs use ODE and PDE arguments to establish some geometric properties of the trajectories of the flow and the orthogonal trajectories of the gradient of the stream function. On the other hand, we also show some comparison results of independent interest for a derived semilinear elliptic equation satisfied by the stream function. These last results, which are based on the method of moving planes, adapted here to some almost circular domains located between some streamlines of the flow, lead with a limiting argument to the radial symmetry of the stream function and the streamlines of the flow.
\end{abstract}

\tableofcontents


\SE{Introduction and main results}\label{intro}

Throughout this paper, $|\ |$ denotes the Euclidean norm in $\R^2$ and, for $0\le a<b\le\infty$, $\Omega_{a,b}$ denotes the two-dimensional domain defined by
$$\Omega_{a,b}=\big\{x\in\R^2:\ a<|x|<b\big\}.$$
When $a<b$ are two positive real numbers, $\Omega_{a,b}$ is a bounded smooth annulus. When $0<a<b=\infty$, $\Omega_{a,\infty}$ is an exterior domain which is the complement of a closed disk. When $0=a<b<\infty$, $\Omega_{0,b}$ is a punctured disk. When $0=a<b=\infty$, $\Omega_{0,\infty}$ is the punctured plane~$\R^2\!\setminus\!\{0\}$, where we denote $0=(0,0)$ with a slight abuse of notation.\par We also denote
$$e_r(x)=\frac{x}{|x|}\ \hbox{ and }\ e_\theta(x)=e_r(x)^\perp=\Big(\!-\frac{x_2}{|x|},\frac{x_1}{|x|}\Big)$$
for $x=(x_1,x_2)\in\R^2\!\setminus\!\{0\}$. Moreover, for $x\in\R^2$ and $r>0$,
$$B(x,r)=\{y\in\R^2: |y-x|<r\}$$
denotes the open Euclidean disk with center $x$ and radius $r$. We also write $B_r=B(0,r)$ and
$$C_r=\partial B_r=\{x\in\R^2: |x|=r\}.$$


\subsection{Liouville-type radial symmetry results for steady Euler flows in~$\Omega_{a,b}$}

In $\Omega_{a,b}$, we consider steady flows
$$v=(v_1,v_2)$$
of an inviscid fluid, solving the system of the Euler equations: 
\begin{equation}\label{1}
\left\{\begin{array}{ll}
v\cdot\nabla\,v +\nabla\,p =0 & \mbox{in $\Omega_{a,b}$},\vspace{3pt}\\
{\div}\ v=0 & \mbox{in $\Omega_{a,b}$},
\end{array} \right.
\end{equation}
where the solutions $v$ and $p$ are always understood in the classical sense, that is, they are (at least) of class $C^1$ in $\Omega_{a,b}$ and therefore satisfy~\eqref{1} everywhere in $\Omega_{a,b}$. We always assume rigid wall boundary conditions, that is, $v$ is (at least) continuous up to the regular parts of $\partial\Omega_{a,b}$ and tangential there:
\beq\label{2}\left\{\baa{l}
v\cdot e_r=0\hbox{ on }C_a\hbox{ if }a>0,\vspace{3pt}\\
v\cdot e_r=0\hbox{ on }C_b\hbox{ if }b<\infty.\eaa\right.
\eeq\par
The Euler equations are an old but still very active research field. The search for qualitative properties of steady solutions is an important aspect of the study of the Euler flows, and the first main motivation of our paper is to understand the effect of the geometry of the underlying domain~$\Omega_{a,b}$ on the properties of steady flows, and more precisely to get some conditions on the flow which guarantee its radially symmetry (see the precise definition below). In other words, how does the solution inherit the geometric radial symmetry properties of the domain~? Our primary goal is thus to establish some Liouville-type results for the Euler equations. But the paper is also concerned with related Serrin-type free boundary problems, for which the underlying domain, simply or doubly connected, is free, but is eventually proved to be circular due to additional boundary conditions. Lastly, the paper contains some new comparison results of independent interest on semilinear elliptic equations in doubly connected domains, which are used to show the rigidity results for the Euler equations in given domains and for the related free boundary problems.\par
A flow $v$ in $\Omega_{a,b}$ is called a circular flow if $v(x)$ is parallel to the vector $e_\theta(x)$ at every point~$x\in\Omega_{a,b}$, that is, $v\cdot e_r=0$ in $\Omega_{a,b}$. The main goal of the results of this subsection is to show that, under some conditions, the flow is circular. We obtain such results in the four cases $0<a<b<\infty$, $0<a<b=\infty$, $0=a<b<\infty$, and $0=a<b=\infty$.

\subsubsection*{The case of bounded smooth annuli $\Omega_{a,b}$ with $0<a<b<\infty$}

The first result is concerned with flows having no stagnation point in the closed annulus $\overline{\Omega_{a,b}}$. Throughout the paper, the stagnation points of a flow~$v$ are the points~$x$ for which $|v(x)|=0$.

\begin{theo}\label{th1}
Assume $0<a<b<\infty$. Let $v$ be a $C^2(\overline{\Omega_{a,b}})$ flow solving~\eqref{1}-\eqref{2} and such that~$|v|>0$ in $\overline{\Omega_{a,b}}$. Then $v$ is a circular flow, and there is a $C^2([a,b])$ function $V$ with constant strict sign such that
$$v(x)=V(|x|)\,e_\theta(x)\ \hbox{ for all }x\in\overline{\Omega_{a,b}}.$$
\end{theo}

It actually turns out that the assumption $|v|>0$ in $\overline{\Omega_{a,b}}$ can be slightly relaxed. Namely, if~$|v|>0$ in the open annulus $\Omega_{a,b}$ and if the set of stagnation points is assumed to be properly included in one of the connected components of $\partial\Omega_{a,b}$, then the same conclusion holds, and then in fact $|v|>0$ in $\overline{\Omega_{a,b}}$. This is the purpose of the following result.

\begin{theo}\label{th1bis}
Assume $0<a<b<\infty$. Let $v$ be a $C^2(\overline{\Omega_{a,b}})$ flow solving~\eqref{1}-\eqref{2} and such that
\beq\label{hypstagnation}
\big\{x\in\overline{\Omega_{a,b}}: |v(x)|=0\big\}\subsetneq C_a\ \hbox{ or }\ \big\{x\in\overline{\Omega_{a,b}}: |v(x)|=0\big\}\subsetneq C_b.\footnote{Throughout the paper, by $E\subsetneq F$, we mean that $E\subset F$ and $E\neq F$.}
\eeq
Then $|v|>0$ in $\overline{\Omega_{a,b}}$ and the conclusion of Theorem~$\ref{th1}$ holds.
\end{theo}

Theorem~\ref{th1bis} is clearly stronger than Theorem~\ref{th1}, but we preferred to state Theorem~\ref{th1} separately since the assumption is simpler to read.\par
Several further comments are in order. First of all, despite the fact that $\Omega_{a,b}$ is not simply connected, the flow $v$ has a stream function $u:\overline{\Omega_{a,b}}\to\R$ of class $C^3(\overline{\Omega_{a,b}})$ defined by
\beq\label{defu}
\nabla^{\perp}u=v,\ \hbox{ that is},\ \ \frac{\partial u}{\partial x_1}=v_2\ \hbox{ and }\ \frac{\partial u}{\partial x_2}=-v_1
\eeq
in $\overline{\Omega_{a,b}}$, since $v$ is divergence free and tangential on $C_a$. Notice that the stream function~$u$ is uniquely defined in $\overline{\Omega_{a,b}}$ up to an additive constant. Theorems~\ref{th1} and~\ref{th1bis} can then be viewed as Liouville-type symmetry results since their conclusion means that the stream function $u$ is radially symmetric (and strictly monotone with respect to $|x|$ in $\overline{\Omega_{a,b}}$). Furthermore, if for $x$ in $\overline{\Omega_{a,b}}$ one calls $\xi_x$ the solution of
\beq\label{defxix}\left\{\baa{l}
\dot\xi_x(t)=v(\xi_x(t)),\vspace{3pt}\\
\xi_x(0)=x,\eaa\right.
\eeq
the conclusion of Theorems~\ref{th1} and~\ref{th1bis} then implies that each function $\xi_x$ is defined in $\R$ and periodic, and that the streamlines $\Xi_x=\xi_x(\R)$ of the flow are concentric circles.\par
Theorems~\ref{th1} and~\ref{th1bis} also mean equivalently that any $C^2(\overline{\Omega_{a,b}})$ non-circular flow for~\eqref{1}-\eqref{2} must either have a stagnation point in the open annulus~$\Omega_{a,b}$, or must have stagnation points in both circles $C_a$ and $C_b$, or in the whole circle~$C_a$, or in the whole circle $C_b$.\par
Without the assumption $|v|>0$ in $\overline{\Omega_{a,b}}$ or the weaker one~\eqref{hypstagnation}, the conclusion of Theorems~\ref{th1} and~\ref{th1bis} obviously does not hold in general, in the sense that there are non-circular flows which do not fulfill~\eqref{hypstagnation}. To construct such flows explicitly, we first point out that, for any continuous function $f:\R\to\R$ and any non-radial $C^2(\overline{\Omega_{a,b}})$ solution $u$ of
\beq\label{equ}
\Delta u+f(u)=0
\eeq
in $\Omega_{a,b}$ which is constant on $C_a$ and on $C_b$ and which has a critical point in $\overline{\Omega_{a,b}}$, the $C^1(\overline{\Omega_{a,b}})$ field
$$v=\nabla^\perp u$$
is a non-circular solution of~\eqref{1}-\eqref{2} with a stagnation point in $\overline{\Omega_{a,b}}$: notice indeed that $v=\nabla^\perp u$ satisfies the boundary condition $v\cdot e_r=-\nabla u\cdot e_\theta=0$ on $\partial\Omega_{a,b}$ since $u$ is constant on~$C_a$ and on $C_b$, and $v$ solves~\eqref{1} with pressure
$$p=-\frac{|v|^2}{2}-F(u)=-\frac{|\nabla u|^2}{2}-F(u),$$
where $F'=f$. As an example, let $\lambda\in\R$ and~$\varphi\in C^{\infty}([a,b])$ be the principal eigenvalue and the principal eigenfunction of the eigenvalue problem
$$-\varphi''(r)-r^{-1}\varphi'(r)+r^{-2}\varphi(r)=\lambda\,\varphi(r)\ \hbox{ in }[a,b]$$
with $\varphi>0$ in $(a,b)$ and Dirichlet boundary condition $\varphi(a)=\varphi(b)=0$ (the principal eigenvalue~$\lambda$ is unique and the principal eigenfunction $\varphi$ is unique up to multiplication by positive constants). The $C^{\infty}(\overline{\Omega_{a,b}})$ function~$u$ defined by $u(x)=\varphi(|x|)\,x_1/|x|$ (that is, $u(x)=\varphi(r)\,\cos(\theta)$ in the usual polar coordinates) satis\-fies
$$\Delta u+\lambda u=0\ \hbox{ in $\overline{\Omega_{a,b}}$}$$
and it has some critical points in $\overline{\Omega_{a,b}}$ (since $\min_{\overline{\Omega_{a,b}}}u<0<\max_{\overline{\Omega_{a,b}}}u$ and~$u=0$ on $\partial\Omega_{a,b}$). Actually, it can easily be seen that $\varphi$ has only one critical point in~$[a,b]$ and that $u$ has exactly 6 critical points in $\overline{\Omega_{a,b}}$ (2 in~$\Omega_{a,b}$, 2 on~$C_a$, and 2 on~$C_b$). Then the $C^\infty(\overline{\Omega_{a,b}})$ flow $v=\nabla^\perp u$ is a non-circular flow solving~\eqref{1}-\eqref{2} and having 2 stagnation points in $\Omega_{a,b}$ and 4 on $\partial\Omega_{a,b}$.\par
However, we do not know whether the hypothesis~\eqref{hypstagnation} could be more relaxed for the conclusion of Theorems~\ref{th1} and~\ref{th1bis} to still hold. For instance, would it be sufficient to assume that~$v$ has no stagnation point in $\Omega_{a,b}$? We refer to the comments after the proof of Theorems~\ref{th1} and~\ref{th1bis} in Section~\ref{sec31} below for further details on this question.\par
On the other hand, we point out that the sufficient conditions~$|v|>0$ in $\overline{\Omega_{a,b}}$ or the more general one~\eqref{hypstagnation} are obviously not equivalent to being a circular flow, in the sense that there are circular flows for~\eqref{1}-\eqref{2} which do not fulfill~\eqref{hypstagnation} (besides the trivial flow $v=(0,0)$!). Actually, any~$C^1(\overline{\Omega_{a,b}})$ circular flow $v(x)=V(|x|)\,e_\theta(x)$ solving~\eqref{1}-\eqref{2} and for which $V\in C^1([a,b])$ does not have a constant strict sign, has a set of stagnation points containing at least a circle. For instance, let $\mu\in\R$ and $\phi\in C^{\infty}([a,b])$ be the principal eigenvalue and the principal eigenfunction of the eigenvalue problem
$$-\phi''(r)-r^{-1}\phi'(r)=\mu\,\phi(r)\ \hbox{ in }[a,b],$$
with $\phi>0$ in $(a,b)$ and Dirichlet boundary condition $\phi(a)=\phi(b)=0$, and let $u=\phi(|\cdot|)$. Then $v=\nabla^\perp u=\phi'(|\cdot|)\,e_\theta$ is a $C^\infty(\overline{\Omega_{a,b}})$ non-trivial circular flow solving~\eqref{1}-\eqref{2} with pressure $p(x)=-\phi'(|x|)^2/2-\mu\phi(|x|)^2/2$ and with a circle of stagnation points in $\overline{\Omega_{a,b}}$: more precisely, if $r^*\in(a,b)$ denotes a real number such that $\phi(r^*)=\max_{[a,b]}\phi$ (it is easy to see that $r^*$ is the only critical point of $\phi$ in~$[a,b]$), then the set of stagnation points of the flow $v$ is equal to the whole circle~$C_{r^*}$.\par
Lastly, the assumption on the $C^2(\overline{\Omega_{a,b}})$ smoothness of $v$ is a technical assumption which is used in the proof. It warrants the $C^1$ smoothness of the vorticity
$$\omega=\frac{\partial v_2}{\partial x_1}-\frac{\partial v_1}{\partial x_2},$$
satisfying $v\cdot\nabla\omega=0$ in~$\Omega_{a,b}$, and the $C^1$ smoothness of the vorticity function $f$ arising in the semilinear elliptic equation of the type~\eqref{equ} satisfied by the stream function $u$. We refer to the proofs of the preliminary results in Section~\ref{sec2} and especially Lemma~\ref{lem24} below for further details.

\subsubsection*{The case of exterior domains $\Omega_{a,\infty}$ with $0<a<\infty$}

\begin{theo}\label{th2}
Assume $0<a<\infty$ and $b=\infty$. Let $v$ be a $C^2(\overline{\Omega_{a,\infty}})$ flow solving~\eqref{1}-\eqref{2} and such that
\beq\label{hypth2}
\big\{x\in\overline{\Omega_{a,\infty}}:|v(x)|=0\big\}\subsetneq C_a\ \hbox{ and }\ \liminf_{|x|\to+\infty}|v(x)|>0.
\eeq
Assume moreover that
\beq\label{hypth22}
v(x)\cdot e_r(x)=o\Big(\frac{1}{|x|}\Big)\ \hbox{ as }|x|\to+\infty.
\eeq
Then $|v|>0$ in $\overline{\Omega_{a,\infty}}$ and $v$ is a circular flow, namely there is a $C^2([a,+\infty))$ function $V$ with constant strict sign such that $v(x)=V(|x|)\,e_\theta(x)$ for all $x\in\overline{\Omega_{a,\infty}}$.
\end{theo}

As for Theorems~\ref{th1} and~\ref{th1bis}, the conclusion of Theorem~\ref{th2} says that the stream function~$u$ is radially symmetric and strictly monotone with respect to $|x|$ in $\overline{\Omega_{a,\infty}}$, and that the streamlines of the flow $v$ are concentric circles.\par
In the proof of Theorem~\ref{th2} given in Section~\ref{sec32}, the $o(1/|x|)$ behavior in~\eqref{hypth22} seems merely optimal. Let us show in this paragraph that without the condition~\eqref{hypth22} the conclusion of Theorem~\ref{th2} does not hold in general. Namely, consider the $C^{\infty}(\overline{\Omega_{a,\infty}})$ function $u$ defined by $u(x)=2(|x|^2/a^2-1)+(|x|/a-a/|x|)x_1/|x|$, that is,
$$u=2\,\Big(\frac{r^2}{a^2}-1\Big)+\Big(\frac{r}{a}-\frac{a}{r}\Big)\,\cos\theta$$
in the usual polar coordinates. The function $u$ satisfies $\Delta u-8/a^2=0$ in $\overline{\Omega_{a,\infty}}$ with Dirichlet boundary condition $u=0$ on $C_a$, and the $C^\infty(\overline{\Omega_{a,\infty}})$ field $v=\nabla^\perp u$ satisfies~\eqref{1}-\eqref{2} with pressure $p=-|v|^2/2+8u/a^2$. In the usual polar coordinates, the field $v$ is given by
\beq\label{defv}
v=\Big[\frac{4r}{a^2}+\Big(\frac{1}{a}+\frac{a}{r^2}\Big)\cos\theta\Big]\,e_\theta+\Big[\Big(\frac{1}{a}-\frac{a}{r^2}\Big)\,\sin\theta\Big]\,e_r.
\eeq
It satisfies condition~\eqref{hypth2} (and even $\inf_{\Omega_{a,\infty}}|v|\ge2/a>0$). But
$$v(x)\cdot e_r(x)=\Big(\frac{1}{a}-\frac{a}{|x|^2}\Big)\,\frac{x_2}{|x|}\neq o\Big(\frac{1}{|x|}\Big)\ \hbox{as }|x|\to+\infty,$$
and $v$ is not a circular flow. However, since $u(x)\to+\infty$ as $|x|\to+\infty$ and $u=0$ on $C_a$ and since $u$ has no critical point, it is easily seen that all solutions $\xi_x$ of~\eqref{defxix} are defined in~$\R$ and periodic and that all streamlines $\Xi_x=\xi_x(\R)$ (which are level sets of $u$) surround the origin.\footnote{Throughout the paper, we say that a Jordan curve $\mathcal{C}$ surrrounds the origin if the bounded connected component of $\R^2\!\setminus\!\mathcal{C}$ contains the origin.} Nevertheless, the streamlines do not converge to any family of circles at infinity since a calculation yields $\max_{y\in\Xi_x}|y|-\min_{y\in\Xi_x}|y|=\max_\R|\xi_x(\cdot)|-\min_\R|\xi_x(\cdot)|\to a/2>0$ as $|x|\to+\infty$. In this counterexample, one actually has $0<\limsup_{|x|\to+\infty}|v(x)\cdot e_r(x)|<+\infty$. Thus, there may be another critical behavior than $o(1/|x|)$ in~\eqref{hypth22} for which the conclusion would still hold, although a different proof would be necessary. The question of the characterization of a critical behavior is left open.\par
We point out that, in Theorem~\ref{th2}, the flow $v$ is not assumed to be bounded. Actually, there are unbounded circular flows satisfying all assumptions of Theorem~\ref{th2}: consider for instance the~$C^{\infty}(\overline{\Omega_{a,\infty}})$ unbounded circular flow $v$ defined by
$$v(x)=|x|\,e_\theta(x),$$
solving~\eqref{1}-\eqref{2} with stream function $u(x)=|x|^2/2$ and pressure~$p(x)=|x|^2/2$, and satisfying $\inf_{\Omega_{a,\infty}}|v|=a>0$.\par
Notice lastly that the condition~\eqref{hypth2} is fulfilled in particular when $\inf_{\Omega_{a,\infty}}|v|>0$. Furthermore, as soon as $|v|>0$ on $C_a$ (that holds if $\inf_{\Omega_{a,\infty}}|v|>0$), the boundary condition~\eqref{2} and the continuity of $v$ imply in particular that $v\cdot e_\theta$ has a constant strict sign on $C_a$. Under the condition~$\inf_{\Omega_{a,\infty}}|v|>0$, the following result then provides some estimates on the infimum or the supremum of the vorticity $\frac{\partial v_2}{\partial x_1}-\frac{\partial v_1}{\partial x_2}$ in $\Omega_{a,\infty}$, in terms of the sign of $v\cdot e_\theta$ on~$C_a$.

\begin{theo}\label{th2bis}
Assume $0<a<\infty$ and $b=\infty$. Let $v$ be a $C^2(\overline{\Omega_{a,\infty}})$ flow solving~\eqref{1}-\eqref{2} and such that~$\inf_{\Omega_{a,\infty}}|v|>0$. If $v\cdot e_\theta>0$ on $C_a$ $($respectively if $v\cdot e_\theta<0$ on $C_a$$)$, then
$$\sup_{\Omega_{a,\infty}}\!\Big(\frac{\partial v_2}{\partial x_1}-\frac{\partial v_1}{\partial x_2}\Big)>0\ \ (\hbox{respectively }\inf_{\Omega_{a,\infty}}\!\Big(\frac{\partial v_2}{\partial x_1}-\frac{\partial v_1}{\partial x_2}\Big)<0).$$
\end{theo}

The flow $v$ given by~\eqref{defv} is an example of a flow satisfying the assumptions of Theorem~\ref{th2bis}, with $v\cdot e_\theta>0$ on $C_a$, and for which the vorticity (namely $\Delta u$) is actually equal to the positive constant  $8/a^2$ everywhere in $\overline{\Omega_{a,\infty}}$.\par
Theorem~\ref{th2bis} can also be viewed as a Liouville-type result. Namely, we show in its proof that, if~$\inf_{\Omega_{a,\infty}}|v|>0$, if $v\cdot e_\theta>0$ on $C_a$, and if the vorticity is nonpositive everywhere in~$\Omega_{a,\infty}$, then $v$ is a circular flow of the type $v=V(|\cdot|)\,e_\theta$ with $V:[a,+\infty)\to[\eta,+\infty)$ for some $\eta>0$. Therefore, the vorticity $\frac{\partial v_2}{\partial x_1}(x)-\frac{\partial v_1}{\partial x_2}(x)=V'(|x|)+V(|x|)/|x|$ can not be nonpositive everywhere (since otherwise the function $r\mapsto r\,V(r)\ (\ge\eta r)$ would be nonincreasing in~$[a,+\infty)$, leading to a contradiction).\footnote{The same arguments do not lead to any contradiction in the case of bounded annuli $\Omega_{a,b}$ and $\Omega_{0,b}$ with~$b<\infty$, see Remark~\ref{remomega0b}.}\par
Notice that Theorem~\ref{th2bis} does not hold good if the assumption $\inf_{\Omega_{a,\infty}}|v|>0$ is dropped. There are actually some circular flows $v$ satisfying~\eqref{1}-\eqref{2} such that $|v|>0$ in $\overline{\Omega_{a,\infty}}$ and $v\cdot e_\theta>0$ on~$C_a$, but $\inf_{\Omega_{a,\infty}}|v|=0$ and for which the vorticity is negative everywhere. Consider for instance the~$C^{\infty}(\overline{\Omega_{a,\infty}})$ circular flow
$$v(x)=\frac{1}{|x|^2}\,e_\theta(x),$$
solving~\eqref{1}-\eqref{2} with stream function $u(x)=-1/|x|$ and pressure $p(x)=-1/(4|x|^2)$: one has~$|v|>0$ in $\overline{\Omega_{a,\infty}}$ and $v\cdot e_\theta>0$ on $C_a$, but $\inf_{\Omega_{a,\infty}}|v|=0$ and $\frac{\partial v_2}{\partial x_1}(x)-\frac{\partial v_1}{\partial x_2}(x)=-1/|x|^3<0$ in~$\overline{\Omega_{a,\infty}}$.

\subsubsection*{The case of punctured disks~$\Omega_{0,b}$ with $0<b<\infty$}

\begin{theo}\label{th3}
Assume $a=0$ and $0<b<\infty$. Let $v$ be a $C^2(\overline{\Omega_{0,b}}\!\setminus\!\{0\})$ flow solving~\eqref{1}-\eqref{2} and such that
\beq\label{hypth31}
\big\{x\in\overline{\Omega_{0,b}}\!\setminus\!\{0\}:|v(x)|=0\big\}\subsetneq C_b
\eeq
and
\beq\label{hypth32}
\int_{C_\epsilon}|v\cdot e_r|\to0\hbox{ as }\epsilon\displaystyle\mathop{\to}^>0.
\eeq
Then $|v|>0$ in $\overline{\Omega_{0,b}}\!\setminus\!\{0\}$ and $v$ is a circular flow, namely there is a~$C^2((0,b])$ function $V$ with constant strict sign such that $v(x)=V(|x|)\,e_\theta(x)$ for all $x\in\overline{\Omega_{0,b}}\!\setminus\!\{0\}$.
\end{theo}

Notice that the condition~\eqref{hypth32} is fulfilled in particular if $v(x)\cdot e_r(x)=o(1/|x|)$ as $|x|\displaystyle\mathop{\to}^>0$. Let us show in this paragraph that without~\eqref{hypth32} the conclusion of Theorem~\ref{th3} does not hold in general. To do so, let us give a counter-example similar to~\eqref{defv} above (which was there defined in~$\overline{\Omega_{a,\infty}}$). More precisely, consider the $C^{\infty}(\overline{\Omega_{0,b}}\!\setminus\!\{0\})$ function $u$ defined by $u(x)=(|x|/b-b/|x|)x_1/|x|$, that~is,
$$u=\Big(\frac{r}{b}-\frac{b}{r}\Big)\cos\theta$$
in the usual polar coordinates. The function $u$ satisfies $\Delta u=0$ in $\overline{\Omega_{0,b}}\!\setminus\!\{0\}$ with Dirichlet boundary condition $u=0$ on $C_b$, and the $C^\infty(\overline{\Omega_{0,b}}\!\setminus\!\{0\})$ field $v=\nabla^\perp u$ satisfies~\eqref{1}-\eqref{2} with pressure~$p=-|v|^2/2$ (and vorticity equal to $0$). In the usual polar coordinates, the field~$v$ is given by
\beq\label{defv2}
v=\Big[\Big(\frac{1}{b}+\frac{b}{r^2}\Big)\cos\theta\Big]\,e_\theta+\Big[\Big(\frac{1}{b}-\frac{b}{r^2}\Big)\,\sin\theta\Big]\,e_r.
\eeq
It has only two stagnation points in $\overline{\Omega_{0,b}}\!\setminus\!\{0\}$ and they both lie on $C_b$. Hence,~\eqref{hypth31} is fulfilled. But $\int_{C_\epsilon}|v\cdot e_r|=4(b/\epsilon-\epsilon/b)\not\to0$ as $\epsilon\displaystyle\mathop{\to}^>0$, and $v$ is not a circular flow. In this counterexample, one actually has $\int_{C_\epsilon}|v\cdot e_r|\sim 4b/\epsilon$ as $\epsilon\displaystyle\mathop{\to}^>0$. Thus, there may be another critical behavior than~$o(1)$ in the condition~\eqref{hypth32} for which the conclusion would still hold, although a different proof would be necessary. The question of the characterization of a critical behavior is still open.\par
Lastly, in Theorem~\ref{th3}, the flow $v$ is not assumed to be bounded. Actually, there are unbounded circular flows satisfying all assumptions of Theorem~\ref{th3}: consider for instance the~$C^{\infty}(\overline{\Omega_{0,b}}\!\setminus\!\{0\})$ unbounded circular flow $v$ defined by
\beq\label{vln}
v(x)=\frac{1}{|x|}\,e_\theta(x)
\eeq
solving~\eqref{1}-\eqref{2} with stream function $u(x)=\ln|x|$ and pressure~$p(x)=-1/(2|x|^2)$, and satis\-fying~$|v|>0$ in $\overline{\Omega_{0,b}}\!\setminus\!\{0\}$ and then~\eqref{hypth31}-\eqref{hypth32}.

\begin{rem}\label{remomega0b}{\rm A result similar to Theorem~$\ref{th2bis}$ does not hold in the punctured disk $\Omega_{0,b}$. For instance, the $C^{\infty}(\overline{\Omega_{0,b}}\!\setminus\!\{0\})$ flow~\eqref{vln} satisfies~\eqref{1}-\eqref{2}, $v\cdot e_\theta>0$ on $C_b$, $\inf_{\Omega_{0,b}}|v|>0$, but $\frac{\partial v_2}{\partial x_1}-\frac{\partial v_1}{\partial x_2}\equiv 0$ in $\Omega_{0,b}$. The same observation holds good in a smooth annulus $\Omega_{a,b}$ with $0<a<b<\infty$.}
\end{rem}

\begin{rem}{\rm{Let us comment here the similar conditions~\eqref{hypth22} and~\eqref{hypth32}. Condition~\eqref{hypth22} implies that $\int_{C_R}|v\cdot e_r|\to0$ as $R\to+\infty$, which would be the dual of~\eqref{hypth32}. But the stronger pointwise asymptotic behavior~\eqref{hypth22} is truly used in the proof of Theorem~\ref{th2}, and in particular in the proof of Lemma~\ref{lem8} below. Replacing~\eqref{hypth22} by the weaker integral condition $\lim_{R\to+\infty}\int_{C_R}|v\cdot e_r|=0$ still yields some intermediary results (see the common preliminary results of Section~\ref{sec2} below), but it is an open question to decide whether the conclusion of Theorem~\ref{th2} would still hold with the integral condition instead of the pointwise one.}}
\end{rem}

\subsubsection*{The case of the punctured plane~$\Omega_{0,\infty}$}

The last geometric configuration considered in the paper is that the punctured plane
$$\Omega_{0,\infty}=\R^2\!\setminus\!\{0\}.$$

\begin{theo}\label{th4}
Let $v$ be a $C^2(\Omega_{0,\infty})$ flow solving~\eqref{1} and such that~$|v|>0$ in $\Omega_{0,\infty}$ and $\liminf_{|x|\to+\infty}|v(x)|>0$. Assume moreover that
\beq\label{hypth4}
v(x)\cdot e_r(x)=o\Big(\frac{1}{|x|}\Big)\hbox{ as }|x|\to+\infty\ \hbox{ and }\ \int_{C_\epsilon}|v\cdot e_r|\to0\hbox{ as }\epsilon\displaystyle\mathop{\to}^>0.
\eeq
Then~$v$ is a circular flow. Furthermore, there is a~$C^2((0,+\infty))$ function $V$ with constant strict sign such that $v(x)=V(|x|)\,e_\theta(x)$ for all $x\in\Omega_{0,\infty}$.
\end{theo}

The conclusion says that, under roughly speaking the absence of stagnation points in the punctured plane and at infinity, and under the same conditions as in Theorems~\ref{th2} and~\ref{th3} on the behavior of the radial component of $v$ at infinity and at the origin, all streamlines are closed and are nothing but concentric circles.

\begin{rem}{\rm Let us mention here other rigidity results for the stationary solutions of~\eqref{1} in various geometrical configurations. The analyticity of the streamlines under a condition of the type~$v_1>0$ in the unit disk was shown in~\cite{kn}. The local correspondence between the vorticities of the solutions of~\eqref{1} and the co-adjoint orbits of the vorticities for the non-stationary version of~\eqref{1} in more general annular domains was investigated in~\cite{cs}. In a previous paper~\cite{hn1} (see also~\cite{hn3}), we considered the case of a two-dimensional strip with bounded section and the case of bounded flows in a half-plane, assuming in both cases that the flows $v$ are tangential on the boundary and that $\inf\,|v|>0$: all streamlines are then proved to be lines which are parallel to the boundary of the domain (in other words the flow is a parallel flow). Compared to~\cite{hn1}, the results of the present paper are concerned with different geometrical situations, and the cases of punctured disks or exterior domains involve specific difficulties. We here also include and prove some new comparison results of independent interest for the solutions of semilinear elliptic equations in doubly connected domains, see Proposition~\ref{promoving} in Section~\ref{sec13} below, not to mention the Serrin-type free boundary problems considered in Section~\ref{sec12}. These types of problems, as well as some methods used here such as the method of moving planes or the Kelvin transform, were not used in~\cite{hn1,hn3,hn2}. Earlier results by Kalisch~\cite{k} were concerned with flows in two-dimensional strips under the additional assumption~$v\cdot e\neq0$, where $e$ is the main direction of the strip. In~\cite{hn2}, we considered the case of the whole plane~$\R^2$ and we showed that any $C^2(\R^2)$ bounded flow $v$ is still a parallel flow under the condition~$\inf_{\R^2}|v|>0$, with completely different tools based on the study of the growth of the argument of the flow at infinity.}
\end{rem}


\subsection{Serrin-type free boundary problems with overdetermined boun\-dary conditions}\label{sec12}

The last main results on the solutions of the Euler equations~\eqref{1} are two Serrin-type results in smooth simply or doubly connected bounded domains whose boundaries are free but on which the flow is assumed to satisfy an additional condition.

\begin{theo}\label{th5}
Let $\Omega$ be a $C^2$ non-empty simply connected bounded domain of $\R^2$. Let $v\in C^2(\overline{\Omega})$ satisfy the Euler equations~\eqref{1} and assume that $v\cdot n=0$ and $|v|$ is constant on~$\partial\Omega$, where $n$ denotes the outward unit normal on $\partial\Omega$. Assume moreover that $v$ has a unique stagnation point in $\overline{\Omega}$. Then, up to translation,
$$\Omega=B_R$$
for some $R>0$. Furthermore, the unique stagnation point of $v$ is the center of the disk and~$v$ is a circular flow, that is, there is a $C^2([0,R])$ function~$V:[0,R]\to\R$ such that $V\neq0$ in~$(0,R]$, $V(0)=0$, and $v(x)=V(|x|)\,e_\theta(x)$ for all $x\in\overline{B_R}\!\setminus\!\{0\}$. 
\end{theo}

In the proof, we will show that the $C^3(\overline{\Omega})$ stream function $u$ defined by~\eqref{defu} satisfies a semilinear elliptic equation $\Delta u+f(u)=0$ in $\overline{\Omega}$. Furthermore, up to normalization, the function~$u$ vanishes on $\partial\Omega$ and is positive in $\Omega$. Lastly, since $|v|$ is assumed to be constant along $\partial\Omega$, the normal derivative $\frac{\partial u}{\partial n}$ of $u$ along $\partial\Omega$ is constant. This problem is therefore an elliptic equation with overdetermined boundary conditions. Since the celebrated paper by Serrin~\cite{s}, it has been known that these overdetermined boundary conditions on $\partial\Omega$ determine the geometry of~$\Omega$, namely, $\Omega$ is then a ball and the function $u$ is radially symmetric (hence, here, $v$ would then be a circular flow). The proof is based on the method of moving planes developed in~\cite{a,bn,gnn,s} and on the maximum principle, and it relies on the Lipschitz continuity of the function~$f$. In our case, the function $f$ is given in terms of the function $u$ itself and it is continuous in $[0,\max_{\overline{\Omega}}u]$, as will be seen in the proof of Theorem~\ref{th5}. But it can be non-Lipschitz-continuous on the whole range $[0,\max_{\overline{\Omega}}u]$. More precisely, it can be non-Lipschitz-continuous in any left neighborhood of the maximal value $\max_{\overline{\Omega}}u$.\footnote{For instance, for any $R>0$, the $C^\infty(\overline{B_R})$ flow $v(x)=-4|x|^2x^\perp$ satisfies~\eqref{1} in $\overline{B_R}$ with $v\cdot e_r=0$ on~$C_R$, and with pressure $p(x)=(8/3)|x|^6$ (up to an additive constant). Furthermore, $|v|$ is constant on $C_R$ and the only stagnation point of $v$ in $\overline{B_R}$ is the center of the disk. The stream function $u(x)=R^4-|x|^4$ (up to an additive constant) satisfies the elliptic equation $\Delta u+f(u)=0$ in $\overline{B_R}$ with $f(s)=16\sqrt{R^4-s}$, and the function~$f$ is not Lipschitz continuous in any left neighborhood of $\max_{\overline{B_R}}u=R^4$.} One therefore has to adapt the proof to this case by removing small neighborhoods of size $\epsilon$ around the maximal point of $u$ (which is the unique stagnation point of $v$): one shows the symmetry of the domain in all directions up to $\epsilon$ and one concludes by passing to the limit as $\epsilon\displaystyle\mathop{\to}^>0$.

\begin{rem}{\rm{Other free boundary problems related to the Euler equations have been considered by G\'omez-Serrano, Park, Shi and Yao in~\cite{gpsy}. Among other things, the authors proved that, if a solution $v$ of the Euler equations~\eqref{1} in $\R^2$ has a vorticity which is the indicator function of a bounded set (a patch) and if $v$ is tangential on the boundary of this set, then $v$ is circular (up to translation) and the patch is a disk (see also~\cite[Chapter~4]{f} for an earlier result when the patch is assumed to be simply connected). It was also shown in~\cite{gpsy} that smooth solutions of~\eqref{1} in $\R^2$ with nonnegative compactly supported vorticity must be radially symmetric (up to translation). Other rigidity results of~\cite{gpsy,h} also deal with non-stationary uniformly-rotating solutions.}}
\end{rem}

In connection with Theorems~\ref{th3} and~\ref{th5}, we state the following conjecture.

\begin{conj}\label{conj1}
Let $D$ be an open non-empty disk and let $z\in D$. Let $v$ be a $C^2(\overline{D}\!\setminus\!\{z\})$ and bounded flow solving~\eqref{1} and $v\cdot n=0$ on $\partial D$, where $n$ denotes the outward unit normal on~$\partial D$. Assume that $|v|>0$ in $\overline{D}\!\setminus\!\{z\}$. Then $z$ is the center of the disk and the flow is circular with respect to~$z$.
\end{conj}

Up to translation, one can assume that $D=B_b$ for some $b\in(0,+\infty)$, hence $n=e_r$ on~$\partial D$. If the point~$z$ is a priori assumed to be the center of the disk, namely the origin, then Theorem~\ref{th3} implies that $v$ is a circular flow. Up to rotation, assume now that $z=(\alpha,0)$ for some $\alpha\in(0,b)$ and, without loss of generality, that the stream function $u$ is positive in $D\!\setminus\!\{z\}$ and vanishes on~$\partial D$. The goal would be to reach a contradiction. As far as Theorem~\ref{th5} is concerned, the method of proof described in the paragraph following the statement shows simultaneously the symmetry of the domain and the symmetry of the function~$u$ (which obeys an equation of the type $\Delta u+f(u)=0$), thanks to the overdetermined boundary conditions satisfied by~$u$. Here in Conjecture~\ref{conj1}, the same technics based on the method of moving method implies for instance on the one hand that the function~$u$ is even in~$x_2$ in~$\overline{\Omega_{0,b}}\!\setminus\!\{z\}$, and on the other hand that $u(x_1,x_2)<u(2\alpha-x_1,x_2)$ for all $(x_1,x_2)\in\overline{\Omega_{0,b}}$ such that $x_1>\alpha$. But, regarding the second property, the Hopf lemma might not apply to the function $(x_1,x_2)\mapsto u(x_1,x_2)-u(2\alpha-x_1,x_2)$ at the point $z=(\alpha,0)$ since the vorticity function~$f$ might not be Lipschitz continuous around the limiting value of~$u$ at~$z$ (see also the comments after the proof of Theorems~\ref{th1} and~\ref{th1bis} in Section~\ref{sec31} below, and notice that $u$ is not differentiable at $z$, unless one further assumes that $|v(x)|\to 0$ as $x\to z$). Therefore, the same arguments as the ones in the proof of Theorem~\ref{th5} do not lead to an obvious contradiction if $z$ is not the center of the disk. However, Conjecture~\ref{conj1} seems natural and will be the purpose of further investigation.\par
A related weaker conjecture (with stronger assumptions) can also be formulated: if $D$ is an open non-empty disk, if $z\in D$, if $v\in C^2(\overline{D})$ solves~\eqref{1}, if $v\cdot n=0$ on $\partial D$ and if $z$ is the only stagnation point of $v$ in $\overline{D}$, then $z$ is the center of the disk and $v$ is circular with respect to $z$. For the same reasons as in the previous paragraph (since the vorticity function $f$ might not be Lipschitz continuous around $u(z)$), the proof of that second conjecture is not clear either.\par
The last main result related to the Euler equations is concerned with the case of doubly connected bounded domains.

\begin{theo}\label{th5bis}
Let $\omega_1$ and $\omega_2$ be two $C^2$ non-empty simply connected bounded domains of $\R^2$ such that $\overline{\omega_1}\subset\omega_2$, and denote
$$\Omega=\omega_2\!\setminus\!\overline{\omega_1}.$$
Let $v\in C^2(\overline{\Omega})$ satisfy the Euler equations~\eqref{1}. Assume that $v\cdot n=0$ on $\partial\Omega=\partial\omega_1\cup\partial\omega_2$, where $n$ denotes the outward unit normal on $\partial\Omega$, and that $|v|$ is constant on $\partial\omega_1$ and on $\partial\omega_2$. Assume moreover that $|v|>0$ in $\overline{\Omega}$. Then $\omega_1$ and $\omega_2$ are two concentric disks and, up to translation,
$$\Omega=\Omega_{a,b}$$
for some $0<a<b<\infty$ and $v$ is a circular flow satisfying the conclusion of Theorem~$\ref{th1}$ in~$\overline{\Omega}=\overline{\Omega_{a,b}}$.
\end{theo}

In this case, by using the arguments of Section~\ref{sec2} below (which also lead to the proof of Theorems~\ref{th1} and~\ref{th1bis} in~$\Omega_{a,b}$ with $0<a<b<\infty$), it follows that the stream function $u$ of the flow $v$ satisfies a semilinear elliptic equation $\Delta u+f(u)=0$ in $\overline{\Omega}$, with $u=c_1$ on $\partial\omega_1$ and~$u=c_2$ on $\partial\omega_2$, for some real numbers $c_1\neq c_2$. Furthermore, $\min(c_1,c_2)<u<\max(c_1,c_2)$ in~$\Omega$ and the normal derivative~$\frac{\partial u}{\partial n}$ is constant along~$\partial\omega_1$ and along $\partial\omega_2$. Since 
$v$ has no stagnation point in $\overline{\Omega}$, the function~$f$ is then shown to be Lipschitz continuous in the whole interval~$[\min(c_1,c_2),\max(c_1,c_2)]$, and known results of Reichel~\cite{r1} and Sirakov~\cite{si} then imply that, up to translation, $\Omega=\Omega_{a,b}$ for some $0<a<b<\infty$, and $u$ is radially symmetric.\par
Further symmetry results have been obtained for nonlinear elliptic equations of the type $\Delta u+f(u)=0$ or more general ones in exterior domains with overdetermined boundary conditions (see e.g.~\cite{ab,r2,si}), or in the whole space (see e.g.~\cite{gnn2,li,sz}), in both cases with further assumptions on the solution $u$ at infinity and on the function $f$. Such conditions are in general not satisfied by the stream function $u$ and the vorticity function $f$ of a flow $v$ that would be defined in the complement of a simply connected bounded domain or in the whole or punctured plane. Lastly, we refer to~\cite{bcn,fv,rrs,rs} for further references on overdetermined boundary value elliptic problems in domains with more complex topology or in unbounded epigraphs.


\subsection{Directional comparison results for semilinear elliptic equations in some doubly connected domains}\label{sec13}

As briefly mentioned after Theorems~\ref{th5} and~\ref{th5bis}, the main strategy of these results, as well as the other ones in the fixed annular domains $\Omega_{a,b}$, is to show that the stream function $u$ satisfies a semilinear elliptic equation $\Delta u+f(u)=0$ in the considered domain and that this stream function is then radially symmetric (and the domain itself is circular if the boundary is free).  The radial symmetry of the stream function $u$ means that the flow $v$ is circular. The proof of the radial symmetry of~$u$ follows from its even symmetry and monotonicity with respect to each direction. The proof of the symmetry and monotonicity relies on some directional comparison results which are themselves based on the method of moving planes in doubly connected domains trapped between two level sets of $u$. We point out that these level sets are not known a priori to be circles and their precise shape is not known. This is why we have to show in this framework a key-proposition containing some new directional comparison results, which we think are of independent interest and which we state and will use for more general semilinear heterogeneous elliptic equations. Proposition~\ref{promoving} below will be used as a key-step in the proof of Theorems~\ref{th2}-\ref{th3},~\ref{th4} and~\ref{th5}. We also point out that the heterogeneity in the considered equations~\eqref{eqvarphi} below is not an artifice, since we will truly deal with heterogeneous equations obtained after a Kelvin transform of some original equations set in exterior domains.\par
To do so, let us first introduce a few notations. For $e\in\mathbb{S}^1=C_1$ and $\lambda\in\R$, we denote
\beq\label{Telambda}
T_{e,\lambda}=\big\{x\in\R^2:x\cdot e=\lambda\big\},\ \ H_{e,\lambda}=\big\{x\in\R^2:x\cdot e>\lambda\big\},
\eeq
and, for $x\in\R^2$,
\beq\label{Relambda}
R_{e,\lambda}(x)=x_{e,\lambda}=x-2(x\cdot e-\lambda)e.
\eeq
In other words, $R_{e,\lambda}$ is the orthogonal reflection with respect to the line $T_{e,\lambda}$.

\begin{pro}\label{promoving}
Let $\Xi$ and $\Xi'$ be two $C^1$ Jordan curves surrounding the origin, and let $\Omega$ and $\Omega'$ be the bounded connected components of $\R^2\!\setminus\!\Xi$ and $\R^2\!\setminus\!\Xi'$, respectively. Assume that $\overline{\Omega'}\subset\Omega$ and let
$$\omega=\Omega\!\setminus\!\overline{\Omega'}$$
be the non-empty and doubly connected domain located between $\Xi$ and $\Xi'$, with boundary
$$\partial\omega=\Xi\cup\Xi'.$$
Call $R'=\min_{x\in\Xi'}|x|>0$ and $R=\max_{x\in\Xi}|x|>R'$. Let $e\in\mathbb{S}^1$, let $\overline{\lambda}=\max_{x\in\Xi}x\cdot e>0$ and let $\epsilon\in[0,\overline{\lambda})$. Let $c_1<c_2\in\R$ and let $\varphi\in C^2(\overline{\omega})$ be a solution of
\beq\label{eqvarphi}\left\{\baa{rcll}
\Delta\varphi+F(|x|,\varphi) & \!\!=\!\! & 0 & \!\!\hbox{in }\overline{\omega},\vspace{3pt}\\
c_1\ <\ \varphi & \!\!<\!\! & c_2 & \!\!\hbox{in }\omega,\vspace{3pt}\\
\varphi\ =\ c_1\hbox{ on }\Xi,\ \ \varphi & \!\!=\!\! & c_2 & \!\!\hbox{on }\Xi',\eaa\right.
\eeq
with a continuous function $F:[R',R]\times[c_1,c_2]\to\R$ that is nonincreasing with respect to its first variable and uniformly Lipschitz continuous with respect to its second variable. Assume that
\beq\label{hypelambda}
R_{e,\lambda}(H_{e,\lambda}\cap\overline{\Omega})\subset\Omega\ \hbox{ for all }\lambda>\epsilon
\eeq
and that
\beq\label{hypelambda2}
R_{e,\lambda}(H_{e,\lambda}\cap\Xi')\subset\Omega'\ \hbox{ for all }\lambda>\epsilon,
\eeq
see Fig.~$1$. Then, for every $\lambda\in[\epsilon,\overline{\lambda})$, there holds
\beq\label{varphielambda}
\varphi(x)\le\varphi_{e,\lambda}(x)=\varphi(x_{e,\lambda})\ \hbox{ for all }x\in\overline{\omega_{e,\lambda}},
\eeq
with
$$\omega_{e,\lambda}=(H_{e,\lambda}\cap\omega)\setminus R_{e,\lambda}(\overline{\Omega'}).\footnote{Notice that $\omega_{e,\lambda}$ is open by definition, and it is non-empty for each $\lambda\in[0,\overline{\lambda})$: indeed, for such $\lambda$, the set~$T_{e,\lambda}\cap\Xi$ is not empty and, for any $x\in T_{e,\lambda}\cap\Xi$ and $r>0$, $\omega_{e,\lambda}\cap B(x,r)\neq\emptyset$. But $\omega_{e,\lambda}$ may not be connected, as in Fig.~1.}$$
\end{pro}

\begin{figure}
\centering\includegraphics[scale=0.7]{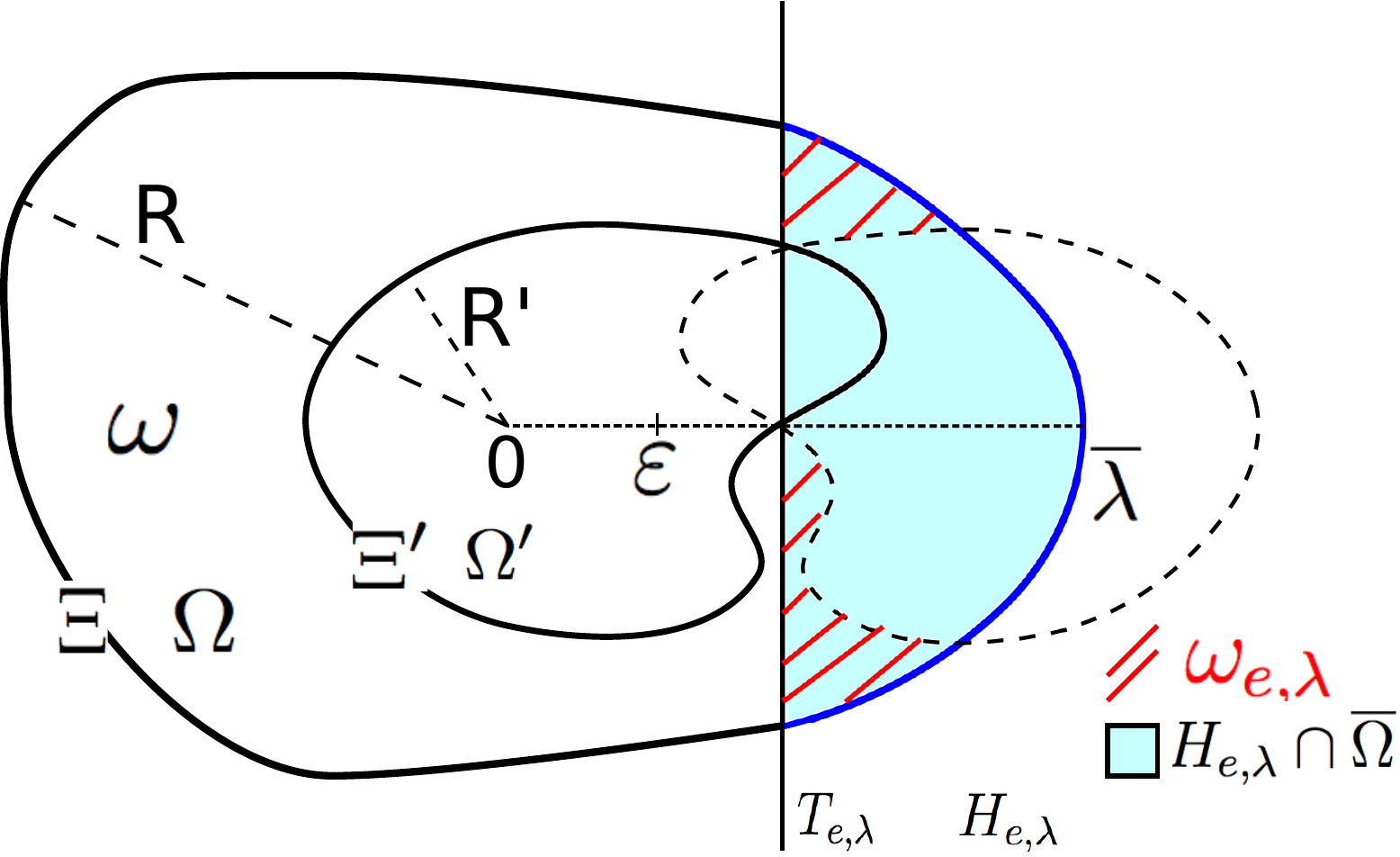}
\caption{The sets $\Omega$, $\Omega'$, $\omega=\Omega\!\setminus\!\overline{\Omega'}$, $H_{e,\lambda}\cap\overline{\Omega}$ (light blue background), $R_{e,\lambda}(\Omega')$ (with dashed boundary), and $\omega_{e,\lambda}$ (dashed red)}
\end{figure}

\subsubsection*{Organization of the paper}

In Section~\ref{sec2}, we show some common preliminaries for the proofs of Theorems~\ref{th1}-\ref{th3}, \ref{th4}, \ref{th5} and~\ref{th5bis} in given or free domains. Namely, we study the properties of the streamlines of the flow and we derive a semilinear elliptic equation for the stream function. Section~\ref{sec3} is devoted to the proof of Theorems~\ref{th1}-\ref{th3} and~\ref{th4}, when the circular domain $\Omega_{a,b}$ is fixed. The cases of the exterior domains $\Omega_{a,\infty}$ and the punctured disks $\Omega_{0,b}$ and punctured plane $\Omega_{0,\infty}$ involve some additional difficulties and require specific additional assumptions. The proof of the Serrin-type  Theorems~\ref{th5} and~\ref{th5bis} is carried out in Section~\ref{sec40}. Lastly, Section~\ref{sec50} is concerned with the proof of Proposition~\ref{promoving}. Proposition~\ref{promoving} and its subsequent limiting argument showing that the considered flow is circular are necessary for the proof of Theorems~\ref{th2}-\ref{th3} and~\ref{th4} in the case of exterior or punctured domains, as well as for the proof of Theorem~\ref{th5} about simply connected domains with free boundary.


\SE{Some common preliminaries}\label{sec2}

In this section, we state and prove some common properties which will be used in the proof of the main results related to the Euler equations~\eqref{1} in a fixed annular domain~$\Omega_{a,b}$ or in simply or doubly connected domains with free boundaries. To cover all possible cases, we consider throughout this section two~$C^1$ Jordan curves~$\mathcal{C}_1$ and~$\mathcal{C}_2$ and we assume that they both surround the origin, that is, the bounded connected components~$\omega_1$ and~$\omega_2$ of~$\R^2\!\setminus\!\mathcal{C}_1$ and~$\R^2\!\setminus\!\mathcal{C}_2$ contain the origin~$0$. We further assume that
$$\overline{\omega_1}\subset\omega_2.$$\par
Our aim here is to study the properties of the stream function $u$ and the streamlines of a divergence-free flow $v$ in the open connected set $\Omega$ in one of the following four possible cases:
$$\Omega=\omega_2\!\setminus\!\overline{\omega_1},\ \hbox{ or }\ \Omega=\R^2\!\setminus\!\overline{\omega_1},\ \hbox{ or }\ \Omega=\omega_2\!\setminus\!\{0\},\ \hbox{ or }\ \Omega=\R^2\!\setminus\!\{0\}.$$
Notice that $\Omega$ is {\it bounded} if $\Omega=\omega_2\!\setminus\!\overline{\omega_1}$ or $\omega_2\!\setminus\!\{0\}$, and {\it unbounded} in the other two cases. We say that $\Omega$ is {\it not punctured} if $\Omega=\omega_2\!\setminus\!\overline{\omega_1}$ or $\R^2\!\setminus\!\overline{\omega_1}$, and {\it punctured} in the other two cases. We point out that the annular domains $\Omega_{a,b}$ with $0\le a<b\le\infty$ fall within these four types of domains $\Omega$. We also set
\beq\label{defD}\left\{\baa{ll}
D=\overline{\Omega} & \hbox{if $\Omega$ is not punctured},\vspace{3pt}\\
D=\overline{\Omega}\!\setminus\!\{0\} & \hbox{if $\Omega$ is punctured}.\eaa\right.
\eeq\par
Throughout this section, $v=(v_1,v_2)$ is a $C^1(D)$ vector field. We point out that the $C^2(D)$ regularity of $v$ will be specified and used only in Lemma~\ref{lem24} below, as well as the Euler equations $v\cdot\nabla v+\nabla p=0$ themselves. We also point out that, if $\Omega$ is punctured, then $v$ is not assumed to be defined at $0$. It is however always assumed that
\beq\label{hypv1}
{\div}\,v=0\hbox{ in }D,\ \ |v|>0\ \hbox{ in }\Omega,
\eeq
and
\beq\label{hypv2}\left\{\baa{ll}
v\cdot n=0\hbox{ on }\partial\omega_1 & \hbox{if $\Omega$ is not punctured},\vspace{3pt}\\
\displaystyle\int_{C_\epsilon}|v\cdot e_r|\to0\hbox{ as }\epsilon\mathop{\to}^>0& \hbox{if $\Omega$ is punctured},\eaa\right.
\eeq
where $n$ denotes the outward unit normal to $\Omega$ (that is, on $\mathcal{C}_1$ and/or $\mathcal{C}_2$ where appropriate). {\it The assumptions~\eqref{hypv1}-\eqref{hypv2} are made throughout Section~$\ref{sec2}$, and are therefore not repeated in the statements.}\par
The first common result is concerned with the existence and some elementary properties of the stream function $u$.

\begin{lem}\label{lem21}
{\rm{(i)}} There is a unique $($up to additive constants$)$ $C^2(D)$ scalar function~$u$ such that
$$\nabla^\perp u=v\hbox{ in }D,$$
and there is a constant $c_{in}$ such that
\beq\label{cin}\left\{\baa{ll}
u=c_{in}\in\R\hbox{ on $\partial\omega_1$} & \hbox{if $\Omega$ is not punctured},\vspace{3pt}\\
u(x)\to c_{in}\in\overline{\R}\hbox{ as $|x|\displaystyle\mathop{\to}^>0$} & \hbox{if $\Omega$ is punctured}.\eaa\right.
\eeq\par
{\rm{(ii)}} If one further assumes that
\beq\label{C2}\left\{\baa{ll}
v\cdot n=0\hbox{ on }\partial\omega_2 & \hbox{if $\Omega$ is bounded},\vspace{3pt}\\
\displaystyle\int_{C_R}|v\cdot e_r|\to0\hbox{ as }R\to+\infty & \hbox{if $\Omega$ is unbounded},\eaa\right.
\eeq
then there is a constant $c_{out}$ such that
\beq\label{cout}\left\{\baa{ll}
u=c_{out}\in\R\hbox{ on }\partial\omega_2 & \hbox{if $\Omega$ is bounded},\vspace{3pt}\\
u(x)\to c_{out}\in\overline{\R}\hbox{ as $|x|\to+\infty$} & \hbox{if $\Omega$ is unbounded}.\eaa\right.
\eeq
Furthermore,
\beq\label{cinout}
c_{in}\neq c_{out}\ \hbox{ and }\ \min(c_{in},c_{out})<u<\max(c_{in},c_{out})\ \hbox{ in }\Omega.\footnote{The index $in$ in $c_{in}$ refers to the ``inner" boundary of $\Omega$, that is $\partial\omega_1$ (resp. $\{0\}$) if $\Omega$ is not punctured (resp. punctured). Similarly, the index $out$ in $c_{out}$ refers to the ``outer" boundary of $\Omega$, that is $\partial\omega_2$ (resp. the infinity) if $\Omega$ is bounded (resp. unbounded).} 
\eeq 
\end{lem}

\noindent{\bf{Proof.}} (i) The existence and uniqueness (up to additive constants) of a stream function $u\in C^2(D)$ is a consequence of the fact that $\Omega$ is doubly connected, $v$ is divergence free and satisfies the conditions~\eqref{hypv2} on the inner boundary of $\Omega$.\par
Assume now that $\Omega$ is not punctured, and consider a parametrization of  its inner boundary $\mathcal{C}_1$ by a $C^1(\R)$ periodic function $\zeta_1$ such that $|\zeta_1'|>0$ in $\R$. From the first condition in~\eqref{hypv2}, the function $u\circ\zeta_1$ is constant in $\R$, that is, there is $c_{in}\in\R$ such that $u=c_{in}$ on $\mathcal{C}_1=\partial\omega_1$.\par
Assume now that $\Omega$ is punctured, and let $\rho>0$ be such that $\overline{B_r}\!\setminus\!\{0\}\subset\Omega$ for all $r\in(0,\rho)$. Since $u$ has no critical point in $\Omega$ (because $|v|>0$ in $\Omega$ by~\eqref{hypv1}), it follows that
$$u>\min\Big(\min_{C_{r_1}}u,\min_{C_{r_3}}u\Big)\hbox{ in $\Omega_{r_1,r_3}$, for all $0<r_1<r_3<\rho$}.$$
In particular, $\min_{C_{r_2}}u>\min(\min_{C_{r_1}}u,\min_{C_{r_3}}u)$ for all $0<r_1<r_2<r_3<\rho$, hence the function $r\mapsto\min_{C_r}u$ is strictly monotone in a right neighborhood of $0$ and has a limit $c_{in}\in\overline{\R}$ at $0$. On the other hand, the integral condition in~\eqref{hypv2} means that $\int_{C_r}|\nabla u\cdot e_\theta|\to0$ as $r\displaystyle\mathop{\to}^>0$, hence the oscillation of $u$ on $C_r$ converges to $0$ as $r\displaystyle\mathop{\to}^>0$, that is:
\beq\label{osc1}
\mathop{{\rm{osc}}}_{C_r}u:=\max_{C_r}u-\min_{C_r}u\to0\ \hbox{ as }r\displaystyle\mathop{\to}^>0.
\eeq
Finally, $\max_{C_r}u\to c_{in}$ as $r\displaystyle\mathop{\to}^>0$, and $u(x)\to c_{in}$ as $|x|\displaystyle\mathop{\to}^>0$.\par
(ii) Similarly, the conditions~\eqref{C2}, together with~\eqref{hypv1}, imply the existence of $c_{out}\in\overline{\R}$ such that~\eqref{cout} holds in both the bounded and unbounded cases. In the unbounded case, the integral condition in~\eqref{C2} also implies, as in~\eqref{osc1} above, that
\beq\label{osc2}
\max_{C_R}u-\min_{C_R}u\to0\ \hbox{ as }R\to+\infty.
\eeq
Finally, in all cases, since $u$ has no critical point in $\Omega$, the properties~\eqref{cin} and~\eqref{cout} immediately yield~\eqref{cinout}.\hfill$\Box$\break 

The second common result is concerned with the trajectories of the gradient flow $\nabla u$, with~$u$ given by Lemma~\ref{lem21}. Namely, for $x\in D$, with $D$ as in~\eqref{defD}, let $\sigma_x$ be the solution of
\beq\label{defsigmax}\left\{\baa{rcl}
\dot\sigma_x(t) & = & \nabla u(\sigma_x(t)),\vspace{3pt}\\
\sigma_x(0) & = & x.\eaa\right.
\eeq
In the sequel, for $y\in\R^2$ and $E\subset\R^2$, we set
$${\rm{dist}}(y,E)=\inf_{z\in E}|y-z|.$$

\begin{lem}\label{lem21bis}
{\rm{(i)}} If $x\in\Omega$, then there exist $-\infty\!\le\!t^-_x\!<\!0\!<\!t^+_x\!\le\!+\infty$ such that $\sigma_x:(t^-_x,t^+_x)\!\to\!\Omega$ is of class $C^1$ with $(u\circ\sigma_x)'>0$ in $(t^-_x,t^+_x)$ and
\beq\label{claimdist}\left\{\baa{lll}
{\rm{dist}}(\sigma_x(t),\partial\Omega)\to0 & \hbox{as }t\to t_x^\pm & \hbox{if $\Omega$ is bounded},\vspace{3pt}\\
{\rm{dist}}(\sigma_x(t),\partial\Omega)\to0\hbox{ or }|\sigma_x(t)|\to+\infty & \hbox{as }t\to t_x^\pm & \hbox{if $\Omega$ is unbounded}.\eaa\right.
\eeq\par
{\rm{(ii)}} If $\Omega$ is not punctured and $x\in\partial\omega_1$, then either $\sigma_x(t)=x$ for all $t\in\R$, or there is $t^+_x\in(0,+\infty]$ such that $\sigma_x:[0,t^+_x)\to\Omega\cup\partial\omega_1$ is of class $C^1$ with $(u\circ\sigma_x)'>0$ in $[0,t^+_x)$, $\sigma_x((0,t^+_x))\subset\Omega$, and
\beq\label{claimdist2}
{\rm{dist}}(\sigma_x(t),\partial\omega_2)\mathop{\longrightarrow}_{t\to t_x^+}0\hbox{ if $\Omega$ is bounded}\ \big(\hbox{resp. }|\sigma_x(t)|\mathop{\longrightarrow}_{t\to t_x^+}+\infty\hbox{ if $\Omega$ is unbounded}\big),
\eeq
or there is $t^-_x\!\in\![-\infty,0)$ such that $\sigma_x:(t^-_x,0]\to\Omega\cup\partial\omega_1$ is of class $C^1$ with $(u\circ\sigma_x)'\!>\!0$ in~$(t^-_x,0]$, $\sigma_x((t^-_x,0))\subset\Omega$, and~\eqref{claimdist2} holds with $t^-_x$ instead of $t^+_x$.\par
{\rm{(iii)}} If $\Omega$ is bounded, if~\eqref{C2} holds and if $x\in\partial\omega_2$, then either $\sigma_x(t)=x$ for all $t\in\R$, or there is $t^+_x\in(0,+\infty]$ such that $\sigma_x:[0,t^+_x)\to\Omega\cup\partial\omega_2$ is of class $C^1$ with $(u\circ\sigma_x)'>0$ in~$[0,t^+_x)$, $\sigma_x((0,t^+_x))\subset\Omega$, and
\beq\label{claimdist3}
{\rm{dist}}(\sigma_x(t),\partial\omega_1)\mathop{\longrightarrow}_{t\to t_x^+}0\hbox{ if $\Omega$ is not punctured}\ \big(\hbox{resp. }|\sigma_x(t)|\mathop{\longrightarrow}_{t\to t_x^+}0\hbox{ if $\Omega$ is punctured}\big),
\eeq
or there is $t^-_x\!\in\![-\infty,0)$ such that $\sigma_x:(t^-_x,0]\to\Omega\cup\partial\omega_2$ is of class $C^1$ with $(u\circ\sigma_x)'\!>\!0$ in~$(t^-_x,0]$, $\sigma_x((t^-_x,0))\subset\Omega$, and~\eqref{claimdist3} holds with $t^-_x$ instead of $t^+_x$.\end{lem}

\noindent{\bf{Proof.}} (i) Consider any $x$ in $\Omega$. Since $\nabla u$ is of class $C^1(D)$, the solution $\sigma_x$ of~\eqref{defsigmax} is defined in a neighborhood of $0$ and the quantities $t^\pm_x$ defined by
\beq\label{deft-+}\left\{\baa{l}
t_x^+=\sup\big\{t>0: \sigma_x((0,t))\subset\Omega\},\vspace{3pt}\\
t_x^-=\inf\big\{t<0: \sigma_x((t,0))\subset\Omega\}\eaa\right.
\eeq
are such that $-\infty\le t^-_x<0<t^+_x\le+\infty$. The functions $\sigma_x:(t^-_x,t^+_x)\to\Omega$ and $u\circ\sigma_x:(t^-_x,t^+_x)\to\R$ are of class~$C^1$, with
$$(u\circ\sigma_x)'(t)=|\nabla u(\sigma_x(t))|^2=|v(\sigma_x(t))|^2>0\ \hbox{ for all }t\in(t^-_x,t^+_x),$$
since $|v|>0$ in $\Omega$ by~\eqref{hypv1}. Let us then show~\eqref{claimdist} as $t\to t^+_x$ (the limit as $t\to t^-_x$ can be treated similarly). Let us assume by way of contradiction that~\eqref{claimdist} (with $t\to t^+_x$) does not hold. So, in all configurations, there exist an increasing sequence $(t_n)_{n\in\N}$ in $(t^-_x,t^+_x)$ converging to $t^+_x$ and a point $y\in\Omega$ such that $\sigma_x(t_n)\to y$ as $n\to+\infty$. Since the continuous field $|\nabla u|=|v|$ does not vanish in $\Omega$ by~\eqref{hypv1}, there are three real numbers $r>0$, $\eta>0$ and $\tau>0$ such that
$$\left\{\baa{l}
\overline{B(y,r)}\subset\Omega,\ \ |\nabla u|\ge\eta\hbox{ in }\overline{B(y,r)},\vspace{3pt}\\
\sigma_z(t)\in\overline{B(y,r)}\hbox{ for all }z\in\overline{B(y,r/2)}\hbox{ and }t\in[-\tau,\tau].\eaa\right.$$
Since $\sigma_x(t_n)\to y$ as $n\to+\infty$, one has $\sigma_x(t_n)\in\overline{B(y,r/2)}$ for all $n$ large enough, hence~$\sigma_x$ is defined in $[t_n-\tau,t_n+\tau]$ with $\sigma_x(t)\in\overline{B(y,r)}\subset\Omega$ for all $t\in[t_n-\tau,t_n+\tau]$ and~$n$ large enough. This implies that $t^+_x=+\infty$. Furthermore, for all $n$ large enough, one has $(u\circ\sigma_x)'(t)=|\nabla u(\sigma_x(t))|^2\ge\eta^2$ for all~$t\in[t_n-\tau,t_n+\tau]$, hence
$$u(\sigma_x(t_n+\tau))\ge u(\sigma_x(t_n-\tau))+2\eta^2\tau.$$
Since $u\circ\sigma_x$ is increasing in~$(t^-_x,t^+_x)$ and since $t_n\to t^+_x=+\infty$ as $n\to+\infty$, one then gets that $u(\sigma_x(t))\to+\infty$ as $t\to t^+_x=+\infty$, contradicting the fact that $\sigma_x(t_n)\to y$, $t_n\to t^+_x$ and the continuity of $u$ at $y$. Therefore,~\eqref{claimdist} has been proved.\par
(ii) Assume now that $\Omega$ is not punctured and consider $x\in\partial\omega_1$. From~\eqref{hypv2}, either $|\nabla u(x)|\!=\!0$ (and then $\sigma_x(t)=x$ for all $t\in\R$), or $|\nabla u(x)|>0$ and
$$\nabla u(x)=\pm|\nabla u(x)|\,n(x),$$
with~$n(x)$ the outward normal vector to $\Omega$ at $x$. If $\nabla u(x)=-|\nabla u(x)|\,n(x)$ (resp. $\nabla u(x)=|\nabla u(x)|\,n(x)$) with~$|\nabla u(x)|>0$, then $\nabla u(x)$ points inward (resp. outward) $\Omega$ at $x$ and the quantity~$t_x^+$ (resp.~$t^-_x$) given in~\eqref{deft-+} is well defined and satisfies $0<t^+_x\le+\infty$ together with $\sigma_x((0,t_x^+))\subset\Omega$ and $\sigma_x([0,t_x^+))\subset\Omega\cup\partial\omega_1$ (resp. $-\infty\le t^-_x<0$, $\sigma_x((t_x^-,0))\subset\Omega$, and~$\sigma_x((t_x^-,0])\subset\Omega\cup\partial\omega_1$). Furthermore, $\sigma_x$ and $u\circ\sigma_x$ are of class $C^1([0,t^+_x))$ (resp. $C^1((t^-_x,0])$) and $(u\circ\sigma_x)'>0$ in $[0,t^+_x)$ (resp. in $(t^-_x,0]$). Lastly, if $\nabla u(x)=-|\nabla u(x)|\,n(x)$ with $|\nabla u(x)|>0$, then, as in~(i),~\eqref{claimdist} still holds with $t\to t_x^+$ and, since the continuous function $u$ is equal to the constant $c_{in}\in\R$ on $\partial\omega_1$ by Lemma~\ref{lem21}-(i), one gets~\eqref{claimdist2}. If $\nabla u(x)=|\nabla u(x)|\,n(x)$ with~$|\nabla u(x)|>0$, one similarly gets~\eqref{claimdist2} with $t^-_x$ instead of $t^+_x$.\par
(iii) If $\Omega$ is bounded, if~\eqref{C2} holds and if $x\in\partial\omega_2$, then either $|\nabla u(x)|=0$ and $\sigma_x(t)\equiv x$ in~$\R$, or $\sigma_x$ is defined in $(t^-_x,0]$ or $[0,t^+_x)$ and ranges in $\Omega\cup\partial\omega_2$, and the conclusion follows as in~(ii).~\hfill$\Box$
 
\begin{rem}{\rm The extremal values $t_x^\pm$ in Lemma~$\ref{lem21bis}$ can be finite or infinite. Consider for instance the case $\Omega=\Omega_{a,b}$ with $0\le a<b\le\infty$ and a circular flow $v(x)=V(|x|)\,e_\theta(x)$ with a~$C^1((a,b))$ {\it positive} scalar function $V$ (the function $V$ can be assumed to be extended in a~$C^1$ fashion at $a$ and $b$ when they are positive real numbers, and therefore $v$ is of class~$C^1(D)$, with~$D$ as in~\eqref{defD}). The assumptions~\eqref{hypv1}-\eqref{hypv2} and~\eqref{C2} are fulfilled, and the stream function~$u$ is given by $u(x)=U(|x|)$ in $D$, with $U'=V$. For any $x\in\Omega_{a,b}$, the solution~$\sigma_x$ of~\eqref{defsigmax} can then be written as $\sigma_x(t)=\varsigma_x(t)\,e_r(x)$, with $\varsigma_x'(t)=V(\varsigma_x(t))$ for all $t\in(t^-_x,t^+_x)$. Therefore, the finiteness of $t^+_x$ $($resp. of $t^-_x$$)$ is equivalent to the integrability of the function~$1/V$ at $b$ $($resp. at $a$$)$.}
\end{rem}

The third common result provides the existence of a $C^1$ curve in $D$ connecting the inner and outer boundaries of $\Omega$, along which $u$ is strictly monotone, and which will be as in Lemma~\ref{lem21bis} a trajectory of the gradient flow.

\begin{lem}\label{lem22}
Call 
$$S=\big\{x\in D:|v(x)|=0\big\}$$
the set of stagnation points of $v$ in $D$, defined by~\eqref{defD} $($notice that $S\subset\partial\Omega\cap D$ by~\eqref{hypv1}$)$.\par
{\rm{(i)}} Assume that
\beq\label{stagnation}\left\{\baa{ll}
S\subsetneq\mathcal{C}_1\hbox{ or }S\subsetneq\mathcal{C}_2,\hbox{ and }v\cdot n=0\hbox{ on }\mathcal{C}_2 & \hbox{if }\Omega=\omega_2\!\setminus\!\overline{\omega_1},\vspace{3pt}\\
S\subsetneq\mathcal{C}_1 & \hbox{if }\Omega=\R^2\!\setminus\!\overline{\omega_1},\vspace{3pt}\\
S\subsetneq\mathcal{C}_2\hbox{ and }v\cdot n=0\hbox{ on }\mathcal{C}_2 & \hbox{if }\Omega=\omega_2\!\setminus\!\{0\},\vspace{3pt}\\
\displaystyle\displaystyle\int_{C_R}|v\cdot e_r|\to0\hbox{ as }R\to+\infty & \hbox{if }\Omega=\R^2\!\setminus\!\{0\}.\eaa\right.
\eeq
Then there exist $-\infty\le t_{in}<t_{out}\le+\infty$, an interval $I\subset\R$ and a $C^1$ function $\sigma:I\to D$ such that $(u\circ\sigma)'\neq0$ in $I$ and
\beq\label{limitsC12}\left\{\baa{llll}
I=[t_{in},t_{out}], & \sigma(t_{in})\in\partial\omega_1, & \sigma(t_{out})\in\partial\omega_2, & \hbox{if }\Omega=\omega_2\!\setminus\!\overline{\omega_1},\vspace{3pt}\\
I=[t_{in},t_{out}), & \sigma(t_{in})\in\partial\omega_1, & \displaystyle|\sigma(t)|\to+\infty\hbox{ as }t\mathop{\to}^<t_{out}, & \hbox{if }\Omega=\R^2\!\setminus\!\overline{\omega_1},\vspace{3pt}\\
I=(t_{in},t_{out}], & \displaystyle|\sigma(t)|\to0\hbox{ as }t\mathop{\to}^>t_{in}, & \sigma(t_{out})\in\partial\omega_2, & \hbox{if }\Omega=\omega_2\!\setminus\!\{0\},\vspace{3pt}\\
I=(t_{in},t_{out}), & \displaystyle|\sigma(t)|\to0\hbox{ as }t\mathop{\to}^>t_{in}, & \displaystyle|\sigma(t)|\to+\infty\hbox{ as }t\mathop{\to}^<t_{out}, & \hbox{if }\Omega=\R^2\!\setminus\!\{0\}.\eaa\right.
\eeq\par
{\rm{(ii)}} If $\Omega=\R^2\!\setminus\!\overline{\omega_1}$ together with
\beq\label{hypC1infty}
S\subsetneq\mathcal{C}_1\ \hbox{ and }\ \liminf_{|x|\to+\infty}|v(x)|>0,
\eeq
then there exist $-\infty<t_{in}<t_{out}\le+\infty$ and a $C^1$ function $\sigma:[t_{in},t_{out})\to D=\overline{\Omega}$ such that $(u\circ\sigma)'\neq0$ in $[t_{in},t_{out})$ and
$$\sigma(t_{in})\in\partial\omega_1,\ \ |\sigma(t)|\to+\infty\hbox{ as }t\mathop{\to}^<t_{out},\ \hbox{ and }\ |u(\sigma(t))|\to+\infty\hbox{ as }t\mathop{\to}^<t_{out}.$$
Furthermore, one can take $t_{in}=0$ without loss of generality.
\end{lem}

\begin{rem}{\rm Notice that the condition~\eqref{stagnation} implies~\eqref{C2} in all possible configurations of~$\Omega$ except when $\Omega=\R^2\!\setminus\!\overline{\omega_1}$: in that case,~\eqref{stagnation} does not assume that $\lim_{R\to+\infty}\int_{C_R}|v\cdot e_r|=0$, whereas~\eqref{C2} does. We also point out that, again if $\Omega=\R^2\!\setminus\!\overline{\omega_1}$, condition~\eqref{hypC1infty} implies~\eqref{stagnation}, but still does not imply~\eqref{C2}.}
\end{rem}

\noindent{\bf{Proof of Lemma~\ref{lem22}.}} (i) Let us assume here~\eqref{stagnation}. Consider first the case $\Omega=\omega_2\setminus\overline{\omega_1}$, and assume that~$S\subsetneq\mathcal{C}_1$ in~\eqref{stagnation} (the case $S\subsetneq\mathcal{C}_2$ can be handled similarly). Then $|v|=|\nabla u|>0$ on~$\mathcal{C}_2=\partial\omega_2$ and there exists a point
$$A\in\mathcal{C}_1=\partial\omega_1$$
such that $|\nabla u(A)|=|v(A)|>0$ and $\nabla u(A)$ is parallel to the normal vector $n(A)$. Assume first that~$\nabla u(A)$ points inward $\Omega$ at $A$, that is, $\nabla u(A)=-|\nabla u(A)|\,n(A)$, and let
$$\sigma=\sigma_A$$
be the solution of~\eqref{defsigmax} with $x=A$. With the notations of Lemma~\ref{lem21bis}-(ii), the function $\sigma$ is of class $C^1([0,t^+_A))$ and $(u\circ\sigma)'>0$ in $[0,t^+_A)$, with $u(\sigma(0))=u(A)=c_{in}$ by Lemma~\ref{lem21}-(i). Furthermore,~\eqref{claimdist2} implies that ${\rm{dist}}(\sigma(t),\partial\omega_2)\to0$ as $t\to t_A^+$, hence $u(\sigma(t))\to c_{out}$ as $t\to t_A^+$, where $c_{out}\in\R$ is given by Lemma~\ref{lem21}-(ii) (since~\eqref{C2} is satisfied by~\eqref{stagnation}, here with~$\Omega=\omega_2\!\setminus\!\overline{\omega_1}$). In particular, one has $c_{in}<c_{out}$ and $c_{in}<u<c_{out}$ in $\Omega$. On the other hand, the function $|\nabla u\circ\sigma|$ is continuous in~$[0,t^+_A)$, positive at $0$ (since $|\nabla u(\sigma(0))|=|v(A)|>0$), positive in $(0,t^+_A)$ (since~$\sigma((0,t^+_A))\subset\Omega$), and
$$\liminf_{t\to t^+_A}|\nabla u(\sigma(t))|=\liminf_{t\to t^+_A}|v(\sigma(t))|>0$$
(since the continuous field~$|v|$ is positive on the compact set $\mathcal{C}_2=\partial\omega_2$ and ${\rm{dist}}(\sigma(t),\partial\omega_2)\to 0$ as~$t\to t^+_A$). As a consequence, there is~$\eta>0$ such that $|\dot\sigma(t)|=|\nabla u(\sigma(t))|\ge\eta$ for all~$t\in[0,t^+_A)$. Therefore, $(u\circ\sigma)'(t)=|\nabla u(\sigma(t))|^2\ge\eta^2$ for all $t\in[0,t^+_A)$ and $t^+_A$ is a positive real number, since~$u$ is bounded in the compact set $\overline{\Omega}$. Moreover, for every $t\in[0,t^+_A)$, there holds
$$c_{out}-c_{in}\ge u(\sigma(t))-u(A)=u(\sigma(t))-u(\sigma(0))=\int_0^t|\nabla u(\sigma(s))|^2ds\ge\eta\int_0^t|\dot\sigma(s)|\,ds,$$
hence the length of the curve $\sigma([0,t^+_A))$ is finite. Finally, there is a point $A^+\in\mathcal{C}_2=\partial\omega_2$ such that~$\sigma(t)\to A^+$ as $t\to t^+_A$. By setting $\sigma(t^+_A)=A^+$ and remembering that the field~$\nabla u$ is (at least) continuous in $\overline{\Omega}$, it follows that the function $\sigma:[0,t^+_A]\to\overline{\Omega}=D$ is then of class~$C^1([0,t^+_A])$ and $(u\circ\sigma)'>0$ in $[0,t^+_A]$. To sum up, if $\nabla u(A)$ points inward $\Omega$ at $A$, then~\eqref{limitsC12} holds in the case $\Omega=\omega_2\setminus\overline{\omega_1}$ with $t_{in}=0$ and $t_{out}=t^+_A$. Similarly, if $\nabla u(A)$ points outward $\Omega$ at $A$, then~\eqref{limitsC12} still holds, with $t_{in}=0$, $t_{out}=-t^-_A\in(0,+\infty)$, $\sigma(t)=\sigma_A(-t)$ for~$t\in[0,-t^-_A]$, and $\sigma_A$ of class $C^1([t^-_A,0])$. \par
Consider now the case $\Omega=\R^2\!\setminus\!\overline{\omega_1}$. By~\eqref{hypv2} and~\eqref{stagnation}, there is a point $A\in\mathcal{C}_1=\partial\omega_1$ such that $\nabla u(A)$ is not zero and parallel to the normal vector $n(A)$. Assume that $\nabla u(A)$ points inward $\Omega$ at $A$ (the other case can be treated similarly) and let $\sigma=\sigma_A$. The function $\sigma:[0,t^+_A)\to\Omega\cup\partial\omega_1=\overline{\Omega}=D$ is of class $C^1([0,t^+_A))$ and $(u\circ\sigma)'>0$ in $[0,t^+_A)$, for some $t^+_A\in(0,+\infty]$. It then follows from Lemma~\ref{lem21bis}-(ii) (with $t\to t^+_A$) that $|\sigma(t)|\to+\infty$ as $t\to t^+_A$. This yields~\eqref{limitsC12} with $t_{in}=0$ and $t_{out}=t^+_A$.\par
Consider then the case $\Omega=\omega_2\!\setminus\!\{0\}$. By~\eqref{stagnation}, there is a point $A\in\mathcal{C}_2=\partial\omega_2$ such that~$\nabla u(A)$ is not zero and parallel to the normal vector $n(A)$. Assume that $\nabla u(A)$ points outward $\Omega$ at $A$ (the other case can be treated similarly) and let $\sigma=\sigma_A$. The function $\sigma:(t^-_A,0]\to\Omega\cup\partial\omega_2=D$ is of class $C^1((t^-_A,0])$ and $(u\circ\sigma)'>0$ in $(t^-_A,0]$, for some $t^-_A\in[-\infty,0)$. Since~\eqref{C2} is satisfied by~\eqref{stagnation} (here with $\Omega=\omega_2\!\setminus\!\{0\}$), it then follows from Lemma~\ref{lem21bis}-(iii) (with $t\to t^-_A$) that $|\sigma(t)|\to0$ as $t\to t^-_A$. This yields~\eqref{limitsC12} with $t_{in}=t^-_A$ and $t_{out}=0$.\par
Lastly, consider the case $\Omega=\R^2\!\setminus\!\{0\}$. By~\eqref{stagnation}, the condition~\eqref{C2} is fulfilled. By Lemma~\ref{lem21}, there are then $c_{in}\neq c_{out}$ in $\overline{\R}$ such that $u(x)\to c_{in}$ as $\displaystyle|x|\mathop{\to}^>0$ and $u(x)\to c_{out}$ as~$|x|\to+\infty$. Pick any point $A\in\Omega$. By Lemma~\ref{lem21bis}-(i), the function $\sigma_A:(t^-_A,t^+_A)\to\Omega=D$ is of class $C^1((t^-_A,t^+_A))$, with $-\infty\le t^-_A<t^+_A\le+\infty$, and $(u\circ\sigma_A)'>0$ in $(t^-_A,t^+_A)$. Together with~\eqref{claimdist}, one gets that either $|\sigma_A(t)|\to0$ as $t\to t^-_A$ and $|\sigma_A(t)|\to+\infty$ as $t\to t^+_A$, or~$|\sigma_A(t)|\to+\infty$ as $t\to t^-_A$ and $|\sigma_A(t)|\to0$ as $t\to t^+_A$. In the former case,~\eqref{limitsC12} holds with~$\sigma=\sigma_A$ and $(t_{in},t_{out})=(t^-_A,t^+_A)$, whereas in the latter case~\eqref{limitsC12} holds with $\sigma=\sigma_A(-\cdot)$ and~$(t_{in},t_{out})=(-t^+_A,-t^-_A)$.\par
(ii) Assume now that $\Omega=\R^2\!\setminus\!\overline{\omega_1}$ and that~\eqref{hypC1infty} is fulfilled. In particular,~\eqref{stagnation} holds and the previous part~(i) yields the existence of $0=t_{in}<t_{out}\le+\infty$ and a $C^1$ function $\sigma:[t_{in},t_{out})\to\overline{\Omega}$ such that $(u\circ\sigma)'\neq0$ in $[t_{in},t_{out})$, $\sigma(t_{in})=\sigma(0)\in\partial\omega_1$, and $|\sigma(t)|\to+\infty$ as~$\displaystyle t\mathop{\to}^<t_{out}$. Furthermore, by construction, either $\sigma=\sigma_A$, with $A=\sigma(0)\in\partial\omega_1$ such that $|v(A)|>0$, or $\sigma=\sigma_A(-\cdot)$. The only thing to be proven is that $|u(\sigma(t))|\to+\infty$ as $\displaystyle t\mathop{\to}^<t_{out}$. By~\eqref{hypC1infty}, there are $R>0$ and~$\eta>0$ such that $\R^2\!\setminus\!B_R\subset\Omega$ and $|v(x)|\ge\eta$ for all $|x|\ge R$. Let $T\in(0,t_{out})$ such that~$|\sigma(s)|\ge R$ for all $s\in[T,t_{out})$. For all $t\in[T,t_{out})$, one has
$$\big|u(\sigma(t))-u(\sigma(T))\big|=\int_T^t|\nabla u(\sigma(s))|^2ds\ge\eta\int_T^t|\dot\sigma(s)|\,ds\ge\eta\,(|\sigma(t)|-|\sigma(T)|).$$
Consequently, $|u(\sigma(t))|\to+\infty$ as $t\displaystyle\mathop{\to}^{<}t_{out}$. The proof of Lemma~\ref{lem22} is thereby complete.\hfill$\Box$\break

The fourth common result states that the streamlines of the flow are $C^1$ Jordan curves surrounding the origin, and that they approach $0$ or infinity where appropriate. For $x$ in $D$ defined by~\eqref{defD}, $\xi_x$ denotes the solution of~\eqref{defxix}, ranging in $D$ and defined in a maximal interval~$I_x$ containing $0$, and
$$\Xi_x=\xi_x(I_x)$$
denotes the streamline of the flow containing $x$. We recall that, by definition, the stream function $u$ given in Lemma~\ref{lem21} is constant along each streamline of the flow.


\begin{lem}\label{lem23} {\rm{(i)}} If~\eqref{C2} holds, then, for every $x\in\Omega$, the function $\xi_x$ is defined in $\R$ and periodic, and the streamline $\Xi_x=\xi_x(\R)$ is a $C^1$ Jordan curve surrounding the origin. Furthermore,
$$\left\{\baa{lll}
\displaystyle\min_{\R}|\xi_x|\to+\infty & \hbox{as $|x|\to+\infty$} & \hbox{if $\Omega$ is unbounded},\vspace{3pt}\\
\displaystyle\max_{\R}|\xi_x|\to0 & \displaystyle\hbox{as }|x|\mathop{\to}^>0 & \hbox{if $\Omega$ is punctured}.\eaa\right.$$\par
{\rm{(ii)}} If~\eqref{C2} holds with $\Omega$ unbounded and $\liminf_{|x|\to+\infty}|v(x)\cdot e_\theta(x)|>0$, then
$$\max_\R|\xi_x|-\min_\R|\xi_x|\to0\hbox{ as $|x|\to+\infty$}.$$\par
{\rm{(iii)}} If $\Omega=\R^2\!\setminus\!\overline{\omega_1}$ and $\inf_\Omega|v|>0$, then, for every $x\in\overline{\Omega}$, the function $\xi_x$ is defined in~$\R$ and periodic, and the streamline $\Xi_x=\xi_x(\R)$ is a $C^1$ Jordan curve surrounding the origin. Furthermore, $\min_\R|\xi_x|\to+\infty$ and $|u(x)|\to+\infty$ as $|x|\to+\infty$, and $u-c_{in}$ has a constant strict sign in $\Omega$, where $u=c_{in}$ on $\partial\omega_1$.
\end{lem}

\noindent{\bf{Proof.}} (i) From Lemma~\ref{lem21}, there are $c_{in}\neq c_{out}$ in $\overline{\R}$ such that the stream function $u$ satisfies~\eqref{cin} and~\eqref{cout}-\eqref{cinout}. Therefore, together with the continuity of $u$ in $D$ and the fact that $u$ is constant along each streamline of the flow, it follows that, for each $x\in\Omega$, $\inf_{t\in I_x}{\rm{dist}}(\xi_x(t),\partial\Omega)>0$ and $\Xi_x=\xi_x(I_x)$ is bounded. Since $|\nabla u|=|v|>0$ in $\Omega$ by~\eqref{hypv1}, it is standard to conclude that, for each $x\in\Omega$, the function $\xi_x$ is periodic and the streamline $\Xi_x=\xi_x(\R)$ (here, $I_x=\R$) is a $C^1$ Jordan curve surrounding the origin.\par
If $\Omega$ is unbounded, then $u(x)\to c_{out}\in\overline{\R}$ as $|x|\to+\infty$ by~\eqref{cout}. Together with~\eqref{cin},~\eqref{cinout}, and the continuity of $u$ in $D$, it easily follows that $\min_{t\in\R}|\xi_x(t)|\to+\infty$ as $|x|\to+\infty$. Similarly, if $\Omega$ is punctured, then $\max_{t\in\R}|\xi_x(t)|\to0$ as $\displaystyle|x|\mathop{\to}^>0$.\par
(ii) In the unbounded case, assume now that $\liminf_{|x|\to+\infty}|v(x)\cdot e_\theta(x)|>0$, in addition to~\eqref{C2}. From~(i), one knows that, for each $x\in\Omega$, the streamline $\Xi_x$ is a $C^1$ Jordan curve surrounding the origin, with $\min_{\R}|\xi_x|\to+\infty$ as $|x|\to+\infty$. Let us now show that $\max_\R|\xi_x|-\min_\R|\xi_x|\to0$ as $|x|\to+\infty$. Let $R>0$ and $\eta>0$ be such that $\R^2\!\setminus\!B_R\subset\Omega$ and $|\nabla u(x)\cdot e_r(x)|=|v(x)\cdot e_\theta(x)|\ge\eta$ for all $|x|\ge R$. Together with the continuity of $\nabla u$, let us only consider the case
\beq\label{nablauer}
\nabla u(x)\cdot e_r(x)\ge\eta\ \hbox{ for all $|x|\ge R$}
\eeq
(the case $\nabla u(x)\cdot e_r(x)\le-\eta$ can be handled similarly). Consider any $\epsilon>0$. From~\eqref{osc2} and~\eqref{nablauer}, there is $R_\epsilon\ge R+\epsilon$ such that
$$\max_{C_{|x|-\epsilon}}u<u(x)-\frac{\epsilon\,\eta}{2}\ \hbox{ and }\ \min_{C_{|x|+\epsilon}}u>u(x)+\frac{\epsilon\,\eta}{2},\ \hbox{ for all }|x|\ge R_\epsilon.$$
Thus, for every $x$ with $|x|\ge R_\epsilon$, one has $\Xi_x\subset\Omega_{|x|-\epsilon,|x|+\epsilon}$ and $\max_\R|\xi_x|-\min_\R|\xi_x|<2\epsilon$.\par
(iii) Consider now the case $\Omega=\R^2\!\setminus\!\overline{\omega_1}$ and assume that
\beq\label{defeta0}
\eta:=\inf_\Omega|v|>0.
\eeq
Notice that we do not assume the condition~\eqref{C2}, so Lemma~\ref{lem21}-(ii) and parts~(i)-(ii) of the present lemma can not be applied. In particular, we do not know yet that $u$ has a limit at infinity or that all streamlines of the flow surround the origin. However, Lemma~\ref{lem22}-(ii) can be applied, since~\eqref{hypC1infty} is fulfilled here. Therefore, there exist $t_{out}\in(0,+\infty]$, a point $A\in\partial\omega_1$ (actually, $A$ can here be arbitrarily chosen on $\partial\omega_1$ since $|v|>0$ on $\partial\omega_1$), and a $C^1$ function $\sigma:[0,t_{out})\to D$ such that $\sigma(0)=A$, $(u\circ\sigma)'\neq0$ in $[0,t_{out})$, and $|\sigma(t)|\to+\infty$ and $|u(\sigma(t))|\to+\infty$ as $\displaystyle t\mathop{\to}^<t_{out}$. Let us assume without loss of generality that $(u\circ\sigma)'>0$ in~$[0,t_{out})$ (the other case can be handled similarly), hence
\beq\label{usigmat}
u(\sigma(t))\to+\infty\ \hbox{ as }\displaystyle t\mathop{\to}^<t_{out}.
\eeq
By construction of $\sigma$ in Lemma~\ref{lem22}-(ii), this case corresponds to the condition $\nabla u(A)\cdot n(A)\!<\!0$.\par
Denote
$$E=\big\{s\in[0,t_{out}):\hbox{the streamline }\Xi_{\sigma(s)}\hbox{ is a $C^1$ Jordan curve surrounding the origin}\big\},$$
and let us show that
\beq\label{claimE}
E=[0,t_{out}).
\eeq
To do so, we prove that $E$ is not empty (it contains~$0$), open relatively to $[0,t_{out})$ and that the largest interval containing $0$ and contained in $E$ is actually equal to $[0,t_{out})$. Note first that, since $v\cdot n=0$ and $|v|>0$ on $\mathcal{C}_1=\partial\omega_1$, the streamline $\Xi_{\sigma(0)}=\Xi_A$ is equal to the $C^1$ Jordan curve~$\mathcal{C}_1$ and it surrounds the origin by assumption. In other words, $0\in E$.\par
Let us now show that $E$ is open relatively to $[0,t_{out})$. Let $s_0\in E$ and denote $x=\sigma(s_0)\in\overline{\Omega}$. By definition, the function $\xi_x$ is periodic, with some period $T_x>0$. Remember also that $u$ is constant along each streamline of the flow. Therefore, since $v$ is (at least) continuous and $|v(x)|=|\nabla u(x)|>0$ in $\overline{\Omega}$, there are some real numbers $r>0$ and $\tau\in(0,T_x)$ such that, for every~$y\in B(x,r)\cap\overline{\Omega}$, there are some real numbers $t^\pm_y$ such that
$$-\tau<t^-_y<0<t^+_y<\tau\ \hbox{ and }\ B(x,r)\cap\Xi_y=B(x,r)\cap u^{-1}(\{u(y)\})=\xi_y((t^-_y,t^+_y)).$$
On the other hand, since $\xi_x(T_x)=\xi_x(0)=x$, the Cauchy-Lipschitz theorem provides the existence of a real number $r'\in(0,r]$ such that, for every $z\in B(x,r')\cap\overline{\Omega}$, the function~$\xi_z$ is defined (and of class $C^1$) at least on the interval $[0,T_x]$ and $\xi_z(T_x)\in B(x,r)\cap\overline{\Omega}$. Furthermore, by continuity of $\sigma$, there is $\varepsilon>0$ such that $s_0+\varepsilon<t_{out}$ and
$$\sigma(s)\in B(\sigma(s_0),r')\cap\overline{\Omega}= B(x,r')\cap\overline{\Omega}\ \hbox{ for all }s\in[\max(0,s_0-\varepsilon),s_0+\varepsilon].$$
As a consequence, for every $s\in[\max(0,s_0-\varepsilon),s_0+\varepsilon]$, the points $z:=\sigma(s)\in B(x,r')\cap\overline{\Omega}$ and~$y:=\xi_z(T_x)\in B(x,r)\cap\overline{\Omega}$ satisfy $u(z)=u(y)$, hence
$$z\in B(x,r')\cap u^{-1}(\{u(y)\})\ \subset\ B(x,r)\cap u^{-1}(\{u(y)\})$$
and $z=\xi_y(t)$ for some $t\in(t^-_y,t^+_y)\,(\subset(-\tau,\tau))$. Thus, $\xi_y(-T_x)=z=\xi_y(t)$ and since $|t|<\tau<T_x$, the function $\xi_y$ is defined in $\R$ and $(T_x+t)$-periodic. So is $\xi_z$ since $z\in\Xi_y$. In other words, for every $s\in[\max(0,s_0-\varepsilon),s_0+\varepsilon]$, the function $\xi_{\sigma(s)}=\xi_z$ is defined in $\R$ and periodic. Since $|\nabla u|=|v|>0$ in $\overline{\Omega}$, one then concludes that $\Xi_{\sigma(s)}$ is a $C^1$ Jordan curve surrounding the origin. Finally, the set~$E$ is open relatively to $[0,t_{out})$.\par
Denote
$$T_*=\sup\big\{t\in[0,t_{out}):[0,t]\subset E\big\}.$$
The previous paragraphs imply that $0<T_*\le t_{out}$. The proof of~\eqref{claimE} will be complete once we show that $T_*=t_{out}$. Assume by way of contradiction that $T_*<t_{out}$ (in particular, $T_*$ is then a positive real number). Consider any increasing sequence $(s_n)_{n\in\N}$ in $(0,T_*)$ and converging to $T_*$. Owing to the definition of $T_*$, each function $\xi_{\sigma(s_n)}$ is periodic and each streamline $\Xi_{\sigma(s_n)}$ surrounds the origin. Furthermore, since each $s_n$ is positive and $u\circ\sigma$ is increasing in $[0,t_{out})$ and $u$ is constant on~$\mathcal{C}_1\,(\ni A=\sigma(0))$, each streamline $\Xi_{\sigma(s_n)}$ is included in the (open) set $\Omega$. Consider now any~$n\in\N$ and any point
$$x\in\Xi_{\sigma(s_n)}.$$
Notice that $u(x)=u(\sigma(s_n))$ and remember that $u\circ\sigma$ is increasing in $[0,t_{out})$, hence
$$u(A)=u(\sigma(0))<u(\sigma(s_n))=u(x)<u(\sigma(T_*)).$$
Moreover, by Lemma~\ref{lem21bis}-(i), there exist $-\infty\le t^-_x\!<0<t^+_x\le+\infty$ such that $\sigma_x:(t^-_x,t^+_x)\to\Omega$ is of class $C^1$ with
\beq\label{dist4}
(u\circ\sigma_x)'>0\hbox{ in $(t^-_x,t^+_x),\ $ and }\ {\rm{dist}}(\sigma_x(t),\partial\omega_1)\to0\hbox{ or }|\sigma_x(t)|\to+\infty\ \hbox{ as }t\to t_x^\pm.
\eeq
The non-zero vector $\dot\sigma_x(0)=\nabla u(\sigma_x(0))=\nabla u(x)$ is orthogonal to~$\Xi_{\sigma(s_n)}$ at $x$ by definition of~$u$. Since $u(\sigma_x(0))=u(x)>u(A)=c_{in}$ with $u=c_{in}$ on $\mathcal{C}_1=\partial\omega_1$, and since $\Xi_{\sigma(s_n)}$ is a $C^1$ Jordan curve surrounding the origin and meeting orthogonally $\sigma_x((t^-_x,t^+_x))$ at the only point~$x$, it then follows from~\eqref{dist4} that ${\rm{dist}}(\sigma_x(t),\partial\omega_1)\to 0$ as $t\to t^-_x$ and $u(\sigma_x(t))>c_{in}$ for all $t\in(t^-_x,t^+_x)$. Then, for any $t\in(t^-_x,0)$, there holds
$$\baa{rcl}
\displaystyle u(\sigma(T_*))>u(x)=u(\sigma_x(0))\!=\!\int_t^0|\nabla u(\sigma_x(s))|^2ds+u(\sigma_x(t)) & \!\!\!\ge\!\!\! & \displaystyle\eta\int_t^0|\dot\sigma_x(s)|\,ds+c_{in}\vspace{3pt}\\
& \!\!\!\ge\!\!\! & \eta\,\big(|x|\!-\!|\sigma_x(t)|\big)+c_{in},\eaa$$
with $\eta>0$ as in~\eqref{defeta0}. At the limit as $t\to t^-_x$, one gets that $|x|\le R_1+(u(\sigma(T_*))-c_{in})/\eta$, with $R_1:=\max_{y\in\mathcal{C}_1}|y|$. This property holds for any~$n\in\N$ and any $x\in\Xi_{\sigma(s_n)}$, hence
$$\sup_{n\in\N}\Big(\max_\R|\xi_{\sigma(s_n)}|\Big)\le R_1+\frac{u(\sigma(T_*))-c_{in}}{\eta}=:M.$$\par
Finally, consider the streamline $\Xi_{\sigma(T_*)}$ parametrized by the function $\xi_{\sigma(T_*)}$. If there is a real number~$t$ such that $|\xi_{\sigma(T_*)}(t)|>M$, then $|\xi_{\sigma(s_n)}(t)|>M$ for all $n$ large enough, by the Cauchy-Lipschitz theorem. Therefore, $\Xi_{\sigma(T_*)}\subset\overline{B_M}$ and, since $|v|>0$ in $\overline{\Omega}$, it follows that~$\xi_{\sigma(T_*)}$ is defined in $\R$ and periodic, and that $\Xi_{\sigma(T_*)}$ is a $C^1$ Jordan curve surrounding the origin. In other words,~$T_*\in E$. Since $E$ has been proved to be open relatively to $[0,t_{out})$, one is led to a contradiction with the definition of $T_*$ if $T_*<t_{out}$. Eventually, $T_*=t_{out}$ and~\eqref{claimE} is thereby proved.\par
We now claim that
\beq\label{claimsigma}
\min_{\R}|\xi_{\sigma(s)}|\to+\infty\ \hbox{ as }s\displaystyle\mathop{\to}^{<}t_{out}.
\eeq
Indeed, for any $R>R_1$, let $C\in[0,+\infty)$ be such that $|u|\le C$ in $\overline{\Omega}\cap\overline{B_R}$. By~\eqref{usigmat}, there is then $\tau\in(0,t_{out})$ such that $u(\sigma(s))>C$ for all $s\in(\tau,t_{out})$, hence $u(\xi_{\sigma(s)}(t))=u(\sigma(s))>C$ and $|\xi_{\sigma(s)}(t)|>R$ for all $s\in(\tau,t_{out})$ and $t\in\R$. Thus, $\min_{\R}|\xi_{\sigma(s)}|>R$ for all $s\in(\tau,t_{out})$. This yields~\eqref{claimsigma}.\par
Consider then any point $x\in\overline{\Omega}$ and let us deduce that $\Xi_x$ surrounds the origin. From~\eqref{usigmat}-\eqref{claimE} and~\eqref{claimsigma}, there is $s\in(0,t_{out})$ such that $u(\sigma(s))>u(x)$ and the streamline $\Xi_{\sigma(s)}$ is a $C^1$ Jordan curve surrounding both $x$ and the origin. Therefore, the streamline $\Xi_x$ is bounded (it belongs to the bounded connected component of $\R^2\!\setminus\!\Xi_{\sigma(s)}$). Using again that $|v|>0$ in $\overline{\Omega}$, one then concludes that the function~$\xi_x$ is periodic and that $\Xi_x$ is a $C^1$ Jordan curve surrounding the origin.\par
From the previous properties, we then easily get that $\min_{\R}|\xi_x|\to+\infty$ as $|x|\to+\infty$. Indeed, for any fixed $R>R_1=\max_{y\in\partial\omega_1}|y|$, there is $s\in(0,t_{out})$ such that $\min_{\R}|\xi_{\sigma(s)}|>R$, by~\eqref{claimsigma}. Then, for any $x$ with $|x|>R':=\max_{\R}|\xi_{\sigma(s)}|>R$, the streamlines~$\Xi_x$ and~$\Xi_{\sigma(s)}$ do not intersect, and both of them are $C^1$ Jordan curves surrounding the origin, hence $\min_\R|\xi_x|>\min_\R|\xi_{\sigma(s)}|>R$. This shows that $\min_{\R}|\xi_x|\to+\infty$ as $|x|\to+\infty$.\par
Finally, let us prove that $u(x)\to+\infty$ as $|x|\to+\infty$. Fix any $B>0$, and, by~\eqref{usigmat}, let~$s\in(0,t_{out})$ such that $u\circ\sigma>B$ in $(s,t_{out})$. Let then $R_B>R_1$ be such that $\min_\R|\xi_x|>\max_{t\in[0,s]}|\sigma(t)|$ for all $|x|\ge R_B$. For every $x$ such that $|x|\ge R_B$, the streamline $\Xi_x$ surrounds the origin and necessarily crosses $\sigma([0,t_{out}))$, at a point $\sigma(s_x)$ with $s_x\in(s,t_{out})$, hence $u(x)=u(\sigma(s_x))>B$. This shows that $u(x)\to+\infty$ as $|x|\to+\infty$. Since $u=c_{in}$ on $\partial\omega_1$ and $u$ has no critical point in $\Omega$, one then concludes that $u>c_{in}$ in $\Omega$. The proof of Lemma~\ref{lem23} is thereby complete.\hfill$\Box$

\begin{rem}{\rm If $|v|>0$ on $\mathcal{C}_1=\partial\omega_1$ (resp. on $\mathcal{C}_2=\partial\omega_2$), then the boundary conditions~\eqref{hypv2} (resp.~\eqref{C2}) imply that, for any $x\in\partial\omega_1$ with $\Omega$ not punctured, (resp. $x\in\partial\omega_2$ with $\Omega$ bounded), $\xi_x$ is still defined and periodic in $\R$ with $\Xi_x=\partial\omega_1$ (resp.~$\Xi_x=\partial\omega_2$). If~$x\in\partial\Omega\cap D$ and $|v(x)|=0$, then $\xi_x(t)=x$ for all $t\in\R$ and~$\Xi_x=\{x\}$. If~$x\in\partial\omega_1$ with $\Omega$ not punctured (resp. $x\in\partial\omega_2$ with $\Omega$ bounded) with $|v(x)|>0$ and if $v$ has some stagnation points on~$\partial\omega_1$ (resp. $\partial\omega_2$), then $\xi_x$ is still defined in $\R$, but it is not periodic anymore and $\Xi_x$ is a proper arc of $\partial\omega_1$ (resp. $\partial\omega_2$) which is open relatively to $\partial\omega_1$ (resp. $\partial\omega_2$).}
\end{rem}

In the last common preliminary result, we derive a semilinear elliptic equation $\Delta u+f(u)=0$ in $D$ for some $C^1$ function $f$ defined in the range of $u$. We recall that $D$ is defined in~\eqref{defD}.

\begin{lem}\label{lem24} Assume that $v$ is of class $C^2(D)$ and that $v$ solves the Euler equations~\eqref{1} in~$D$, still together with~\eqref{hypv1}-\eqref{hypv2}. Let $u$, in $C^3(D)$, be given by Lemma~$\ref{lem21}$ and let $J$ its range defined by
$$J=\{u(x):x\in D\}.$$\par
{\rm{(i)}} If~\eqref{C2} and~\eqref{stagnation} are fulfilled, then there is a $C^1$ function $f:J\to\R$ such that
$$\Delta u+f(u)=0\ \hbox{ in $D$}.$$\par
{\rm{(ii)}} If $\Omega=\R^2\!\setminus\!\overline{\omega_1}$ and $\inf_\Omega|v|>0$, then the same conclusion holds.
\end{lem}

\noindent{\bf{Proof.}} First of all, since $v$ is of class $C^2(D)$, the stream function $u$ given by Lemma~\ref{lem21} is now of class $C^3(D)$. By continuity of $u$ and connectedness of $D$, the range $J$ of $u$ is an interval.\par
(i) Assume here~\eqref{C2} and~\eqref{stagnation}. It follows from Lemma~\ref{lem21} and its notations that the interior of $J$ is equal to $(\min(c_{in},c_{out}),\max(c_{in},c_{out}))$. Furthermore, $J$ is open at $c_{in}$ if and only if $\Omega$ is punctured, while $J$ is open at $c_{out}$ if and only if $\Omega$ is unbounded (for instance, if~$\Omega=\R^2\setminus\overline{\omega_1}$, then $J=[c_{in},c_{out})$ or $(c_{out},c_{in}]$). It then follows from Lemmas~\ref{lem21} and~\ref{lem22}-(i) that, with the same notations as there, the function $g:=u\circ\sigma:I\to J$ is a $C^1$ diffeomorphism from $I$ onto $J$. Let $g^{-1}:J\to I$ be its~$C^1$ reciprocal diffeomorphism, and define
\beq\label{deff}
f(\tau)=-\Delta u(\sigma(g^{-1}(\tau)))\ \hbox{ for }\tau\in J.
\eeq
By the chain rule, the function $f$ is of class $C^1(J)$ (remember that $\Delta u$ is now of class $C^1(D)$). The above formula means that the equation $\Delta u+f(u)=0$ is satisfied along the curve $\sigma(I)$. Let us now check it in the whole set $D$. Consider first any point $x\in\Omega$. From Lemmas~\ref{lem22}-(i) and~\ref{lem23}-(i), the streamline $\Xi_x$ surrounds the origin and meets the curve $\sigma(I)$. Hence, there is $s\in I$ such that $\sigma(s)\in\Xi_x$. On the one hand, the stream function $u$ is constant along the streamline $\Xi_x$. On the other hand, the $C^1(D)$ vorticity $\frac{\partial v_2}{\partial x_1}-\frac{\partial v_1}{\partial x_2}=\Delta u$ satisfies $v\cdot\nabla(\Delta u)=0$ in $D$ from the Euler equations~\eqref{1}, hence~$\Delta u$ is constant along the streamline~$\Xi_x$ too. As a consequence,~\eqref{deff} yields
$$\Delta u(x)+f(u(x))=\Delta u(\sigma(s))+f(u(\sigma(s)))=\Delta u(\sigma(s))+f(g(s))=0.$$
Therefore, $\Delta u+f(u)=0$ in $\Omega$. Finally, since both functions $\Delta u$ and $f\circ u$ are (at least) continuous in $D$, one concludes that $\Delta u+f(u)=0$ in $D$.\par
(ii) Assume now that $\Omega=\R^2\!\setminus\!\overline{\omega_1}$ and $\inf_\Omega|v|>0$. From Lemmas~\ref{lem21}-(i) and~\ref{lem23}-(iii), one then has $J=[c_{in},+\infty)$ or $(-\infty,c_{in}]$. Together with Lemma~\ref{lem22}-(ii), there are then~$t_{out}$ in~$(0,+\infty]$ and a $C^1$ curve $\sigma:[0,t_{out})\to D=\overline{\Omega}$ such that the function $g:=u\circ\sigma$ is a $C^1$ diffeomorphism from $[0,t_{out})$ onto $J$. Let $g^{-1}:J\to[0,t_{out})$ be its~$C^1$ reciprocal diffeomorphism, and define $f\in C^1(J)$ as in~\eqref{deff}. Since, by Lemma~\ref{lem23}-(iii), each streamline~$\Xi_x$ (for each $x\in D=\overline{\Omega}$) surrounds the origin, it then follows as in~(i) above that the equation $\Delta u+f(u)=0$ holds, here directly in $D$. The proof of Lemma~\ref{lem24} is thereby complete.\hfill$\Box$


\SE{Proof of the main results in fixed annular domains $\Omega_{a,b}$}\label{sec3}

This section is devoted to the proof of Theorems~\ref{th1}-\ref{th3} and~\ref{th4} on the Euler flows in the fixed annular domains $\Omega_{a,b}$ with $0\le a<b\le\infty$. The proofs rely on the common properties proved in Section~\ref{sec2}, as well as on various applications of Proposition~\ref{promoving} and further specific arguments in the unbounded and punctured cases.


\subsection{The case of bounded annuli $\Omega_{a,b}$: proof of Theorems~\ref{th1} and~\ref{th1bis}}\label{sec31}

This section is devoted to the proof of Theorem~\ref{th1bis} (we recall that Theorem~\ref{th1} is a particular case of Theorem~\ref{th1bis}). Throughout this section, we consider two positive real numbers $a<b$ and a~$C^2(\overline{\Omega_{a,b}})$ solution~$v$ of~\eqref{1}-\eqref{2} satisfying~\eqref{hypstagnation}, namely
$$\big\{x\in\overline{\Omega_{a,b}}: |v(x)|=0\big\}\subsetneq C_a\ \hbox{ or }\ \big\{x\in\overline{\Omega_{a,b}}: |v(x)|=0\big\}\subsetneq C_b.$$
This situation falls within the general framework of Section~\ref{sec2}, with $\mathcal{C}_1=C_a$, $\omega_1=B_a$, $\mathcal{C}_2=C_b$, $\omega_2=B_b$, $\Omega=\omega_2\!\setminus\!\overline{\omega_1}=\Omega_{a,b}$, and $D=\overline{\Omega_{a,b}}$. Notice also that the conditions~\eqref{hypv1}-\eqref{hypv2},~\eqref{C2}, and~\eqref{stagnation} are fulfilled by assumption. Therefore, the flow has a $C^3(\overline{\Omega_{a,b}})$ stream function~$u$ and, by Lemmas~\ref{lem21} and~\ref{lem24}-(i), there are two real numbers $c_{in}\neq c_{out}$ and a~$C^1([\min(c_{in},c_{out}),\max(c_{in},c_{out})])$ function $f$ such that
\beq\label{eqannulus}\left\{\baa{l}
\Delta u+f(u)=0\hbox{ in }\overline{\Omega_{a,b}},\vspace{3pt}\\
\min(c_{in},c_{out})<u<\max(c_{in},c_{out})\hbox{ in }\Omega_{a,b},\vspace{3pt}\\
u=c_{in}\hbox{ on }C_a,\ \ u=c_{out}\hbox{ on }C_b.\eaa\right.
\eeq
It then follows from~\cite[Theorem~5]{si}\footnote{Notice that this result holds in any dimension $n\ge2$. It is similar to the classical radial symmetry property proved in~\cite{gnn} in the case where $u$ is a positive solution of the equation $\Delta u+f(u)=0$ in a ball, with Dirichlet condition $u=0$ on the boundary.} that $u$ is radially symmetric and strictly monotone with respect to $|x|$ in $\overline{\Omega_{a,b}}$. Therefore, there is a $C^3([a,b])$ strictly monotone function $U:[a,b]\to\R$ such that $u(x)=U(|x|)$ for all $x\in\overline{\Omega_{a,b}}$. The flow $v=\nabla^{\perp}u$ is then given by
$$v(x)=V(|x|)\,e_\theta(x)$$
for all $x\in\overline{\Omega_{a,b}}$, with $V=U'\in C^2([a,b])$. Lastly, since $|v|$ is continuous in $\overline{\Omega_{a,b}}$ and does not vanish in $\Omega_{a,b}$ nor in the whole circle $C_a$ nor in the whole circle $C_b$, the function $V$ then has a constant strict sign in $[a,b]$. The proof of Theorem~\ref{th1bis} is thereby complete.\hfill$\Box$

\subsubsection*{A related open question}

For a $C^2(\overline{\Omega_{a,b}})$ flow $v$ solving~\eqref{1}-\eqref{2}, could the assumption~\eqref{hypstagnation} be slightly relaxed for $v$ still to be necessarily a circular flow? As we mentioned in the introduction, the conclusion does not hold in general if $v$ has stagnation points in $\Omega_{a,b}$. So a natural question is the following one:
\beq\label{q1}
\hbox{if $|v|>0$ in $\Omega_{a,b}$, then is $v$ a circular flow~?}
\eeq
We first point out that Lemmas~\ref{lem21},~\ref{lem21bis} and~\ref{lem23}-(i), with $\Omega=\Omega_{a,b}$, still hold since they do not use the whole assumption~\eqref{stagnation} (more precisely, they do not use~\eqref{hypstagnation}), but only $|v|>0$ in $\Omega_{a,b}$. Consider then any point~$y\in\Omega_{a,b}$. With the same notations as in Lemma~\ref{lem21bis}, and assuming without loss of generality that $c_{in}<c_{out}$ (after possibly changing $v$ into $-v$ and $u$ into $-u$), there are some quantities~$t^\pm_y$ such that $-\infty\le t^-_y<0<t^+_y\le+\infty$ and the solution~$\sigma_y$ of~\eqref{defsigmax} with $y$ instead of $x$ is of class~$C^1((t^-_y,t^+_y))$ and ranges in~$\Omega_{a,b}$, with
\beq\label{convab}\left\{\baa{l}
|\sigma_y(t)|\to a\hbox{ and }u(\sigma_y(t))\to c_{in}\hbox{ as }t\to t^-_y,\vspace{3pt}\\
|\sigma_y(t)|\to b\hbox{ and }u(\sigma_y(t))\to c_{out}\hbox{ as }t\to t^+_y.\eaa\right.
\eeq
The $C^1((t^-_y,t^+_y))$ function $g:=u\circ\sigma_y$ is increasing (we recall that $(u\circ\sigma_y)'(t)=|\nabla u(\sigma_y(t))|^2=|v(\sigma_y(t))|^2>0$ for all $t\in(t^-_y,t^+_y)$), and $g$ is then an increasing homeomorphism from $(t^-_y,t^+_y)$ onto $(c_{in},c_{out})$. The function $f:(c_{in},c_{out})\to\R$ defined by
\beq\label{deff2}
f(\tau)=-\Delta u(\sigma_y(g^{-1}(\tau)))\ \hbox{ for }\tau\in(c_{in},c_{out})
\eeq
is of class $C^1((c_{in},c_{out}))$ and, since for every $x\in\Omega_{a,b}$ the streamline $\Xi_x$ intersects $\sigma_y((t^-_y,t^+_y))$ by Lemma~\ref{lem23}-(i), the same arguments as in the proof of Lemma~\ref{lem24}-(i) imply that
$$\Delta u+f(u)=0\hbox{ in }\Omega_{a,b}.$$ 
Furthermore, remembering from Lemma~\ref{lem23}-(i) that, for each $x\in\Omega_{a,b}$, the $C^1$ solution $\xi_x$ of~\eqref{defxix} is periodic and ranges in $\Omega_{a,b}$, it then follows from the continuity of $u$ in the compact set $\overline{\Omega_{a,b}}$ and the facts that $u=c_{in}$ on $C_a$, $u=c_{out}$ on $C_b$ and $c_{in}<u<c_{out}$ in $\Omega_{a,b}$, that
$$\max_{t\in\R}|\xi_x(t)|\to a\hbox{ as }|x|\mathop{\to}^>a\ \hbox{ and }\ \min_{t\in\R}|\xi_x(t)|\to b\hbox{ as }|x|\mathop{\to}^<b.$$
Since the function $\Delta u$ is constant along any streamline of the flow from the Euler equations~\eqref{1} and since $\Delta u$ is uniformly continuous in $\overline{\Omega_{a,b}}$, it then follows from the previous observations that $\Delta u$ is constant on $C_a$ and constant on $C_b$. Call~$d_1$ and~$d_2$ the values of $\Delta u$ on~$C_a$ and~$C_b$, respectively, and set $f(c_{in})=-d_1$ and $f(c_{out})=-d_2$. One then infers from~\eqref{convab}-\eqref{deff2} that $f:[c_{in},c_{out}]\to\R$ is continuous in $[c_{in},c_{out}]$ and that the equation~$\Delta u+f(u)=0$ holds in the closed annulus $\overline{\Omega_{a,b}}$ ($u$ is then a classical $C^2(\overline{\Omega_{a,b}})$ solution of~\eqref{eqannulus}). However, since
$$f'(\tau)=-\frac{\nabla(\Delta u)(\sigma_y(g^{-1}(\tau)))\cdot\nabla u(\sigma_y(g^{-1}(\tau)))}{|\nabla u(\sigma_y(g^{-1}(\tau)))|^2}\ \hbox{ for all }\tau\in(c_{in},c_{out})$$
and since $|\nabla u(\sigma_y(g^{-1}(\tau)))|$ can converge to $0$ as $\tau\to c_{in}$ or $c_{out}$ (this happens if $|v|=0$ on $C_a$ or if $|v|=0$ on $C_b$), the function~$f'$ can be unbounded in $(c_{in},c_{out})$.\footnote{For instance, the smooth flow $v(x)=(|x|-a)\,e_\theta(x)$ solves~\eqref{1}-\eqref{2} with pressure $p(x)=|x|^2/2-2a|x|+a^2\ln|x|$ and stream function $u(x)=(|x|-a)^2/2$ (up to additive constants), while $|v|=0$ on $C_a$ and $|v|>0$ in~$\Omega_{a,b}$. Here, $c_{in}=0$, $c_{out}=(b-a)^2/2$ and $f(s)=-2+a/(a+\sqrt{2s})$ for $s\in[c_{in},c_{out}]=[0,(b-a)^2/2]$, hence~$f'$ is not bounded in $(c_{in},c_{out})$. Notice that this example is a circular flow, which makes question~\eqref{q1} still relevant.} The argument used in the proof of Theorem~\ref{th1bis} to conclude that the solution~$u$ of~\eqref{eqannulus} is radially symmetric relies on~\cite[Theorem~5]{si}, which itself uses the Lipschitz-continuity of~$f$ over the range of $u$. Thus, the same argument can not be applied as such in general in the case where $v$ is just assumed to have no stagnation point in $\Omega_{a,b}$, without the assumption~\eqref{hypstagnation}. Other arguments should then be used to prove that $v$ is circular or to disprove this property in general. We leave this question open for a further work.


\subsection{The case of unbounded annuli $\Omega_{a,\infty}$: proof of Theo\-rems~\ref{th2} and~\ref{th2bis}}\label{sec32}

This section is devoted to the proof of Theorems~\ref{th2} and~\ref{th2bis}. Throughout this section, we fix a positive real number $a$ and we consider a $C^2(\overline{\Omega_{a,\infty}})$ flow $v$ solving~\eqref{1}-\eqref{2} and such that
\beq\label{defeta}
\big\{x\in\overline{\Omega_{a,\infty}}:|v(x)|=0\big\}\subsetneq C_a\ \hbox{ and }\ |v|\ge\eta>0\hbox{ in }\overline{\Omega_{a+1,\infty}}
\eeq
for some positive real number $\eta>0$ (these conditions are fulfilled in both Theorems~\ref{th2} and~\ref{th2bis}). We also assume that either $v(x)\cdot e_r(x)=o(1/|x|)$ as $|x|\to+\infty$ (that is,~\eqref{hypth22}, for Theorem~\ref{th2}) or $\inf_{\Omega_{a,\infty}}|v|>0$ (for Theorem~\ref{th2bis}). This situation fits into the framework of Section~\ref{sec2}, with $\mathcal{C}_1=C_a$, $\omega_1=B_a$, $\Omega=\R^2\!\setminus\!\overline{\omega_1}=\Omega_{a,\infty}$, and $D=\overline{\Omega_{a,\infty}}$. Notice also that the conditions~\eqref{hypv1}-\eqref{hypv2} and~\eqref{stagnation} are fulfilled by assumption. So is~\eqref{hypC1infty}, hence Lemma~\ref{lem22}-(ii) can be applied. Furthermore, either~\eqref{C2} is fulfilled (from~\eqref{hypth22}, for Theorem~\ref{th2}) and Lemmas~\ref{lem21}-(ii) and~\ref{lem24}-(i) can be applied, or the conditions of Lemmas~\ref{lem23}-(iii) and~\ref{lem24}-(ii) are fulfilled. Therefore, from Lemmas~\ref{lem21},~\ref{lem22}-(ii),~\ref{lem23}-(iii) and~\ref{lem24}, the flow has a $C^3(\overline{\Omega_{a,\infty}})$ stream function~$u$, and there exist a real number $c_{in}$ (which can be taken to be $0$ without loss of generality, since $u$ is unique up to additive constants) and $c_{out}=\pm\infty$ (we can assume that~$c_{out}=+\infty$ without loss of generality, even if it means changing $v$ into $-v$ and $u$ into $-u$), together with a $C^1([0,+\infty))$ function $f$ such that
\beq\label{eqannulus2}\left\{\baa{l}
\Delta u+f(u)=0\hbox{ in }\overline{\Omega_{a,\infty}},\vspace{3pt}\\
u>0\hbox{ in }\Omega_{a,\infty},\vspace{3pt}\\
u=0\hbox{ on }C_a,\ \ u(x)\to+\infty\hbox{ as $|x|\to+\infty$}.\eaa\right.
\eeq

\subsubsection{Proof of Theorem~\ref{th2}}

In addition to~\eqref{defeta}, we further assume that
\beq\label{hypver}
v(x)\cdot e_r(x)=o\Big(\frac{1}{|x|}\Big)\ \hbox{ as }|x|\to+\infty.
\eeq
As a consequence,~\eqref{C2} holds, and there is $R_0\ge a+1$ such that
\beq\label{vetheta}
|\nabla u\cdot e_r|=|v\cdot e_\theta|\ge\frac{\eta}{2}\ \hbox{ in }\overline{\Omega_{R_0,\infty}}.
\eeq
Hence, Lemma~\ref{lem23}-(i)-(ii) can be applied and the streamlines $\Xi_x=\xi_x(\R)$ (with $x\in\Omega_{a,\infty}$) surround the origin and are such that
\beq\label{xixinfty}
\min_\R|\xi_x|\to+\infty\ \hbox{ and }\ \max_\R|\xi_x|-\min_\R|\xi_x|\to0\ \hbox{ as }|x|\to+\infty
\eeq
(actually, the second conclusion is stronger than the first one, since $\xi_x(0)=x$). Together with the normalization of $u$ (that is, $u>0$ in $\Omega_{a,\infty}$), properties~\eqref{vetheta}-\eqref{xixinfty} yield the existence of~$R_1\ge R_0\ge a+1$ such that $\nabla u(x)\cdot e_r(x)\ge\eta/2$ for all $|x|\ge R_0$ and $\min_\R|\xi_x|\ge R_0$ for all~$|x|\ge R_1$. For any $x$ with $|x|\ge R_1$, it then follows that, for every $\theta\in\R$, there is a unique~$\varrho_x(\theta)\ge R_0$ such that $(\varrho_x(\theta)\cos\theta,\varrho_x(\theta)\sin\theta)\in\Xi_x$, and moreover
\beq\label{defvarrho}
\Xi_x=\big\{(\varrho_x(\theta)\cos\theta,\varrho_x(\theta)\sin\theta):\theta\in\R\big\}.
\eeq
Notice also that the $2\pi$-periodic function $\varrho_x$ is of class $C^3(\R)$ from the implicit function theorem.\par
If some streamlines were true circles centered at the origin, then~\cite[Theorem~5]{si} would imply that the stream function~$u$ is radially symmetric in the bounded region between $C_a$ and these streamlines. To circumvent the fact that the streamlines are not known to be true circles a priori, we use Lemma~\ref{lem8} below and Proposition~\ref{promoving} to compare the stream function~$u$ with its reflection with respect to some lines approximating any line containing the origin. We then proceed by passing to the limit as the approximation parameter goes to $0$. With Proposition~\ref{promoving}, it then easily follows that $u$ is radially symmetric and that all streamlines are truly circular, thus completing the proof of Theorem~\ref{th2}.\par
To apply this strategy, let us now introduce some additional notations which will be used in this section, as well as in the proof of Theorems~\ref{th2bis},~\ref{th3},~\ref{th4} and~\ref{th5} in the following sections. For $x\in\Omega_{a,\infty}$, let~$\Omega_x$ denote the bounded connected component of $\R^2\!\setminus\!\Xi_x$. Notice that $\Omega_x$ is well defined and contains the origin, by Lemma~\ref{lem23}-(i). Notice also that $u$ is equal to the positive constant~$u(x)$ along $\Xi_x$, while~$u$ vanishes along $C_a$ and has no critical point in $\Omega_{a,\infty}$. Hence,
\beq\label{strict}
0<u(y)<u(x)\ \hbox{ for all }y\in\Omega_x\cap\Omega_{a,\infty},
\eeq
where $\Omega_x\cap\Omega_{a,\infty}$ is the bounded domain located between $\Xi_x$ and $C_a$. As a consequence, $\nabla u(z)$ points outwards $\Omega_x$ at each point~$z\in\Xi_x$.\par
After recalling that the sets $T_{e,\lambda}$ and $H_{e,\lambda}$ and the reflection $R_{e,\lambda}$ have been defined in~\eqref{Telambda}-\eqref{Relambda}, the following lemma says that, for any $\epsilon>0$, the set $\Omega_x\cap H_{e,\lambda}$ will be an admissible set for the method of moving planes for any $e\in\mathbb{S}^1$ and $\lambda>\epsilon>0$, provided $|x|$ is large enough. 

\begin{figure}
\centering\includegraphics[scale=0.7]{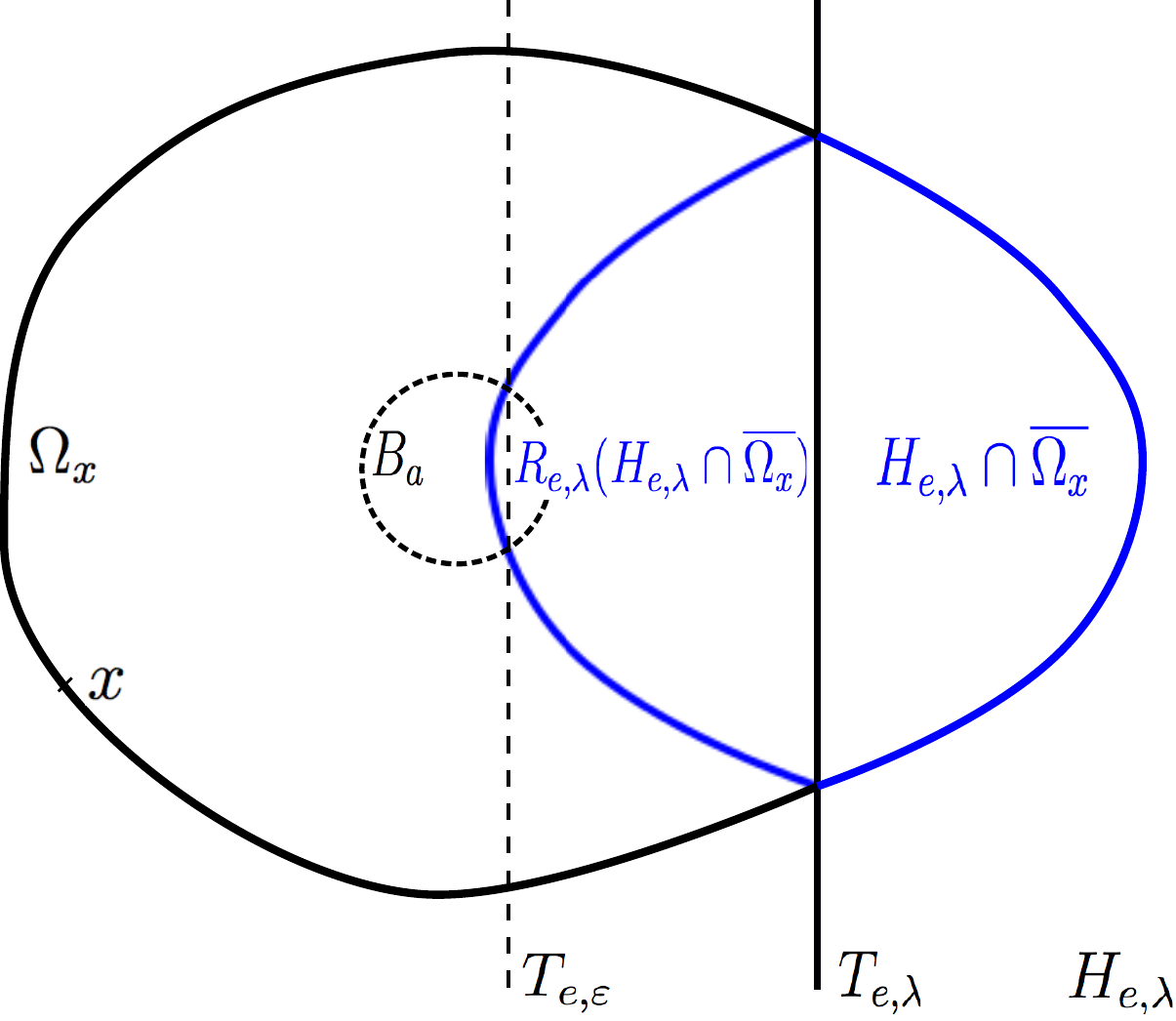}
\caption{The sets $\Omega_x$, $H_{e,\lambda}\cap\overline{\Omega_x}$ and $R_{e,\lambda}(H_{e,\lambda}\cap\overline{\Omega_x})$}
\end{figure}

\begin{lem}\label{lem8}
For each $\epsilon>0$, there exists $R_\epsilon>a$ such that
$$R_{e,\lambda}\big(H_{e,\lambda}\cap\overline{\Omega_x}\big)\subset\Omega_x$$
for all $e\in\mathbb{S}^1$, $\lambda>\epsilon$ and $|x|\ge R_\epsilon$ $($see Fig.~$2)$.
\end{lem}

\noindent{\bf{Proof.}} Fix $\epsilon>0$, and assume by way of contradiction that the conclusion of the lemma does not hold. Then there are some sequences $(x_n)_{n\in\N}$ in $\Omega_{a,\infty}$, $(e_n)_{n\in\N}$ in $\mathbb{S}^1$, $(\lambda_n)_{n\in\N}$ in $(\epsilon,+\infty)$ and $(y_n)_{n\in\N}$ such that
$$\lim_{n\to+\infty}|x_n|=+\infty,\  \hbox{ and }\ y_n\in H_{e_n,\lambda_n}\cap\overline{\Omega_{x_n}}\hbox{ and }z_n:=R_{e_n,\lambda_n}(y_n)\not\in\Omega_{x_n}\hbox{ for all }n\in\N.$$
From~\eqref{xixinfty}, there is a sequence $(r_n)_{n\in\N}$ of positive real numbers converging to $0$ such that~$B_{|x_n|-r_n}\subset\Omega_{x_n}\subset B_{|x_n|+r_n}$ for all $n\in\N$, hence $|y_n|\le|x_n|+r_n$.  On the other hand, since $y_n\cdot e_n>\lambda_n>\epsilon>0$, one has
$$|y_n|^2-|z_n|^2=|y_n|^2-|R_{e_n,\lambda_n}(y_n)|^2=4\lambda_n(y_n\cdot e_n-\lambda_n)>0,$$
hence $|y_n|>|z_n|\ge|x_n|-r_n$ since $z_n\not\in\Omega_{x_n}$. As a consequence, $|x_n|-r_n\le|z_n|<|y_n|\le|x_n|+r_n$ for all $n\in\N$, and $\lim_{n\to+\infty}(|y_n|-|x_n|)=\lim_{n\to+\infty}(|y_n|-|z_n|)=0$. The inequality $|y_n|^2-|z_n|^2=4\lambda_n(y_n\cdot e_n-\lambda_n)>4\epsilon(y_n\cdot e_n-\lambda_n)>0$ then yields $\lim_{n\to+\infty}(y_n\cdot e_n-\lambda_n)=0$. Hence,
$${\rm{dist}}(y_n,\Xi_{x_n}\cap T_{e_n,\lambda_n})\to0\ \hbox{ and }\ |y_n-z_n|\to0\ \hbox{ as }n\to+\infty.$$\par
For each $n\in\N$, let $\varphi_n\in\R$ be such that $e_n=(\cos\varphi_n,\sin\varphi_n)$. Since $y_n\cdot e_n>\lambda_n>\epsilon>0$, there is a unique $\theta_n\in(-\pi/2,\pi/2)$ such that
$$\frac{y_n}{|y_n|}=(\cos(\varphi_n+\theta_n),\sin(\varphi_n+\theta_n)).$$
Similarly, since $(z_n-y_n)\cdot e_n\to0$ as $n\to+\infty$, one has $z_n\cdot e_n>\epsilon/2$ for all large $n$ and there is a unique $\theta'_n\in(-\pi/2,\pi/2)$ such that
$$\frac{z_n}{|z_n|}=(\cos(\varphi_n+\theta'_n),\sin(\varphi_n+\theta'_n)).$$
Since $\lim_{n\to+\infty}|y_n-z_n|=0$ and $\lim_{n\to+\infty}|y_n|=\lim_{n\to+\infty}|z_n|=\lim_{n\to+\infty}|x_n|=+\infty$, one also infers that $\theta_n-\theta'_n\to0$ as $n\to+\infty$. We also recall that~\eqref{defvarrho} holds with $x_n$ instead of $x$, for all $n$ large enough. It then follows from Lemma~\ref{lem23}-(i) and from the assumptions on $y_n$ and~$z_n$ that $|y_n|\le\varrho_{x_n}(\varphi_n+\theta_n)$ and $|z_n|\ge\varrho_{x_n}(\varphi_n+\theta'_n)$ for all $n$ large enough. Denote, for $n$ large enough,
$$\left\{\baa{l}
y'_n=(\varrho_{x_n}(\varphi_n+\theta_n)\cos(\varphi_n+\theta_n),\varrho_{x_n}(\varphi_n+\theta_n)\sin(\varphi_n+\theta_n))\in\Xi_{x_n},\vspace{3pt}\\
z'_n=(\varrho_{x_n}(\varphi_n+\theta'_n)\cos(\varphi_n+\theta'_n),\varrho_{x_n}(\varphi_n+\theta'_n)\sin(\varphi_n+\theta'_n))\in\Xi_{x_n},\eaa\right.$$
and observe that $y_n\in(0,y'_n]$ and $z'_n\in(0,z_n]$.\par
We now claim that $\theta'_n\neq\theta_n$ for all $n$ large enough. Indeed, otherwise, up to extraction of a subsequence, $y'_n=z'_n$ and the four points $0$, $y_n$, $y'_n=z'_n$ and $z_n$ would be aligned in that order. But since $y_n-z_n=2(y_n\cdot e_n-\lambda_n)e_n$ with $y_n\cdot e_n-\lambda_n>0$, the vectors $y_n$ and $z_n$ would be parallel to $e_n$. Hence, $y_n=(y_n\cdot e_n)e_n$ with $y_n\cdot e_n>\lambda_n>\epsilon>0$ and $z_n=(z_n\cdot e_n)e_n$ with~$z_n\cdot e_n=2\lambda_n-y_n\cdot e_n<\lambda_n<y_n\cdot e_n$. This contradicts the fact that $0$, $y_n$ and $z_n$ lie on the half-line~$\R_+e_n$ in that order. Thus, $\theta'_n\neq\theta_n$ for all $n$ large enough, thus for all $n$ without loss of generality. Notice that the same arguments also imply that $\theta_n\neq0$ and $\theta'_n\neq0$ for all~$n$ large enough (since otherwise in either case one would have $\theta_n=\theta'_n=0$ up to extraction of a subsequence), thus for all $n$ without loss of generality. In particular, either $0<\theta_n<\pi/2$ or~$-\pi/2<\theta_n<0$.\par
\begin{figure}
\centering\includegraphics[scale=0.7]{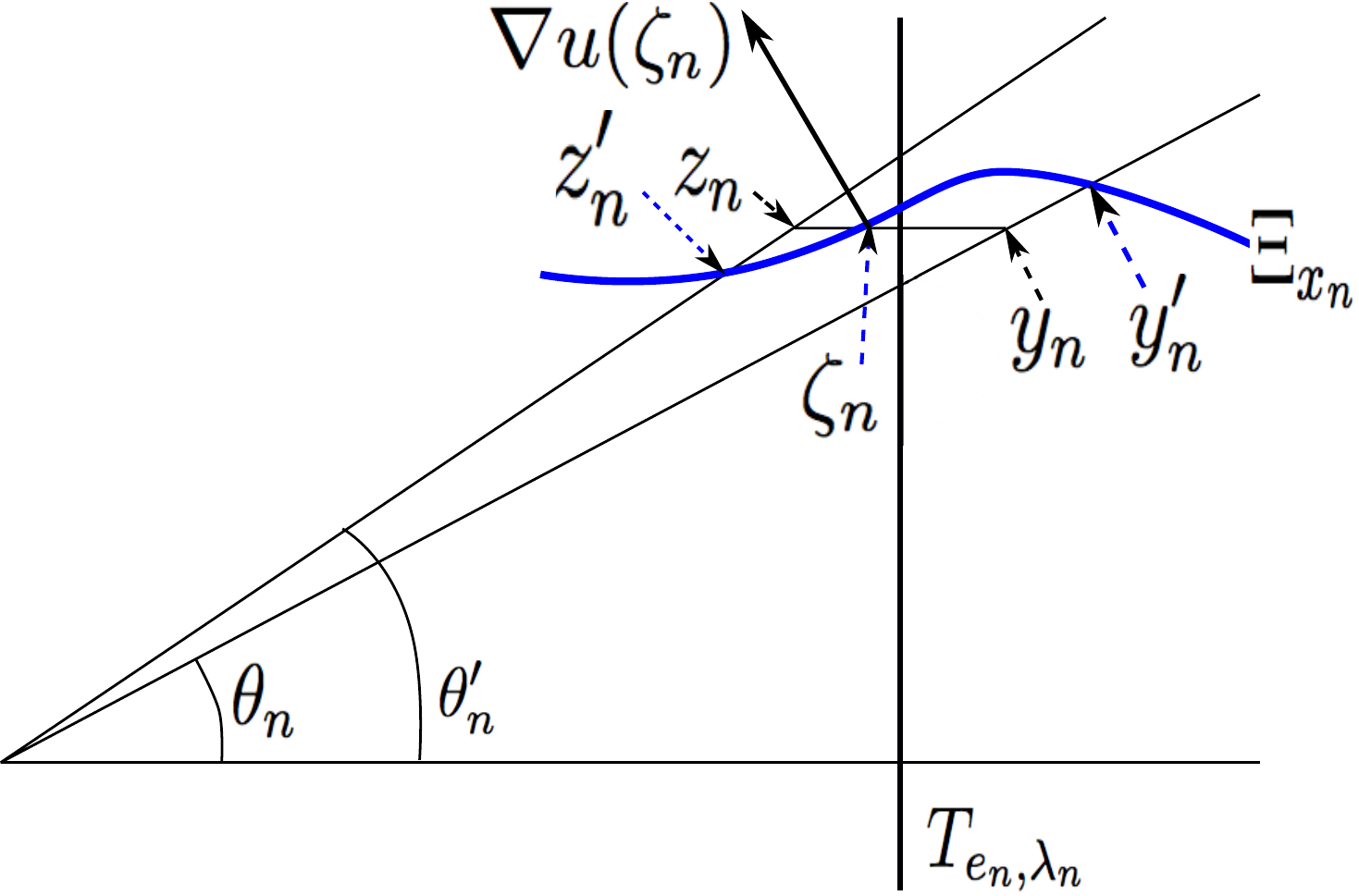}
\caption{The points $y_n$, $y'_n$, $z_n$, $z'_n$, and $\zeta_n$ (with here $e_n=(1,0)$ and $\varphi_n=0$)}
\end{figure}
Assume first that, up to extraction of a subsequence, $0<\theta_n<\pi/2$ for all $n$. One then infers from the definition of $z_n=R_{e_n\lambda_n}(y_n)$ and the previous paragraph that
$$0<\theta_n<\theta'_n<\frac{\pi}{2}.$$
Remember now that $0<u(y)<u(x_n)$ for every $y$ in the domain $\Omega_{x_n}\cap\Omega_{a,\infty}$ between $\Xi_{x_n}$ and~$C_a$, and $\nabla u(z)$ points outwards $\Omega_{x_n}$ at each point $z\in\Xi_{x_n}$. For each $n\in\N$, since $y_n\in(0,y'_n]$, $z'_n\in(0,z_n]$ and since $y_n-z_n=\varsigma_ne_n$ with $\varsigma_n:=2(y_n\cdot e_n-\lambda_n)>0$, there is then an angle $\phi_n\in[\theta_n,\theta'_n]\subset(0,\pi/2)$ such that
$$\zeta_n:=(\varrho_{x_n}(\varphi_n+\phi_n)\cos(\varphi_n+\phi_n),\varrho_{x_n}(\varphi_n+\phi_n)\sin(\varphi_n+\phi_n))\in\Xi_{x_n}\cap[y_n,z_n]$$
and $\nabla u(\zeta_n)\cdot e_n\le 0$, see Fig.~3. The point $\zeta_n$ can be defined as the first point on~$\Xi_{x_n}\cap[y_n,z_n]$ when going from~$y'_n$ to~$z'_n$ along~$\Xi_{x_n}$ with increasing angle. Notice that $|\zeta_n|\to+\infty$ since $|y_n|\to+\infty$ and $|y_n-z_n|\to0$ as $n\to+\infty$. Call $v_{1,n}=v(\zeta_n)\cdot e_n$ and $v_{2,n}=v(\zeta_n)\cdot e_n^\perp$. The inequality $\nabla u(\zeta_n)\cdot e_n\le 0$ means that~$v_{2,n}\le0$. Therefore,
\beq\label{zetan}
v(\zeta_n)\cdot e_r(\zeta_n)=v_{1,n}\cos\phi_n+v_{2,n}\sin\phi_n\le v_{1,n}\cos\phi_n,
\eeq
while $0<\eta/2\le|v(\zeta_n)\cdot e_\theta(\zeta_n)|=|-v_{1,n}\sin\phi_n+v_{2,n}\cos\phi_n|$ for all $n$ large enough, from~\eqref{vetheta} and $\lim_{n\to+\infty}|\zeta_n|=+\infty$. But since the continuous function $v\cdot e_\theta=\nabla u\cdot e_r$ has a constant strict positive sign at infinity, it follows that
$$\frac{\eta}{2}\le v(\zeta_n)\cdot e_\theta(\zeta_n)=-v_{1,n}\sin\phi_n+v_{2,n}\cos\phi_n$$
for all $n$ large enough. Since $v_{2,n}\le0$ and $0<\phi_n<\pi/2$, one gets that $-v_{1,n}\sin\phi_n\ge\eta/2$, hence~$v_{1,n}\le-\eta/2$. Together with~\eqref{zetan}, it follows that $v(\zeta_n)\cdot e_r(\zeta_n)\!\le\!-(\eta/2)\cos\phi_n$ for all~$n$ large enough. On the other hand, since $\zeta_n\!\in\![y_n,z_n]$ and $\lim_{n\to+\infty}|z_n\!-\!y_n|\!=\!\lim_{n\to+\infty}(y_n\cdot e_n\!-\!\lambda_n)\!=\!0$, there holds~$\zeta_n\cdot e_n-\lambda_n\to0$, hence $\zeta_n\cdot e_n\ge\epsilon/2$ for all $n$ large enough (since~$\lambda_n>\epsilon>0$ for all $n$). Finally,
$$\cos\phi_n=\frac{\zeta_n\cdot e_n}{|\zeta_n|}\ge\frac{\epsilon}{2|\zeta_n|}\ \hbox{ and }\ v(\zeta_n)\cdot e_r(\zeta_n)\le-\frac{\eta\,\epsilon}{4|\zeta_n|}$$
for all $n$ large enough. That last inequality contradicts the assumption~\eqref{hypver} and the limit $\lim_{n\to+\infty}|\zeta_n|=+\infty$.\par
The second case, for which, up to extraction of a subsequence, $-\pi/2<\theta_n<0$ for all $n$ (and then $-\pi/2<\theta'_n<\theta_n<0$) can be handled similarly and leads to a contradiction as well. The proof of Lemma~\ref{lem8} is thereby complete.\hfill$\Box$\break

\noindent{\bf{Proof of Theorem~\ref{th2}.}} We shall show that the stream function $u$ is radially symmetric in~$\overline{\Omega_{a,\infty}}$. Notice that we already know that $u=0$ on $C_a$. Let then $x\neq y\in\Omega_{a,\infty}$ be such that
$$|x|=|y|\ (>a).$$
Call
\beq\label{defe}
e=\frac{y-x}{|y-x|}\in\mathbb{S}^1.
\eeq
Consider an arbitrary real number $\epsilon$ such that $0<\epsilon<a$. Let $R_\epsilon>a$ be as in Lemma~\ref{lem8}. From~\eqref{xixinfty}, there is a point $x_\epsilon\in\Omega_{a,\infty}$ such that $|x_\epsilon|\ge R_\epsilon$ and $\min_\R|\xi_{x_\epsilon}|>|x|=|y|$. Lemma~\ref{lem8} then yields
\beq\label{reflection}
R_{e,\lambda}(H_{e,\lambda}\cap\overline{\Omega_{x_\epsilon}})\subset\Omega_{x_\epsilon}\ \hbox{ for all }\lambda>\epsilon.
\eeq\par
We are now going to apply Proposition~\ref{promoving} with
$$\left\{\baa{l}
\Xi=\Xi_{x_\epsilon},\ \Omega=\Omega_{x_\epsilon},\ \Xi'=C_a,\ \Omega'=B_a,\ \omega=\Omega_{x_\epsilon}\!\setminus\!\overline{B_a},\vspace{3pt}\\
\displaystyle R'=a,\ R=\max_\R|\xi_{x_\epsilon}|>a,\ \overline{\lambda}=\max_{z\in\Xi_{x_\epsilon}}z\cdot e>a>\epsilon>0,\vspace{3pt}\\
\varphi=-u\in C^3(\overline{\omega}),\ c_1=-u(x_\epsilon)=-u_{|\Xi_{x_\epsilon}}<0,\ c_2=0=-u_{|C_a},\vspace{3pt}\\
F(r,s)=F(s)=-f(-s)\hbox{ for }(r,s)\in[a,R]\times[-u(x_\epsilon),0].\eaa\right.$$
Notice immediately that assumption~\eqref{hypelambda2} is automatically satisfied. The function~$F$ clearly satisfies the assumptions of Proposition~\ref{promoving} since $f$ is of class $C^1([0,+\infty))$. The function~$\varphi$ satisfies $\Delta\varphi+F(\varphi)=0$ in~$\overline{\omega}$, with $c_1<\varphi<c_2$ in $\omega$ (since $0<u<u(x_\epsilon)$ in $\Omega_{x_\epsilon}\cap\Omega_{a,\infty}$ by~\eqref{strict}). Together with~\eqref{reflection}, all assumptions of Proposition~\ref{promoving} are satisfied. Proposition~\ref{promoving} applied with $\lambda=\epsilon$ then implies that $\varphi\le\varphi_{e,\epsilon}$, namely $u\ge u_{e,\epsilon}$, in $\overline{\omega_{e,\epsilon}}$ with
$$\omega_{e,\epsilon}=(H_{e,\epsilon}\cap\omega)\!\setminus\! R_{e,\epsilon}(\overline{\Omega'})=\big(H_{e,\epsilon}\cap(\Omega_{x_\epsilon}\!\setminus\!\overline{B_a})\big)\!\setminus\! R_{e,\epsilon}(\overline{B_a}).$$
Observe now that $y\cdot e=(|y|^2-x\cdot y)/|y-x|>0$ since $|x|=|y|$ and $x\neq y$, and remember that $a<|y|<\min_\R|\xi_{x_\epsilon}|$, hence $y\in\omega$. Therefore, $y\in\omega_{e,\epsilon}$ for all $\epsilon>0$ small enough, and
$$u(y)\ge u_{e,\epsilon}(y)=u(y_{e,\epsilon})=u(y-2(y\cdot e-\epsilon)e)$$
for all $\epsilon>0$ small enough. By passing to the limit $\epsilon\displaystyle\mathop{\to}^>0$ and using the definition of $e$ and the assumption $|x|=|y|$, one infers that
$$u(y)\ge u(y-2(y\cdot e)e)=u(x).$$
Since the last inequality holds for all $x\neq y\in\Omega_{a,\infty}$ such that $|x|=|y|$ (and also for all $x,\,y\in C_a$), the $C^3(\overline{\Omega_{a,\infty}})$ function $u$ is radially symmetric in $\overline{\Omega_{a,\infty}}$. Together with~\eqref{defeta}-\eqref{eqannulus2}, there is then a $C^3([a,+\infty))$ function $U$ such that $U(a)=0$, $U'>0$ in $[a,+\infty)$ and $u(x)=U(|x|)$ for all $x\in\overline{\Omega_{a,\infty}}$. This means that $v(x)=V(|x|)\,e_\theta(x)$ for all $x\in\overline{\Omega_{a,\infty}}$ with $V=U'\in C^2([a,+\infty))$ and $V>0$ in $[a,+\infty)$. The proof of Theorem~\ref{th2} is thereby complete.~\hfill$\Box$

\subsubsection{Proof of Theorem~\ref{th2bis}}

This section is devoted to the proof of Theorem~\ref{th2bis}. In addition of~\eqref{defeta}, we actually assume the stronger condition
\beq\label{defeta3}
|v|=|\nabla u|\ge\eta>0\ \hbox{ in }\overline{\Omega_{a,\infty}}
\eeq
for some $\eta>0$. Since from our normalization the function $u$ is positive in $\Omega_{a,\infty}$ and vanishes on $C_a$, the conditions~\eqref{2} and~\eqref{defeta3} imply that
\beq\label{vetheta2}
v\cdot e_\theta=\nabla u\cdot e_r\ge\eta>0\ \hbox{ on }C_a.
\eeq\par
To prove Theorem~\ref{th2bis}, we then have to show that the supremum of the vorticity is positive, namely
\beq\label{vorticity}
\sup_{\Omega_{a,\infty}}\!\Big(\frac{\partial v_2}{\partial x_1}-\frac{\partial v_1}{\partial x_2}\Big)>0.
\eeq
To do so, let us assume by way of contradiction that
\beq\label{assumvort}
\Delta u=\frac{\partial v_2}{\partial x_1}-\frac{\partial v_1}{\partial x_2}\le0\ \hbox{ in }\overline{\Omega_{a,\infty}}.
\eeq
We will show that $u$ is radially symmetric, and this will easily lead to a contradiction. To prove the radially symmetry of $u$, let us use a Kelvin transform of the variables by setting
$$w(x)=u\Big(\frac{x}{|x|^2}\Big)\ \hbox{ for }x\in\overline{\Omega_{0,1/a}}\!\setminus\!\{0\},$$
and let us show that the $C^3(\overline{\Omega_{0,1/a}}\!\setminus\!\{0\})$ function $w$ is radially symmetric in $\overline{\Omega_{0,1/a}}\!\setminus\!\{0\}$. From~\eqref{eqannulus2}, one has
$$w=0\hbox{ on }C_{1/a},\ \ w>0\hbox{ in }\Omega_{0,1/a}\ \hbox{ and }\ w(x)\to+\infty\hbox{ as }|x|\mathop{\to}^>0,$$
and a straightforward calculation yields
$$\Delta w(x)+\frac{1}{|x|^4}\,f(w(x))=0\ \hbox{ for all }x\in\overline{\Omega_{0,1/a}}\!\setminus\!\{0\},$$
that is, $\Delta w(x)+F(|x|,w(x))=0$ in $\overline{\Omega_{0,1/a}}\!\setminus\!\{0\}$ with
$$\baa{rcl}
F:(0,1/a]\times[0,+\infty) & \to & \R\vspace{3pt}\\
(r,s) & \mapsto & F(r,s)=r^{-4}f(s).\eaa$$
The function $F$ is of class $C^1((0,1/a]\times[0,+\infty))$. Furthermore, the range of $u$ is equal to the whole interval $[0,+\infty)$ by~\eqref{eqannulus2}, and $f\ge0$ in $[0,+\infty)$ by~\eqref{eqannulus2} and~\eqref{assumvort}. Therefore, the function $F$ is nonincreasing with respect to its first variable in~$(0,1/a]\times[0,+\infty)$.\par
Consider now any two points $x\neq y\in\overline{\Omega_{0,1/a}}\!\setminus\!\{0\}$ with $|x|=|y|$. As in the proof of Theorem~\ref{th2}, denote
$$e=\frac{y-x}{|y-x|}\in\mathbb{S}^1$$
and consider an arbitrary real number $\epsilon$ such that
$$0<\epsilon<|x|=|y|\le\frac{1}{a}.$$
By Lemma~\ref{lem23}-(iii), there is a point $x_\epsilon\in\Omega_{a,\infty}$ such that $\min_\R|\xi_{x_\epsilon}|>1/\epsilon>a$. One knows that the streamline $\Xi_{x_\epsilon}$ surrounds the origin and that $u=u(x_\epsilon)>0$ along $\Xi_{x_\epsilon}$. Furthermore, as in~\eqref{strict}, one has $0<u<u(x_\epsilon)$ in the domain $\Omega_{x_\epsilon}\cap\Omega_{a,\infty}$ between~$\Xi_{x_\epsilon}$ and $C_a$, since $u=0$ on $C_a$ and $u$ has no critical point in $\Omega_{a,\infty}$.\par
Denote $\Xi=C_{1/a}$ and
$$\Xi'=\Big\{z\in\R^2:\frac{z}{|z|^2}\in\Xi_{x_\epsilon}\Big\}.$$
Notice that the Jordan curve $\Xi'$ surrounds the origin and $\Xi'\subset B_\epsilon\,(\subset B_{1/a})$ by definition of $x_\epsilon$. Call $\Omega=B_{1/a}$, let $\Omega'$ be the bounded connected component of $\R^2\!\setminus\!\Xi'$, and let
$$\omega=\Omega\!\setminus\!\overline{\Omega'}=B_{1/a}\!\setminus\!\overline{\Omega'}\ (\supset\Omega_{\epsilon,1/a}).$$
Denote $R=1/a$,
$$0<R'=\min_{z\in\Xi'}|z|=\frac{1}{\displaystyle\max_\R|\xi_{x_\epsilon}|}<\epsilon<R,$$
and $\overline{\lambda}=\max_{z\in\Xi}z\cdot e=1/a>0$. One has $0<\epsilon<1/a$, hence $\epsilon\in[0,\overline{\lambda})$. The function $\varphi=w$ is of class $C^3(\overline{\omega})$ with
$$\hbox{$\varphi=c_1=0$ on $\Xi=C_{1/a}$, $\varphi=c_2=u(x_\epsilon)>0$ on $\Xi'$, and $0<\varphi<u(x_\epsilon)$ in $\omega$}$$
(since $0<u<u(x_\epsilon)$ in $\Omega_{x_\epsilon}\cap\Omega_{a,\infty}$). Furthermore, $\varphi$ satisfies $\Delta\varphi+F(|x|,\varphi)=0$ in $\overline{\omega}$, with~$F(r,s)=r^{-4}f(s)$ and, here, $(r,s)\in[R',1/a]\times[0,u(x_\epsilon)]$. The function $F$ then satisfies the conditions of Proposition~\ref{promoving}. Lastly, the condition~\eqref{hypelambda} is immediately satisfied since $\Omega=B_{1/a}$ and the condition~\eqref{hypelambda2} also holds since $H_{e,\lambda}\cap\Xi'=\emptyset$ for all $\lambda>\epsilon$ (because $\Xi'\subset B_\epsilon$). To sum up, all assumptions of Proposition~\ref{promoving} are fulfilled. Its conclusion with $\lambda=\epsilon$ yields $w\le w_{e,\epsilon}$ in $\overline{\omega_{e,\epsilon}}$, with
$$\omega_{e,\epsilon}=\big(H_{e,\epsilon}\cap(B_{1/a}\!\setminus\!\overline{\Omega'})\big)\!\setminus\! R_{e,\epsilon}(\overline{\Omega'}).$$
Since $y\cdot e>0$ and since $\overline{\Omega'}\subset B_{\epsilon}$ and $R_{e,\epsilon}(\overline{\Omega'})\subset B_{3\epsilon}$, it follows that $y\in\omega_{e,\epsilon}$ for all $\epsilon>0$ small enough. As a consequence, $w(y)\le w_{e,\epsilon}(y)=w(y_{e,\epsilon})=w(y-2(y\cdot e-\epsilon)e)$ for all $\epsilon>0$ small enough and the passage to the limit as $\epsilon\to0$ yields
$$w(y)\le w(y-2(y\cdot e)e)=w(x)$$
by definition of $e$. Since this holds for all $x\neq y\in\overline{\Omega_{0,1/a}}\!\setminus\!\{0\}$ with $|x|=|y|$, this means that $w$ is radially symmetric in $\overline{\Omega_{0,1/a}}\!\setminus\!\{0\}$, hence $u$ is radially symmetric in $\overline{\Omega_{a,\infty}}$. Together with~\eqref{vetheta2}, there is then a $C^3([a,+\infty))$ function $U$ such that $u(x)=U(|x|)$ and $U'\ge\eta>0$ in $[a,+\infty)$. But $\Delta u\le0$ in $\overline{\Omega_{a,\infty}}$ by~\eqref{assumvort}. Hence $U''(r)+r^{-1}U'(r)\le0$ in $[a,+\infty)$ and the function $r\mapsto rU'(r)$ is nonincreasing in $[a,+\infty)$, a contradiction with $U'\ge\eta>0$.\par
As a conclusion,~\eqref{assumvort} can not hold, that is,~\eqref{vorticity} holds and the proof of Theorem~\ref{th2bis} is thereby complete.\hfill$\Box$


\subsection{The case of punctured disks $\Omega_{0,b}$: proof of Theorem~\ref{th3}}\label{sec33}

This section is devoted to the proof of Theorem~\ref{th3}. Throughout this section, we fix a positive real number $b$ and we consider a $C^2(\overline{\Omega_{0,b}}\!\setminus\!\{0\})$ flow $v$ solving~\eqref{1}-\eqref{2} and such that
\beq\label{hypth3bis}
\big\{x\in\overline{\Omega_{0,b}}\!\setminus\!\{0\}:|v(x)|=0\big\}\subsetneq C_b\ \hbox{ and }\ \int_{C_\epsilon}|v\cdot e_r|\to0\hbox{ as }\epsilon\displaystyle\mathop{\to}^>0.
\eeq
This situation falls within the framework of Section~\ref{sec2}, with $\mathcal{C}_2=C_b$, $\omega_2=B_b$, $\Omega\!=\!\omega_2\!\setminus\!\{0\}\!=\!\Omega_{0,b}$, and $D=\overline{\Omega_{0,b}}\!\setminus\!\{0\}$, and the conditions~\eqref{hypv1}-\eqref{hypv2},~\eqref{C2}, and~\eqref{stagnation} are fulfilled by assumption. Therefore, from Lemmas~\ref{lem21} and~\ref{lem24}-(i), the flow has a $C^3(\overline{\Omega_{0,b}}\!\setminus\!\{0\})$ stream function~$u$, and there exist a real number $c_{out}$ (which can be taken to be $0$ without loss of generality, since $u$ is unique up to additive constants) and $0\neq c_{in}\in\overline{\R}$ (we can assume that~$0<c_{in}\le+\infty$ without loss of generality, even if it means changing $v$ into $-v$ and~$u$ into~$-u$), together with a $C^1([0,c_{in}))$ function $f$ such that
\beq\label{eqannulus3}\left\{\baa{l}
\Delta u+f(u)=0\hbox{ in }\overline{\Omega_{0,b}}\!\setminus\!\{0\},\vspace{3pt}\\
u>0\hbox{ in }\Omega_{0,b},\vspace{3pt}\\
u=0\hbox{ on }C_b,\ \ u(x)\to c_{in}\hbox{ as }|x|\displaystyle\mathop{\to}^>0.\eaa\right.
\eeq\par
One goal is to show that $u$ is radially symmetric in $\overline{\Omega_{0,b}}\!\setminus\!\{0\}$. Consider any two points $x\neq y\in\overline{\Omega_{0,b}}\!\setminus\!\{0\}$ with $|x|=|y|$. As in the proof of Theorems~\ref{th2} and~\ref{th2bis}, denote
$$e=\frac{y-x}{|y-x|}\in\mathbb{S}^1$$
and consider an arbitrary real number $\epsilon$ such that
$$0<\epsilon<|x|=|y|\le b.$$
By Lemma~\ref{lem23}-(i), there is a point $x_\epsilon\in\Omega_{0,b}$ such that $\max_\R|\xi_{x_\epsilon}|<\epsilon<b$, and the streamline~$\Xi_{x_\epsilon}$ surrounds the origin with $u=u(x_\epsilon)>0$ along $\Xi_{x_\epsilon}$. Since $u=0$ on $C_b$ and~$u$ has no critical point in $\Omega_{0,b}$, one infers that $0<u<u(x_\epsilon)$ in the domain $\omega$ between $\Xi_{x_\epsilon}$ and~$C_b$. Denote $\Xi=C_b$, $\Xi'=\Xi_{x_\epsilon}$, $\Omega=B_b$, let $\Omega'=\Omega_{x_\epsilon}$ be the bounded connected component of~$\R^2\!\setminus\!\Xi_{x_\epsilon}$, and notice that
$$\omega=\Omega\!\setminus\!\overline{\Omega'}.$$
Set
$$R=b,\ \ 0<R'=\min_{z\in\Xi'}|z|<\epsilon<R\ \hbox{ and }\ \overline{\lambda}=\max_{z\in\Xi}z\cdot e=b>0.$$
One has $\epsilon\in(0,\overline{\lambda})$. The function $\varphi=u$ is of class $C^3(\overline{\omega})$ with
$$\varphi=c_1=0\hbox{ on }\Xi=C_b,\ \ \varphi=c_2=u(x_\epsilon)>0\hbox{ on }\Xi',\ \hbox{ and }\ 0<\varphi<u(x_\epsilon)\hbox{ in }\omega.$$
Furthermore, $\varphi$ satisfies $\Delta\varphi+F(\varphi)=0$ in $\overline{\omega}$, with $F:[R',b]\times[0,u(x_\epsilon)]\ni(r,s)\mapsto F(r,s)=f(s)$ satisfying all conditions of Proposition~\ref{promoving}. Lastly, the condition~\eqref{hypelambda} is immediately satisfied since $\Omega=B_b$ and the condition~\eqref{hypelambda2} also holds since $H_{e,\lambda}\cap\Xi'=\emptyset$ for all $\lambda>\epsilon$ (because $\Xi'\subset B_\epsilon$). To sum up, all assumptions of Proposition~\ref{promoving} are fulfilled. Its conclusion with $\lambda=\epsilon$ yields $u\le u_{e,\epsilon}$ in $\overline{\omega_{e,\epsilon}}$, with
$$\omega_{e,\epsilon}=\big(H_{e,\epsilon}\cap(B_b\!\setminus\!\overline{\Omega'})\big)\setminus R_{e,\epsilon}(\overline{\Omega'}).$$
Since $y\cdot e>0$ and since $\overline{\Omega'}\subset B_{\epsilon}$ and $R_{e,\epsilon}(\overline{\Omega'})\subset B_{3\epsilon}$, it follows that $y\in\omega_{e,\epsilon}$ for all $\epsilon>0$ small enough. As a consequence, $u(y)\le u(y_{e,\epsilon})=u(y-2(y\cdot e-\epsilon)e)$ for all $\epsilon>0$ small enough and the passage to the limit as $\epsilon\displaystyle\mathop{\to}^>0$ yields
$$u(y)\le u(y-2(y\cdot e)e)=u(x)$$
by definition of $e$. Since this holds for all $x\neq y\in\overline{\Omega_{0,b}}\!\setminus\!\{0\}$ with $|x|=|y|$, this means that $u$ is radially symmetric in $\overline{\Omega_{0,b}}\!\setminus\!\{0\}$. Together with~\eqref{hypth3bis}-\eqref{eqannulus3}, there is then a~$C^3((0,b])$ function~$U$ such that $u(x)=U(|x|)$ and $U'<0$ in $(0,b]$. Hence, $v(x)=U'(|x|)e_\theta(x)$ for all~$x\in\overline{\Omega_{0,b}}\!\setminus\!\{0\}$, and the proof of Theorem~\ref{th3} is thereby complete.\hfill$\Box$


\subsection{The case of the punctured plane $\Omega_{0,\infty}=\R^2\!\setminus\!\{0\}$: proof of Theorem~\ref{th4}}\label{sec34}

This section is devoted to the proof of Theorem~\ref{th4}. Let $v$ be a $C^2(\Omega_{0,\infty})$ flow solving~\eqref{1} and such that~$|v|>0$ in $\Omega_{0,\infty}$, $\liminf_{|x|\to+\infty}|v(x)|>0$, and~\eqref{hypth4} holds. This situation fits into the framework of Section~\ref{sec2}, with $\Omega=D=\R^2\!\setminus\!\{0\}=\Omega_{0,\infty}$, and the conditions~\eqref{hypv1}-\eqref{hypv2},~\eqref{C2}, and~\eqref{stagnation} are fulfilled by assumption. Therefore, from Lemmas~\ref{lem21} and~\ref{lem24}-(i), the flow has a $C^3(\Omega_{0,\infty})$ stream function~$u$, and there exist $c_{in}\neq c_{out}\in\overline{\R}$ (we can assume that~$c_{in}>c_{out}$ without loss of generality, even if it means changing $v$ into $-v$ and~$u$ into~$-u$), together with a $C^1((c_{out},c_{in}))$ function $f$ such that
\beq\label{eqannulus4}\left\{\baa{l}
\Delta u+f(u)=0\hbox{ in }\Omega_{0,\infty},\vspace{3pt}\\
c_{out}<u<c_{in}\hbox{ in }\Omega_{0,\infty},\vspace{3pt}\\
u(x)\to c_{in}\hbox{ as }|x|\displaystyle\mathop{\to}^>0,\ \ u(x)\to c_{out}\hbox{ as }|x|\to+\infty.\eaa\right.
\eeq
Moreover, the conditions~\eqref{hypth4} and $\liminf_{|x|\to+\infty}|v(x)|\!>\!0$ yield $\liminf_{|x|\to+\infty}|v(x)\cdot e_\theta(x)|\!>\!0$. It then follows from Lemma~\ref{lem23}-(i)-(ii) that each streamline $\Xi_x=\xi_x(\R)$ is a $C^1$ Jordan curve surrounding the origin, with
$$\max_\R|\xi_x|-\min_\R|\xi_x|\to0\hbox{ as }|x|\to+\infty\ \hbox{ and }\ \max_\R|\xi_x|\to0\hbox{ as~$|x|\displaystyle\mathop{\to}^>0$}.$$
Since $u$ is unique up to additive constants, one can also assume without loss of generality that $c_{in}>0>c_{out}$. Pick then a point $X\in\Omega_{0,\infty}$ such that $u(X)=0$, and let $\Omega_X$ be the bounded connected component of $\R^2\!\setminus\!\Xi_X$ (then $0\in\Omega_X$). Together with the inequalities $c_{in}>0>c_{out}$ and the assumption $|\nabla u|>0$ in $\Omega_{0,\infty}$, it follows that $u>0$ in $\Omega_X\!\setminus\!\{0\}$ and $u<0$ in $\R^2\!\setminus\!\overline{\Omega_X}$. Furthermore, by definition of $u$, $v$ is orthogonal to the normal vector to $\Omega_X$. As a consequence, Lemma~\ref{lem23}-(iii) can also be applied with $\omega_1=\Omega_X$, hence $c_{out}=-\infty$ in~\eqref{eqannulus4}.\par
Now, still using the notations~\eqref{Telambda}-\eqref{Relambda}, it then follows as in Lemma~\ref{lem8}, from the assumptions $\liminf_{|x|\to+\infty}|v(x)|>0$ and $v(x)\cdot e_r(x)=o(1/|x|)$ as $|x|\to+\infty$, that, for every~$\epsilon>0$, there is $R_\epsilon>0$ such that
$$R_{e,\lambda}(H_{e,\lambda}\cap\overline{\Omega_x})\subset\Omega_x$$
for all $e\in\mathbb{S}^1$, $\lambda>\epsilon$ and $|x|>R_\epsilon$.\par
Lastly, consider two points $x\neq y\in\Omega_{0,\infty}$ such that $|x|=|y|$. Let $e\in\mathbb{S}^1$ be defined as in~\eqref{defe}. Consider an arbitrary real number $\epsilon$ such that $0<\epsilon<|x|=|y|$. As in the proofs of Theorems~\ref{th2} and~\ref{th3}, there are two points $x_\epsilon\in\R^2\!\setminus\!\overline{\Omega_X}$ and $x'_\epsilon\in\Omega_X$ such that~$\min_\R|\xi_{x_\epsilon}|>|x|=|y|>\epsilon$,
\beq\label{reflection2}
R_{e,\lambda}(H_{e,\lambda}\cap\overline{\Omega_{x_\epsilon}})\subset\Omega_{x_\epsilon}\ \hbox{ for all }\lambda>\epsilon,
\eeq
and $\max_\R|\xi_{x'_\epsilon}|<\epsilon<|x|=|y|$. The streamlines $\Xi=\Xi_{x_\epsilon}$ and $\Xi'=\Xi_{x'_\epsilon}$ surround the origin, and $u$ is equal to $c_1=u(x_\epsilon)<0$ along $\Xi$ and to $c_2=u(x'_\epsilon)>0$ along $\Xi'$. Furthermore, $u(x_\epsilon)<u<u(x'_\epsilon)$ in the domain
$$\omega=\Omega_{x_\epsilon}\!\setminus\!\overline{\Omega_{x'_\epsilon}}$$
located between $\Xi_{x_\epsilon}$ and $\Xi_{x'_\epsilon}$. Denote $R'=\min_{z\in\Xi'}|z|=\min_\R|\xi_{x'_\epsilon}|\in(0,\epsilon)$, $R=\max_{z\in\Xi}|z|=\max_\R|\xi_{x_\epsilon}|>|x|=|y|>\epsilon>R'$, and $\overline{\lambda}=\max_{z\in\Xi}z\cdot e>\min_\R|\xi_{x_\epsilon}|>|x|=|y|>\epsilon>0$. The~$C^3(\overline{\omega})$ function $\varphi=u$ satisfies~\eqref{eqvarphi} with $[R',R]\times[c_1,c_2]\ni(r,s)\mapsto F(r,s)=f(s)$ satisfying the assumptions of Proposition~\ref{promoving} since $f$ is of class $C^1((-\infty,c_{in}))$. Together with~\eqref{reflection2} and the fact that $H_{e,\lambda}\cap\Xi'=\emptyset$ for all $\lambda>\epsilon$ (since $\Xi'\subset B_\epsilon$), the assumptions~\eqref{hypelambda}-\eqref{hypelambda2} are satisfied. All assumptions of Proposition~\ref{promoving} are therefore fulfilled. Proposition~\ref{promoving} applied with $\lambda=\epsilon$ then implies that $u\le u_{e,\epsilon}$ in $\overline{\omega_{e,\epsilon}}$ with
$$\omega_{e,\epsilon}=\big(H_{e,\epsilon}\cap(\Omega_{x_\epsilon}\!\setminus\!\overline{\Omega_{x'_\epsilon}})\big)\setminus R_{e,\epsilon}(\overline{\Omega_{x'_\epsilon}}).$$
As in the proof of Theorem~\ref{th3}, one has $y\in\omega_{e,\epsilon}$ for all $\epsilon>0$ small enough, hence
$$u(y)\le u_{e,\epsilon}(y)=u(y_{e,\epsilon})=u(y-2(y\cdot e-\epsilon)e)$$
for all $\epsilon>0$ small enough. By passing to the limit as $\epsilon\displaystyle\mathop{\to}^>0$ and using the definition of $e$ and the assumption $|x|=|y|$, one infers that $u(y)\le u(y-2(y\cdot e)e)=u(x)$. Since the last inequality holds for any $x\neq y\in\Omega_{0,\infty}$ such that $|x|=|y|$, the $C^3(\Omega_{0,\infty})$ function $u$ is radially symmetric in $\Omega_{0,\infty}$. Together with the fact that $|\nabla u|=|v|>0$ in $\Omega_{0,\infty}$ and $c_{in}>c_{out}=-\infty$, there is then a $C^3((0,+\infty))$ function $U$ such that $U'<0$ in $(0,+\infty)$ and $u(x)=U(|x|)$ for all~$x\in\Omega_{0,\infty}$. This means that $v(x)=V(|x|)\,e_\theta(x)$ for all $x\in\Omega_{0,\infty}$ with $V=U'\in C^2((0,+\infty))$ and $V<0$ in $(0,+\infty)$. The proof of Theorem~\ref{th4} is thereby complete.\hfill$\Box$


\SE{Proof of the Serrin-type Theorems~\ref{th5} and~\ref{th5bis}}\label{sec40}

We start in Section~\ref{sec41} with the proof of Theorem~\ref{th5bis} dealing with the case of doubly connected bounded domains, since the proof follows easily from the arguments used in the proof of Theorems~\ref{th1}-\ref{th1bis} and on some known results of Reichel~\cite{r1} and Sirakov~\cite{si} on elliptic overdetermined boundary value problems. Section~\ref{sec42} is then devoted to the proof of Theorem~\ref{th5}.


\subsection{Proof of Theorem~\ref{th5bis}}\label{sec41}

Let $\omega_1$, $\omega_2$, $\Omega=\omega_2\!\setminus\!\overline{\omega_1}$ and $v$ be as in Theorem~\ref{th5bis}. Up to translation, one can assume without loss of generality that $0\in\omega_1$, and the assumptions of Theorem~\ref{th5bis} then fall within the framework of Section~\ref{sec2}, with $\mathcal{C}_1=\partial\omega_1$, $\mathcal{C}_2=\partial\omega_2$, $\Omega=\omega_2\!\setminus\!\overline{\omega_1}$, $D=\overline{\Omega}=\overline{\omega_2}\!\setminus\!\omega_1$, and the conditions~\eqref{hypv1}-\eqref{hypv2},~\eqref{C2}, and~\eqref{stagnation} are fulfilled by assumption. Therefore, the flow has a~$C^3(\overline{\Omega})$ stream function~$u$ and, by Lemmas~\ref{lem21} and~\ref{lem24}-(i), there are two real numbers~$c_{in}\neq c_{out}$ (without loss of generality, one can assume that $c_{in}>c_{out}$, even if it means changing~$v$ into~$-v$ and~$u$ into~$-u$) and a~$C^1([c_{out},c_{in}])$ function $f$ such that
\beq\label{eqannulus5}\left\{\baa{l}
\Delta u+f(u)=0\hbox{ in }\overline{\Omega},\vspace{3pt}\\
c_{out}<u<c_{in}\hbox{ in }\Omega,\vspace{3pt}\\
u=c_{in}\hbox{ on }\partial\omega_1,\ \ u=c_{out}\hbox{ on }\partial\omega_2.\eaa\right.
\eeq
Since $|\frac{\partial u}{\partial n}|=|\nabla u|=|v|>0$ along $\partial\Omega$, where $n$ denotes the outward unit normal on $\partial\Omega$, and since $|v|$ is constant along $\partial\omega_1$ and along~$\partial\omega_2$, it follows that $\frac{\partial u}{\partial n}$ is constant too along~$\partial\omega_1$ and along $\partial\omega_2$ (and $\frac{\partial u}{\partial n}>0$ on $\partial\omega_1$ and $\frac{\partial u}{\partial n}<0$ on $\partial\omega_2$). One concludes from~\cite{r1,si} (see also~\cite{al,wgs}) that, up to translation, $\Omega=\Omega_{a,b}$ for some $0<a<b<\infty$ and $u$ is radially symmetric and decreasing with respect to $|x|$ in $\overline{\Omega}=\overline{\Omega_{a,b}}$. The assumptions and the conclusion of Theorem~\ref{th1} are then satisfied and the proof of Theorem~\ref{th5bis} is thereby complete.\hfill$\Box$


\subsection{Proof of Theorem~\ref{th5}}\label{sec42}

Let $\omega$ be a $C^2$ non-empty simply connected bounded domain of $\R^2$ (we here call this domain~$\omega$ instead of $\Omega$ to differentiate it from the notations of Section~\ref{sec2}, which are used below). Let $v\in C^2(\overline{\omega})$ satisfy the Euler equations~\eqref{1} in $\overline{\omega}$. We assume that $v\cdot n=0$ and $|v|$ is constant on~$\partial\omega$, where $n$ denotes the outward unit normal on $\partial\omega$, and that $v$ has a unique stagnation point in~$\overline{\omega}$. Since $\omega$ is simply connected and $v$ is divergence free, there is a $C^3(\overline{\omega})$ stream function $u$ satisfying~\eqref{defu}. Furthermore, $u$ is constant along $\partial\omega$ since $v\cdot n=0$ on $\partial\omega$. Up to normalization, one can assume without loss of generality that
\beq\label{norm6}
u=0\ \hbox{ on }\partial\omega.
\eeq
By assumption, the function $u$ has a unique critical point in $\overline{\omega}$, and $|\nabla u|=|v|$ is constant along $\partial\omega$. Then $|\frac{\partial u}{\partial n}|=|\nabla u|=|v|>0$ on $\partial\omega$. Up to changing $v$ into $-v$ and $u$ into $-u$, one can assume without loss of generality that
\beq\label{overdet}
\frac{\partial u}{\partial n}=\gamma<0\hbox{ on }\partial\omega
\eeq
for some negative constant $\gamma$. Hence, $u$ has a unique maximum point in $\overline{\omega}$ (which is actually in $\omega$) and this point is the unique critical point of $u$ in $\overline{\omega}$. Up to translation, one can assume without loss of generality that this critical point is the origin $0$. One also infers from the uniqueness of the critical point of $u$ that
\beq\label{defL2}
0<u<u(0)\hbox{ for all $x\in\omega\!\setminus\!\{0\}$}.
\eeq
This situation then falls within the framework of Section~\ref{sec2}, with $\mathcal{C}_2=\partial\omega$, $\omega_2=\omega$, $\Omega=\omega\!\setminus\!\{0\}$, and $D=\overline{\omega}\!\setminus\!\{0\}$, and the conditions~\eqref{hypv1}-\eqref{hypv2},~\eqref{C2}, and~\eqref{stagnation} are fulfilled by assumption (in particular, the integral condition~\eqref{hypv2} in the present punctured case is satisfied since here~$v$ is continuous at $0$ and $|v(0)|=0$). Therefore, from Lemma~\ref{lem24}-(i) and the previous notations, there exists a $C^1([0,u(0)))$ function $f$ such that $\Delta u+f(u)=0$ in $\overline{\omega}\!\setminus\!\{0\}$. Furthermore, by setting $f(u(0))=-\Delta u(0)$, it follows from~\eqref{norm6},~\eqref{defL2} and the continuity of $u$ and $\Delta u$ in~$\overline{\omega}$, that~$f$ is continuous at $u(0)$, and then in the whole interval $[0,u(0)]$. Together with the continuity of $u$ and $\Delta u$ at $0$, the equation
$$\Delta u+f(u)=0$$
holds in $\overline{\omega}$. Our goal is to show that $\omega$ is then a ball centered at the origin and that $u$ is radially symmetric and decreasing with respect to $|x|$ in $\overline{\omega}$.\par
Remembering that $u$ satisfies~\eqref{norm6}-\eqref{defL2}, it would then follow from~\cite{s} that $\omega=B_R$ for some $R>0$ and $u$ is radially symmetric and decreasing with respect to $|x|$ in $\overline{\omega}$, if the function~$f$ were known to be Lipschitz continuous in $[0,u(0)]$. However, $f'$ is not bounded in a neighborhood of $u(0)$ in general (see the comments after Theorem~\ref{th5} in Section~\ref{sec12}). We will nevertheless still be able to show the desired symmetry of $\omega$ and of $u$ by taking off from $\omega$ small neighborhoods of $0$ and applying Serrin's strategy and the method of moving planes in punctured domains. The images by $u$ of the closure of these punctured domains are intervals of the type $[0,L]$, with $0<L<u(0)$, and thus $f$ is Lipschitz continuous in $[0,L]$.\par
More precisely, let first $\rho>0$ be such that
$$\overline{B_\rho}\subset\omega$$
and let $e$ be any unit vector. Let $\eta$ be any real number in $(0,\rho)$, and denote
$$\overline{\lambda}_e=\max_{x\in\partial\omega}x\cdot e>\rho>\eta.$$
Using the same notations $T_{e,\lambda}$, $H_{e,\lambda}$ and $R_{e,\lambda}$ as in~\eqref{Telambda}-\eqref{Relambda}, it follows from~\cite{af} that there is $\widetilde{\lambda}\in(\rho,\overline{\lambda}_e)$ such that
\beq\label{deflambda0}
R_{e,\lambda}(H_{e,\lambda}\cap\overline{\omega})\subset\omega\ \hbox{ for all }\lambda\in(\widetilde{\lambda},\overline{\lambda}_e).
\eeq
Since $\max_\R|\xi_x|\to0$ as $|x|\displaystyle\mathop{\to}^>0$ by Lemma~\ref{lem23}-(i) (or here, more simply, because of~\eqref{norm6},~\eqref{defL2} and the continuity of $u$ in $\overline{\omega}$), there is $x_\eta\in\omega\!\setminus\!\{0\}$ such that $\Xi_{x_\eta}\subset B_\eta$. Let then $\Omega'$ be the bounded connected component of $\R^2\!\setminus\!\Xi_{x_\eta}$ (notice that $\overline{\Omega'}\subset B_\eta\subset\omega$) and let
$$\widetilde{\omega}=\omega\!\setminus\!\overline{\Omega'}$$
be the doubly connected bounded domain located between $\Xi_{x_\eta}$ and $\partial\omega$. Notice that $\partial\widetilde{\omega}=\Xi_{x_\eta}\cup\partial\omega$, that $0\not\in\overline{\widetilde{\omega}}$ and that
\beq\label{inequ}
0<u<u(x_\eta)\ \hbox{ in }\widetilde{\omega}
\eeq
since $u$ has no critical point in $\widetilde{\omega}$.\par
From~\eqref{deflambda0}, two cases can occur: either
\begin{itemize}
\item (case~a)
\beq\label{casea}
R_{e,\lambda}(H_{e,\lambda}\cap\overline{\omega})\subset\omega\ \hbox{ for all }\lambda\in[\eta,\overline{\lambda}_e),
\eeq
\item or (case~b) there is $\lambda^*\in[\eta,\widetilde{\lambda}]$ such that $R_{e,\lambda}(H_{e,\lambda}\cap\overline{\omega})\subset\omega$ for all $\lambda\in(\lambda^*,\overline{\lambda}_e)$, and either
\begin{itemize}
\item (internal tangency) there is a point $x^*\in H_{e,\lambda^*}\cap\partial\omega$ such that $x^*_{e,\lambda^*}=R_{e,\lambda^*}(x^*)\in\partial\omega$,
\item or (orthogonality) $T_{e,\lambda^*}$ meets $\partial\omega$ orthogonally, at some point $p^*$.
\end{itemize}
\end{itemize}
We will prove that only case~a occurs.\par
Assume by way of contradiction that case~b occurs. Denote
$$\left\{\baa{l}
\displaystyle\Xi=\partial\Omega,\ \Xi'=\Xi_{x_\eta},\ R'=\min_{x\in\Xi'}|x|\in(0,\rho),\ R=\max_{x\in\Xi}|x|>\rho>R',\ \epsilon=\lambda^*\in[\eta,\widetilde{\lambda}]\subset[0,\overline{\lambda}_e),\vspace{3pt}\\
c_1=0=u_{|\partial\omega},\ c_2=u(x_\eta)=u_{|\Xi_{x_\eta}}\in(0,u(0)).\eaa\right.$$
The $C^3(\overline{\widetilde{\omega}})$ function $\varphi=u$ satisfies $c_1<\varphi<c_2$ in $\widetilde{\omega}$ by~\eqref{inequ} and $\Delta\varphi+F(\varphi)=0$ in $\overline{\widetilde{\omega}}$, where $F(s)=f(s)$ for $s\in[c_1,c_2]\subset[0,u(0))$. The function $F$ is therefore $C^1$ in $[c_1,c_2]$. The condition~\eqref{hypelambda} holds by definition of $\epsilon$, $\lambda^*$ and $\overline{\lambda}_e$, and the condition~\eqref{hypelambda2} is automatically fulfilled since $\Xi'\subset B_\eta$ and $\epsilon=\lambda^* \ge\eta$. Therefore, all assumptions of Proposition~\ref{promoving} are satisfied and it follows from the conclusion applied with $\lambda=\epsilon=\lambda^*$ that
$$u\le u_{e,\lambda^*}\ \hbox{ in }\overline{\omega_{e,\lambda^*}},$$
with
$$\omega_{e,\lambda^*}=(H_{e,\lambda^*}\cap\widetilde{\omega})\setminus R_{e,\lambda^*}(\overline{\Omega'}).$$
Denote
$$w=u_{e,\lambda^*}-u.$$
Since $F=f$ is of class $C^1$ in $[c_1,c_2]\subset[0,u(0))$, the nonnegative $C^3(\overline{\omega_{e,\lambda^*}})$ function $w$ satisfies an equation of the type $\Delta w+cw=0$ in $\overline{\omega_{e,\lambda^*}}$, for some function $c\in C(\overline{\omega_{e,\lambda^*}})$. Thus, the strong maximum principle implies that, for each connected component $\omega'$ of $\omega_{e,\lambda^*}$, either $w>0$ in~$\omega'$, or $w\equiv 0$ in $\overline{\omega'}$. We shall now consider separately the internal tangency case and the orthogonality case.\par
Consider first the case of internal tangency: there is a point $x^*\in H_{e,\lambda^*}\cap\partial\omega$ such that $x^*_{e,\lambda^*}=R_{e,\lambda^*}(x^*)\in\partial\omega$, hence
$$u(x^*)=u(x^*_{e,\lambda^*})=0\ \hbox{ and }\ w(x^*)=0.$$
Since $\overline{\Omega'}\cap\partial\omega=\emptyset$, one has $x^*\not\in\overline{\Omega'}\cup R_{e,\lambda^*}(\overline{\Omega'})$. There is a connected component $\omega^*$ of $\omega_{e,\lambda^*}$ such that $x^*\in\partial\omega^*$, and $B(x^*,r)\cap\omega=B(x^*,r)\cap\omega^*$ for all $r>0$ small enough (in particular, the interior sphere condition in $\omega^*$ is satisfied at the point $x^*\in\partial\omega^*$). Let $n(\zeta)$ be the generic notation for the outward normal to $\omega$ at a point $\zeta\in\partial\omega$. Owing to the definitions of~$\lambda^*$ and~$x^*$, one has $R_{e,\lambda^*}(n(x^*))=n(x^*_{e,\lambda^*})$, while $\nabla u_{e,\lambda^*}(x^*)=R_{e,\lambda^*}(\nabla u(x^*_{e,\lambda^*}))$ owing to the definition of~$u_{e,\lambda^*}$. Hence,
$$\baa{rcl}
\nabla w(x^*)\cdot n(x^*) & = & \nabla u_{e,\lambda^*}(x^*)\cdot n(x^*)-\nabla u(x^*)\cdot n(x^*)\vspace{3pt}\\
& = & R_{e,\lambda^*}(\nabla u(x^*_{e,\lambda^*}))\cdot R_{e,\lambda^*}(n(x^*_{e,\lambda^*}))-\nabla u(x^*)\cdot n(x^*)\vspace{3pt}\\
& = & \nabla u(x^*_{e,\lambda^*})\cdot n(x^*_{e,\lambda^*})-\nabla u(x^*)\cdot n(x^*)=0\eaa$$
since $\nabla u\cdot n$ is equal to the constant $\gamma$ on $\partial\omega$ by~\eqref{overdet}. It then follows from Hopf lemma applied to the function $w$ at the point $x^*$, together with the strong maximum principle, that
$$w\equiv 0\hbox{ in }\overline{\omega^*},\hbox{ that is, }u\equiv u_{e,\lambda^*}\hbox{ in }\overline{\omega^*}.$$
On the other hand, as for formula~\eqref{partials} in the proof of Proposition~\ref{promoving}, one has
$$\partial\omega^*\subset\partial\omega_{e,\lambda^*}\subset\underbrace{\big((T_{e,\lambda^*}\!\cap\overline{\widetilde{\omega}})\!\setminus\!R_{e,\lambda^*}(\Omega')\big)}_{=:\partial_1\omega_{e,\lambda^*}}\,\cup\,\underbrace{\big((H_{e,\lambda^*}\!\cap\partial\omega)\!\setminus\!R_{e,\lambda^*}(\Omega')\big)}_{=:\partial_2\omega_{e,\lambda^*}}\,\cup\,\underbrace{\big(H_{e,\lambda^*}\!\cap\widetilde{\omega}\cap R_{e,\lambda^*}(\Xi_{x_\eta})\big)}_{=:\partial_3\omega_{e,\lambda^*}}.$$
Since $u_{e,\lambda^*}=u(x_\eta)$ on $R_{e,\lambda^*}(\Xi_{x_\eta})$ and $u<u(x_\eta)$ in $\widetilde{\omega}$, one has $w=u_{e,\lambda^*}-u>0$ on $\partial_3\omega_{e,\lambda^*}$, hence $\partial_3\omega_{e,\lambda^*}\cap\partial\omega^*=\emptyset$ and
$$\partial\omega^*\ \subset\ \big((T_{e,\lambda^*}\!\cap\overline{\widetilde{\omega}})\!\setminus\!R_{e,\lambda^*}(\Omega')\big)\cup\big((H_{e,\lambda^*}\!\cap\partial\omega)\!\setminus\!R_{e,\lambda^*}(\Omega')\big)\ \subset\ (T_{e,\lambda^*}\!\cap\overline{\omega})\cup(H_{e,\lambda^*}\!\cap\partial\omega).$$
Therefore, $\omega^*$ is a connected component of $H_{e,\lambda^*}\cap\omega$. Since $w\equiv 0$ in $\overline{\omega^*}$, the arguments of Reichel~\cite{r2} (see also~\cite{ab,si}) imply that $\overline{\omega}=\overline{\omega^*\cup R_{e,\lambda^*}(\omega^*)}$. Hence $\omega$ symmetric with respect to the line $T_{e,\lambda^*}$ and, moreover, $u$ is itself symmetric with respect to $T_{e,\lambda^*}$, which is impossible since $0\not\in T_{e,\lambda^*}$ and $0$ is the only maximum point of $u$. As a consequence, the case of internal tangency is ruled out.\par
Consider now the case of orthogonality, that is, $T_{e,\lambda^*}$ meets $\partial\omega$ orthogonally, at some point~$p^*$. By definition of $u_{e,\lambda^*}$, one has $u(p^*)=u_{e,\lambda^*}(p^*)$, thus $w(p^*)=0$. Notice also, as in the case of internal tangency, that $p^*\not\in\overline{\Omega'}\cup R_{e,\lambda^*}(\overline{\Omega'})$. There is a connected component $\omega^*$ of~$\omega_{e,\lambda^*}$ such that $p^*\in\partial\omega^*$, and $B(p^*,r)\cap\omega\cap H_{e,\lambda^*}=B(p^*,r)\cap\omega^*$ for all $r>0$ small enough. Since $u$ and $\frac{\partial u}{\partial n}$ are constant on $\partial\omega$ and since $T_{e,\lambda^*}$ meets $\partial\omega$ orthogonally at $p^*$, it follows as in~\cite{r2} that all first and second order derivatives of $w$ vanish at $p^*$. Serrin's corner lemma~\cite{s} and the strong maximum principle then yield $w\equiv 0$ in $\overline{\omega^*}$. One is then led to a contradiction as in the previous paragraph.\par
As a consequence, only case~a occurs. Thus,~\eqref{casea} holds. By arguing as in the beginning of the study of case~b and applying Proposition~\ref{promoving} with this time $\epsilon=\eta$, one infers that
\beq\label{casea2}
u\le u_{e,\lambda}\ \hbox{ in }\overline{\omega_{e,\lambda}}\ \hbox{ for all }\lambda\in[\eta,\overline{\lambda}_e),
\eeq
with $\omega_{e,\lambda}=(H_{e,\lambda}\cap\widetilde{\omega})\setminus R_{e,\lambda}(\overline{\Omega'})$. Since~\eqref{casea} and~\eqref{casea2} hold for every direction $e\in\mathbb{S}^1$ and for every $\eta\in(0,\rho)$, one finally concludes that
$$\omega=B_R$$
for some $R>0$ and, as in the proof of Theorem~\ref{th3}, that $u$ is radially symmetric in $\overline{\omega}=\overline{B_R}$. Since $0$ is the unique critical point of the $C^3(\overline{B_R})$ function $u$ and since $u=0$ on $\partial B_R$ with~$u>0$ in $B_R$, there is then a $C^3([0,R])$ function $U:[0,R]\to\R$ such that $u(x)=U(|x|)$ in $\overline{B_R}$, with~$U'(0)=0$ and $U'<0$ in $(0,R]$. Therefore,
$$v(x)=\nabla ^\perp u(x)=U'(|x|)e_\theta(x)\hbox{ for all }x\in\overline{B_R}\!\setminus\!\{0\}$$
and the $C^2([0,R])$ function $V=U'$ satisfies the desired conclusion. The proof of Theorem~\ref{th5} is thereby complete.\hfill$\Box$


\SE{Proof of Proposition~\ref{promoving}}\label{sec50}

It is based on the method of moving planes developed in~\cite{a,bn,gnn,s}, though it has to be adapted here to our geometrical configuration. The idea is to compare the function $\varphi$ to its reflection $\varphi_{e,\lambda}$ in $\overline{\omega_{e,\lambda}}$ by moving the lines $T_{e,\lambda}$ and decreasing $\lambda$ from the value $\overline{\lambda}$ to the value~$\epsilon$. We recall that
$$\omega_{e,\lambda}=(H_{e,\lambda}\cap\omega)\setminus R_{e,\lambda}(\overline{\Omega'}).$$
Notice in particular that $R'\le|x|\le R$ for all $\lambda\in[\epsilon,\overline{\lambda})$ and $x\in\overline{\omega_{e,\lambda}}$, since $\omega_{e,\lambda}\subset\omega=\Omega\!\setminus\!\overline{\Omega'}$.\par
Consider first any $\lambda\in(\epsilon,\overline{\lambda})$. For each $x\in\omega_{e,\lambda}$, there holds
$$x_{e,\lambda}=R_{e,\lambda}(x)\in R_{e,\lambda}(H_{e,\lambda}\cap\omega)\subset R_{e,\lambda}(H_{e,\lambda}\cap\overline{\Omega})\subset\Omega$$
by~\eqref{hypelambda}, and $x_{e,\lambda}\not\in\overline{\Omega'}$, hence, $x_{e,\lambda}\in\omega$. Thus
$$R_{e,\lambda}(\overline{\omega_{e,\lambda}})\subset\overline{\omega}$$
and the function $\varphi_{e,\lambda}$ given in~\eqref{varphielambda} is well defined and of class $C^2$ in $\overline{\omega_{e,\lambda}}$. Furthermore, $\Delta\varphi_{e,\lambda}+F(|x_{e,\lambda}|,\varphi_{e,\lambda})=0$ in $\overline{\omega_{e,\lambda}}$. Since $|x|\ge|x_{e,\lambda}|$ for all $x\in\overline{\omega_{e,\lambda}}$ (remember that $\lambda>\epsilon\ge0$) and since $F$ is nonincreasing with respect to its first variable, it follows that
$$\Delta\varphi_{e,\lambda}+F(|x|,\varphi_{e,\lambda})\le0\ \hbox{ in }\overline{\omega_{e,\lambda}}.$$
Let
$$\Phi_{e,\lambda}=\varphi_{e,\lambda}-\varphi,$$
which is well defined and of class $C^2$ in $\overline{\omega_{e,\lambda}}$. There holds
\beq\label{eqPhi}
\Delta \Phi_{e,\lambda}+c_{e,\lambda}\Phi_{e,\lambda}\le0\ \hbox{ in }\overline{\omega_{e,\lambda}},
\eeq
where, say,
$$c_{e,\lambda}(x)=\left\{\baa{ll}
\displaystyle\frac{F(|x|,\varphi_{e,\lambda}(x))-F(|x|,\varphi(x))}{\varphi_{e,\lambda}(x)-\varphi(x)} & \hbox{if }\varphi_{e,\lambda}(x)\neq\varphi(x),\vspace{3pt}\\
0 & \hbox{if }\varphi_{e,\lambda}(x)=\varphi(x).\eaa\right.$$
Since the function $F$ is assumed to be Lipschitz continuous with respect to its second variable, uniformly with respect to the first one, the function $c_{e,\lambda}$ is in $L^\infty(\omega_{e,\lambda})$ and, moreover, there is a constant $M\ge0$ such that
\beq\label{celambda}
|c_{e,\lambda}(x)|\le M\ \hbox{ for all }\lambda\in(\epsilon,\overline{\lambda})\hbox{ and for all }x\in\overline{\omega_{e,\lambda}}.
\eeq\par 
Consider again any $\lambda\in(\epsilon,\overline{\lambda})$ and let us decompose the boundary of $\omega_{e,\lambda}$ into three parts. More precisely, since
$$\partial(A\cap B\cap C)\ \subset\ \big(\partial A\cap\overline{B}\cap\overline{C}\big)\,\cup\,\big(A\cap\partial B\cap\overline{C}\big)\,\cup\,\big(A\cap B\cap\partial C\big)$$
for any three sets $A$, $B$ and $C$, since $\partial\omega=\Xi\cup\Xi'$ and since $(H_{e,\lambda}\cap\Xi')\setminus R_{e,\lambda}(\Omega')=\emptyset$ by assumption~\eqref{hypelambda2}, one has (with $A=H_{e,\lambda}$, $B=\omega$ and $C=\R^2\setminus R_{e,\lambda}(\overline{\Omega'})$)
\beq\label{partials}
\partial\omega_{e,\lambda}\subset\underbrace{\big((T_{e,\lambda}\cap\overline{\omega})\!\setminus\!R_{e,\lambda}(\Omega')\big)}_{=:\partial_1\omega_{e,\lambda}}\ \cup\ \underbrace{\big((H_{e,\lambda}\cap\Xi)\!\setminus\!R_{e,\lambda}(\Omega')\big)}_{=:\partial_2\omega_{e,\lambda}}\ \cup\ \underbrace{\big(H_{e,\lambda}\cap\omega\cap R_{e,\lambda}(\Xi')\big)}_{=:\partial_3\omega_{e,\lambda}},
\eeq
see Fig.~4. Notice that, since $T_{e,\lambda}\cap\Xi=T_{e,\lambda}\cap\partial\Omega$ is not empty (because $\lambda\in(\epsilon,\overline{\lambda})\subset[0,\overline{\lambda})$), both sets $\partial_1\omega_{e,\lambda}$ and $\partial_2\omega_{e,\lambda}$ are not empty (however, $\partial_3\omega_{e,\lambda}$ may be empty). Furthermore, even if~$\omega_{e,\lambda}$ may not be connected (as in Fig.~4), the boundary of each connected component of~$\omega_{e,\lambda}$ intersects $\partial_2\omega_{e,\lambda}\cup\partial_3\omega_{e,\lambda}$.\par
\begin{figure}
\centering\includegraphics[scale=0.7]{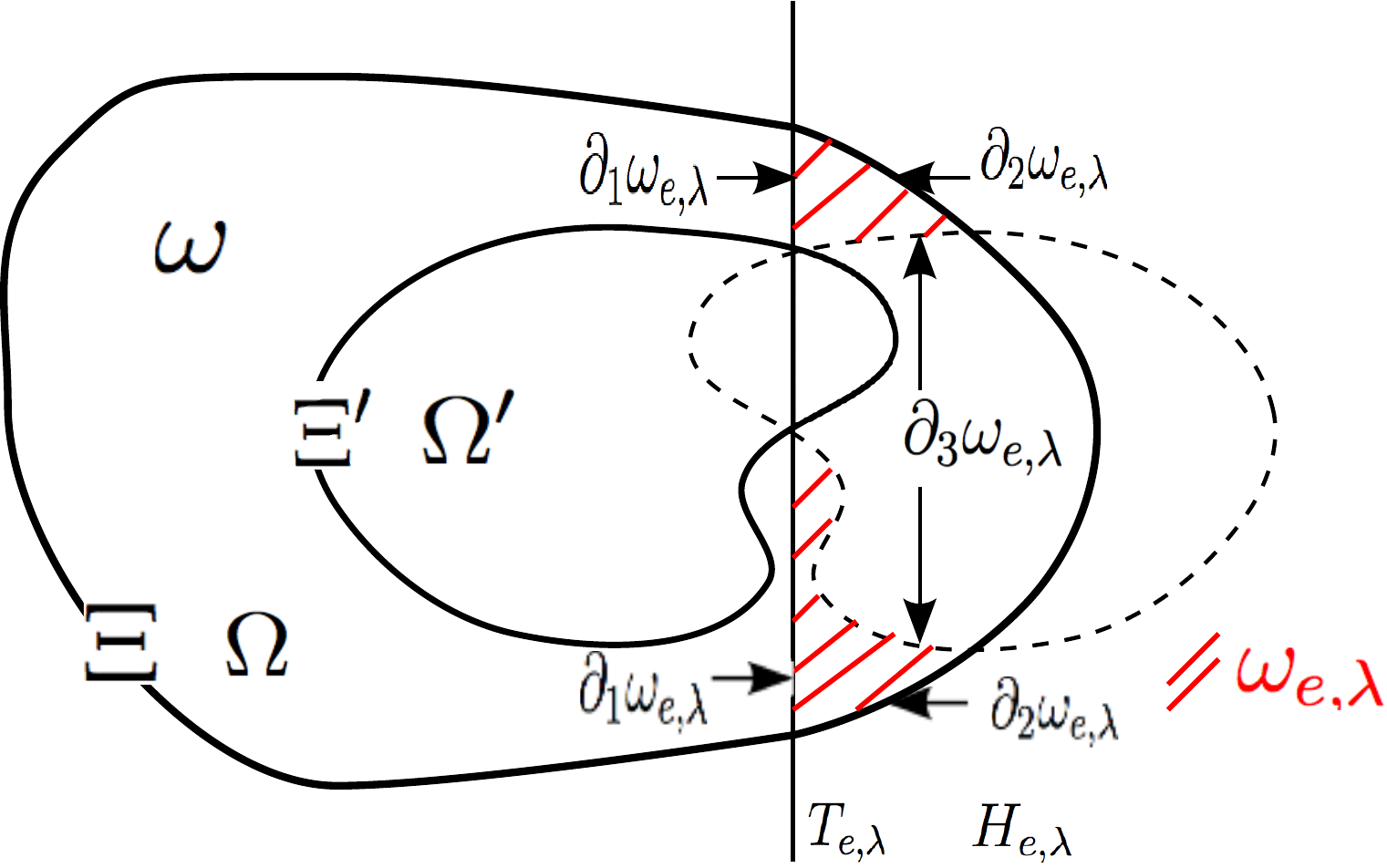}
\caption{The three parts $\partial_1\omega_{e,\lambda}$, $\partial_2\omega_{e,\lambda}$ and $\partial_3\omega_{e,\lambda}$ of the boundary of the set $\omega_{e,\lambda}$ (dashed red)}
\end{figure}
Let us now study the sign of $\Phi_{e,\lambda}$ on $\partial\omega_{e,\lambda}$, for any $\lambda\in(\epsilon,\overline{\lambda})$. Firstly, on $\partial_1\omega_{e,\lambda}\,(\subset T_{e,\lambda})$, one has $\varphi_{e,\lambda}=\varphi$, hence $\Phi_{e,\lambda}=0$. Secondly, for each $x\in\partial_2\omega_{e,\lambda}$, one has
$$x_{e,\lambda}\in R_{e,\lambda}(H_{e,\lambda}\cap\Xi)\subset R_{e,\lambda}(H_{e,\lambda}\cap\overline{\Omega})\subset\Omega$$
by~\eqref{hypelambda}, hence $x_{e,\lambda}\in\omega\cup\Xi'$ and $\varphi_{e,\lambda}(x)=\varphi(x_{e,\lambda})>c_1$ by~\eqref{eqvarphi}, while $x\in\Xi$ and $\varphi(x)=c_1$. Thus, $\Phi_{e,\lambda}(x)=\varphi_{e,\lambda}(x)-\varphi(x)>0$ for each $x\in\partial_2\omega_{e,\lambda}$. Thirdly, for each $x\in\partial_3\omega_{e,\lambda}$, one has $x_{e,\lambda}\in\Xi'$ and $\varphi_{e,\lambda}(x)=\varphi(x_{e,\lambda})=c_2$, while $x\in\omega$ and $\varphi(x)<c_2$, by~\eqref{eqvarphi}. Thus, $\Phi_{e,\lambda}(x)=\varphi_{e,\lambda}(x)-\varphi(x)>0$ for each $x\in\partial_3\omega_{e,\lambda}$. As a consequence, $\Phi_{e,\lambda}\ge0$ on $\partial\omega_{e,\lambda}$ and even $\Phi_{e,\lambda}>0$ on $\partial_2\omega_{e,\lambda}\cup\partial_3\omega_{e,\lambda}\ (\neq\emptyset)$, hence
\beq\label{boundaries}
\hbox{$\Phi_{e,\lambda}\ge\not\equiv0$ on the boundary of each connected component of $\omega_{e,\lambda}$}.
\eeq\par
Let us now consider $\lambda\simeq\overline{\lambda}$ with $\lambda<\overline{\lambda}$. Since the functions $\Phi_{e,\lambda}$ satisfy~\eqref{eqPhi}-\eqref{boundaries}, since the sets $\omega_{e,\lambda}$ are all included in the given bounded domain $\Omega$ and since the Lebesgue measure $|\omega_{e,\lambda}|$ of $\omega_{e,\lambda}$ goes to $0$ as $\lambda\displaystyle\mathop{\to}^<\overline{\lambda}$ owing to the definition of $\overline{\lambda}$, it follows for instance from the maximum principle in sets with bounded diameter and small Lebesgue measure and from the strong maximum principle~\cite{bnv}, that there is $\lambda_0\in(\epsilon,\overline{\lambda})$ such that $\Phi_{e,\lambda}>0$ in~$\omega_{e,\lambda}$ for all $\lambda\in(\lambda_0,\overline{\lambda})$.\par
Let us finally define
$$\lambda_*=\inf\big\{\lambda\in(\epsilon,\overline{\lambda}):\Phi_{e,\lambda'}>0\hbox{ in }\omega_{e,\lambda'}\hbox{ for all }\lambda'\in(\lambda,\overline{\lambda})\big\},$$
and notice that $\epsilon\le\lambda_*\le\lambda_0<\overline{\lambda}$. Our goal is to show that $\lambda_*=\epsilon$. Assume by way of contradiction that $\lambda_*>\epsilon$. Notice that $\Phi_{e,\lambda_*}\ge0$ in $\overline{\omega_{e,\lambda_*}}$ by continuity (indeed, for each $x\in\omega_{e,\lambda_*}$, there holds $x\in\omega_{e,\lambda}$ for $\lambda-\lambda_*>0$ small, hence $\varphi(x)<\varphi_{e,\lambda}(x)$ for $\lambda-\lambda_*>0$ small, and $\varphi(x)\le\varphi_{e,\lambda_*}(x)$ by passing to the limit $\lambda\displaystyle\mathop{\to}^>\lambda_*$ and by continuity of $\varphi$; therefore, $\varphi\le\varphi_{e,\lambda_*}$ in $\overline{\omega_{e,\lambda_*}}$ again by continuity of $\varphi$). On the other hand, $\Phi_{e,\lambda_*}\ge\not\equiv0$ on the boundary of each connected component of $\omega_{e,\lambda_*}$, because $\lambda_*\in(\epsilon,\overline{\lambda})$. Hence, $\Phi_{e,\lambda_*}>0$ in $\omega_{e,\lambda_*}$ from the strong maximum principle. As in the previous paragraph, from~\cite{bnv}, there exists $\delta>0$ such that the weak maximum principle holds in any open set $\omega'\subset\omega$ for the solutions $\Phi\in C^2(\omega')\cap C(\overline{\omega'})$ of $\Delta\Phi+c\Phi\le0$ in $\omega'$ with $\Phi\ge0$ on $\partial\omega'$ and $\|c\|_{L^\infty(\omega')}\le M$, as soon as $|\omega'|\le\delta$. Let then $K$ be a compact subset of $\omega_{e,\lambda_*}$ such that
$$|\omega_{e,\lambda_*}\!\setminus\!K|<\frac{\delta}{2}.$$
Since $\min_K\Phi_{e,\lambda_*}>0$, it follows from the continuity of $\varphi$ in $\overline{\omega}$ that there exists $\underline{\lambda}\in(\epsilon,\lambda_*)$ such that, for all $\lambda\in[\underline{\lambda},\lambda_*]$,
$$\min_K\Phi_{e,\lambda}>0,\ \ \partial(\omega_{e,\lambda}\!\setminus\!K)=\partial\omega_{e,\lambda}\cup\partial K\ \hbox{ and }\ |\omega_{e,\lambda}\!\setminus\!K|<\delta.$$
For any such $\lambda\in[\underline{\lambda},\lambda_*]$, one then has $\Phi_{e,\lambda}\ge\not\equiv0$ on the boundary of each connected component of $\omega_{e,\lambda}\!\setminus\!K$ and one then infers from the choice of $\delta$ and from the strong maximum principle that~$\Phi_{e,\lambda}>0$ in $\omega_{e,\lambda}\!\setminus\!K$, and finally $\Phi_{e,\lambda}>0$ in $\omega_{e,\lambda}$. This last property contradicts the definition of $\lambda_*$.\par
As a conclusion, $\lambda_*=\epsilon$. Therefore, for every $\lambda\in(\epsilon,\overline{\lambda})$, one has $\Phi_{e,\lambda}>0$ in $\omega_{e,\lambda}$, namely~$\varphi<\varphi_{e,\lambda}$ in $\omega_{e,\lambda}$ and $\varphi\le\varphi_{e,\lambda}$ in $\overline{\omega_{e,\lambda}}$ by continuity of $\varphi$. As in the previous paragraph, it also follows by continuity that $\varphi\le\varphi_{e,\epsilon}$ in $\overline{\omega_{e,\epsilon}}$. The proof of Proposition~\ref{promoving} is thereby complete.~\hfill$\Box$\break

\noindent{\bf{Acknowledgements.}} The authors are grateful to the anonymous reviewers for pointing out the references~\cite{f,gpsy} and for very valuable comments which improved the readability and quality of the paper.


\end{document}